\documentstyle[11pt]{article}

\newtheorem{theorem}{$~~~~$ Theorem}[section]
\newtheorem{corollary}{Corollary}

\newtheorem{example}[theorem]{$~~~~$ Example}
\newtheorem{lemma}[theorem]{$~~~~$ Lemma}
\newtheorem{remark}[theorem]{$~~~~$Remark}
\newtheorem{definition}[theorem]{$~~~~$Definition}

\def\nor1{Normed$\{~2^{ \zzz \theta  \, )} ~$,$~\sqrt{~2^{ \zzz \theta  \, )}}~\}$}

\def\xor2{Normed$\{ ~\sqrt{~2^{ \zzz \theta  \, )}}~,~2~ \} $}

\def\zhz{H }
\def\appD{Appendix D }

\def\gggen{$( L^\xi  ,  \Delta_0^\xi   ,  B^\xi  ,  d  ,  G  )~$}
\def\gggcp{$( L^\xi  ,  \Delta_0^\xi   ,  B^\xi  ,  d  ,  G  )$}
\def\peta{\sigma}

\def\zzz{~ \sharp ( ~ }

\def\mheta{\theta^\bullet}
\def\mxi{\xi^\bullet}
\def\xxi{$\, \xi^* \,$}

\def\tftt{~ \frac{1}{2}~ }
\def\sss{ }

\def\Uxp{\Upsilon}

\def\f55{ \normalsize  \baselineskip = 1.8 \normalbaselineskip }

\def\f55{  \baselineskip = 1.1 \normalbaselineskip } 
\def\g55{  \baselineskip = 1.0 \normalbaselineskip } 
\def\s55{ \baselineskip = 1.0 \normalbaselineskip }

\newcommand{\thx}[1]{Theorem \ref{#1}}
\newcommand{\phx}[1]{Theorem \ref{#1}}
\newcommand{\lem}[1]{Lemma \ref{#1}}
\newcommand{\eq}[1]{(\ref{#1})}
\newcommand{\ep}[1]{Equation (\ref{#1})}

\begin{document}

 \title{A Detailed Examination
of Methods for Unifying, Simplifying and Extending 
Several
Results About Self-Justifying Logics}

\def\beq{\begin{equation}}
\def\enq{\end{equation}}

\def\bel{\begin{lemma}}
\def\enl{\end{lemma}}

\def\bec{\begin{corollary}}
\def\enc{\end{corollary}}

\def\bed{\begin{description}}
\def\ennd{\end{description}}
\def\bee{\begin{enumerate}}
\def\ene{\end{enumerate}}

\def\bxbxd{\begin{definition}}
\def\bxbxdd{\begin{definition}}
\def\eedd{\end{definition}}
\def\bxbxdr{\begin{definition} \rm}
\def\bel{\begin{lemma}}
\def\enl{\end{lemma}}
\def\ent{\end{theorem}}

\author{  Dan E.Willard\thanks{This research 
was partially supported
by the NSF Grant CCR  0956495.
Email = dew@cs.albany.edu.}}

\date{State University of New York at Albany}

\maketitle

\setcounter{page}{0}
\thispagestyle{empty}

\normalsize

\baselineskip = 1.5\normalbaselineskip

\normalsize

\baselineskip = 1.0 \normalbaselineskip 
\def\bbint{\large \baselineskip = 1.6 \normalbaselineskip } 
\def\bbint{\large \baselineskip = 1.6 \normalbaselineskip }
\def\bbint{\normalsize \baselineskip = 1.3 \normalbaselineskip }

\def\bbint{\normalsize \baselineskip = 1.27 \normalbaselineskip }

\def\bbint{\large \baselineskip = 2.0 \normalbaselineskip }

\def\bbint{\normalsize \baselineskip = 1.25 \normalbaselineskip }
\def\bbina{\normalsize \baselineskip = 1.24 \normalbaselineskip }

\def\bbint{\large \baselineskip = 2.0 \normalbaselineskip }

\def\bbing{ }
\def\bbins{ }
\def\bbinm{ }

\def\bbint{\normalsize \baselineskip = 1.95 \normalbaselineskip }

\def\bbing{ }
\def\bbins{ }
\def\bbinm{ }

\def\bbint{\large \baselineskip = 2.3 \normalbaselineskip } 
\def\bbing{ }
\def\bbins{ }
\def\bbinm{ }

\def\bbint{\normalsize \baselineskip = 1.7 \normalbaselineskip } 

\def\bbint{\large \baselineskip = 2.3 \normalbaselineskip } 
\def\bbinm{ \baselineskip = 1.18 \normalbaselineskip }

\def\bbint{\large \baselineskip = 2.0 \normalbaselineskip } 
\def\bbing{ }
\def\bbins{ }
\def\bbinm{ }
\def\bbinr{ }

\def\bbint{\normalsize \baselineskip = 1.25 \normalbaselineskip }
\def\bbina{\normalsize \baselineskip = 1.24 \normalbaselineskip }
\def\bbinr{ \baselineskip = 1.3 \normalbaselineskip }
\def\bbing{ \baselineskip = 1.28 \normalbaselineskip }
\def\bbins{ \baselineskip = 1.21 \normalbaselineskip }
\def\bbinm{  }

\bbint

\parskip 5 pt

\noindent

\begin{abstract}
 \baselineskip = 1.5 \normalbaselineskip \large 

This paper will develop a single framework for unifying, simplifying
and extending our prior results about axiom systems that retain a
partial knowledge of their own consistency, via an axiomatic
declaration of self-consistency.  Its perhaps single most surprising
new result will be its exploration of a viable alternative to
conventional reflection principles.
\end{abstract}

\normalsize

\parskip 8pt

\baselineskip = 1.4 \normalbaselineskip 

\setcounter{page}{0}

\bigskip
\bigskip
\bigskip


\bigskip
\bigskip


\bigskip
\bigskip

% {\bf Please Cite this Paper as:}
% {\rm http://arxiv.org/abs/1108.6330}, 
%  appearing in Cornell Archives 

\newpage

\setcounter{page}{1}

\def\ww22{\normalsize \baselineskip = 1.21\normalbaselineskip \parskip 4 pt}
\def\bb22{\normalsize \baselineskip = 1.19\normalbaselineskip \parskip 4 pt}
\def\zz22z{\normalsize \baselineskip = 1.19 \normalbaselineskip \parskip 3 pt}
\def\xx22{\normalsize \baselineskip = 1.17\normalbaselineskip \parskip 4 pt}
\def\vx22s{\normalsize \baselineskip = 1.16 \normalbaselineskip \parskip 3 pt} 
\def\vv22{\normalsize \baselineskip = 1.17 \normalbaselineskip \parskip 3 pt} 
\def\aa22{\normalsize \baselineskip = 1.15 \normalbaselineskip \parskip 3 pt} 
\def\f55{  \baselineskip = 1.1 \normalbaselineskip } 
\def\h55{  \baselineskip = 1.08 \normalbaselineskip } 
\def\g55{  \baselineskip = 1.0 \normalbaselineskip } 
\def\s55{ \baselineskip = 1.0 \normalbaselineskip } 
\def\sm55{ \baselineskip = 0.9 \normalbaselineskip }

\vspace*{- 1.0 em}

\def\waw11{\normalsize \baselineskip = 1.72\normalbaselineskip}
\def\waw11{\normalsize \baselineskip = 1.12\normalbaselineskip}
\def\waw11{\normalsize \baselineskip = 1.85\normalbaselineskip}

\def\waw11{\normalsize \baselineskip = 1.45\normalbaselineskip}

\def\waw11{\normalsize \baselineskip = 1.7\normalbaselineskip}

\def\waw11{\normalsize \baselineskip = 1.12\normalbaselineskip}

\def\f55{  \baselineskip = 1.59 \normalbaselineskip } 
\def\g55{  \baselineskip = 1.50 \normalbaselineskip } 
\def\s55{ \baselineskip = 1.50 \normalbaselineskip } 
\def\sm55{ \baselineskip = 1.5 \normalbaselineskip }

\def\f55{  \baselineskip = 1.59 \normalbaselineskip } 
\def\g55{  \baselineskip = 1.50 \normalbaselineskip } 
\def\s55{ \baselineskip = 1.50 \normalbaselineskip } 
\def\sm55{ \baselineskip = 0.9 \normalbaselineskip }

\def\aa22{\normalsize  \waw11 \parskip 6 pt} 
\def\bb22{\normalsize  \waw11 \parskip 5 pt}
\def\ww22{\normalsize \waw11 \parskip 4 pt}
\def\vv22{\normalsize  \waw11 \parskip 3 pt} 
\def\tt22{\normalsize  \waw11 \parskip 2 pt} 

\def\f55{  \baselineskip = 1.1 \normalbaselineskip } 
\def\g55{  \baselineskip = 1.0 \normalbaselineskip } 
\def\b55{  \baselineskip = 1.0 \normalbaselineskip } 
\def\s55{ \baselineskip = 1.0 \normalbaselineskip } 
\def\sm55{ \baselineskip = 0.9 \normalbaselineskip }

\def\mal{ \bf  }
\def\nal{\mathcal}
\def\cvt{ \baselineskip = 0.98 \normalbaselineskip }
\def\cv9{ \baselineskip = 0.99 \normalbaselineskip }
\def\cvs{ \baselineskip = 1.0 \normalbaselineskip }
\def\cvl{ \baselineskip = 1.0 \normalbaselineskip }
\def\cvh{ \baselineskip = 1.03 \normalbaselineskip }
\def\cvg{ \baselineskip = 1.00 \normalbaselineskip }

\def\cvt{ \baselineskip = 1.6 \normalbaselineskip }
\def\cv9{ \baselineskip = 1.6 \normalbaselineskip }
\def\cvs{ \baselineskip = 1.6 \normalbaselineskip }
\def\cvl{ \baselineskip = 1.6 \normalbaselineskip }
\def\cvh{ \baselineskip = 1.6 \normalbaselineskip }
\def\cvg{ \baselineskip = 1.6 \normalbaselineskip }
\def\cvb{ \baselineskip = 1.6 \normalbaselineskip }
\def\cvnew{ \baselineskip = 1.6 \normalbaselineskip }
\def\cvmew{ \baselineskip = 1.6 \normalbaselineskip }
\def\cvwew{ \baselineskip = 1.6 \normalbaselineskip \parskip 5pt }
\def\cvrew{ \baselineskip = 1.6 \normalbaselineskip \parskip 3pt }

\def\cvt{ \baselineskip = 1.22 \normalbaselineskip }
\def\cv9{ \baselineskip = 1.22 \normalbaselineskip }
\def\cvs{ \baselineskip = 1.22 \normalbaselineskip }
\def\cvl{ \baselineskip = 1.22 \normalbaselineskip }
\def\cvh{ \baselineskip = 1.22 \normalbaselineskip }
\def\cvg{ \baselineskip = 1.22 \normalbaselineskip }
\def\cvb{ \baselineskip = 1.22 \normalbaselineskip }
\def\cvnew{ \baselineskip = 1.4 \normalbaselineskip }
\def\cvmew{ \baselineskip = 1.35 \normalbaselineskip }
\def\cvwew{ \baselineskip = 1.4 \normalbaselineskip \parskip 5pt }
\def\cvrew{ \baselineskip = 1.22 \normalbaselineskip \parskip 3pt }

\def\fend{ 

\medskip -------------------------------------------------------}

\def\f55{  \baselineskip = 1.1 \normalbaselineskip } 
\def\g55{  \baselineskip = 1.0 \normalbaselineskip } 
\def\s55{ \baselineskip = 1.0 \normalbaselineskip } 
\def\sm55{ \baselineskip = 1.0 \normalbaselineskip } 
\def\h55{  \baselineskip = 1.08 \normalbaselineskip } 
\def\b55{  \baselineskip = 1.1 \normalbaselineskip } 

\def\nop{ }
\def\nop{\newpage}

\cvl
\cvnew

\def\nop{\newpage}

\cvnew
\cvl

\parskip 3pt

\def\hgskip{ \medskip }

\def\nop{ }
\def\njp{\newpage}
\def\njp{ }
\def\nop{\newpage}

\cvl
\cvnew

\def\nskip{\bigskip}

\cvl
\cvnew

\section{Introduction}

%% CHANGES REMOVE NEW PAGE IN BIB before willard 2001

%% CHANGES = - 3 footnote -ACKN -APPF -TABLE NUMBERS

%% CHANGES + ww11 Reminder-to-reader Footnote 24 state TABLE

\label{secc1}
\label{B1-lem}
\label{D1-def}

\parskip 2pt

Let
$~\alpha~$ 
denote an axiom system, 
and $~d~$  
denote
 a
deduction method.
The ordered pair
 $~(  \alpha  , d  )$
will 
be called  {\bf Self Justifying} when:
\begin{description}
  \item[  i   ] one of $ \, \alpha \,$'s  theorems
states that the deduction method $ \, d, \, $ applied to the
system $ \, \alpha, \, $ will produce a consistent set of theorems, and
\item[  ii   ]
     the axiom system $ \, \alpha  \, $ is in fact consistent.
\end{description}
For any  $\,(\alpha,d) \,$, 
it is 
easy 
to construct a second
axiom system $ \, \alpha^d \, \supseteq  \,  \alpha  \, $
 that  satisfies
Part-i of 
this definition.
For instance,  $ \, \alpha^d \, $  could
consist of all of $~\alpha \,$'s axioms plus the following added
sentence,  that we call
{\bf SelfRef$(\alpha,d)~$}:
\topsep -3pt
\begin{quote}
$\bullet~~~$ 
There is no proof 
(using 
$d$'s deduction method)
of  $0=1$
from the  {\it union} of
the
 axiom system $\, \alpha \, $
with {\it this}
sentence  ``SelfRef$(\alpha,d) \,$'' (looking at itself).
\end{quote}
Kleene 
\cite{Kl38} 
discussed
how
to
encode
approximate
 analogs of
SelfRef$(\alpha,d)$'s
 self-referential statement.
Each of
Kleene, 
Rogers and Jeroslow 
 \cite{Kl38,Ro67,Je71}
 noted
$\alpha ^d$ 
may,
however,  be inconsistent
(despite SelfRef$(\alpha,d)$'s assertion),
thus causing 
it
to violate Part-ii of   self-justification's
definition. 

\smallskip

This problem arises in
settings
more general than 
 G\"{o}del's
paradigm,
where $\alpha$  was an extension of Peano Arithmetic.
There 
are 
many
settings 
where the Second Incompleteness Theorem does
generalize
\cite{Ad2,AB1,AZ1,BS76,Bu86,BI95,Fe60,Go31,HP91,HB39,Ko6,KT74,Lo55,PD83,PW81,Pu84,Pu85,Pu96,Ro67,Sa11,Sm85,So94,Sv7,Ta0,VV94,Vi92,Vi5,WP87,ww2,wwapal,wwlogos,ww7}.
Each such result formalizes a 
paradigm where 
self-justification is infeasible,
due to a diagonalization issue.
Most
logicians 
have
hesitated to
thus
 employ 
a SelfRef$(\alpha,d)$
 axiom
because
$\alpha+$SelfRef$(\alpha,d)  $
is usually 
 inconsistent     \footnote{ \baselineskip = 1.3 \normalbaselineskip  \label{troub} 
     Typical ordered pairs $(\alpha,d)$
     will have the property that
     the broader axiom system
     $~\alpha^d~=~\alpha \,+ \,$SelfRef$(\alpha,d)$ will
     be inconsistent, even
     when  $~\alpha~$ is consistent. This is because
     a  
     standard
      G\"{o}del-like self-referencing 
     %diagonalization 
     construction
  will 
     %usually
     typically 
     % be able to 
     produce a proof of $0=1$ from
     $~\alpha^d\,$, irregardless of whether or not $~\alpha$ is 
%     formally
     consistent.}.

\smallskip

Our research
explored special
circumstances
\cite{ww1,ww5,ww6,wwapal}
where it is feasible to
construct self-justifying formalisms.
These paradigms involved weakening
the properties a system can prove about
addition and/or 
multiplication
(to avoid the preceding
difficulties).
To be more precise, let
 $Add(x,y,z)$ and    $Mult(x,y,z)$ 
denote 
two 
3-way predicates
 indicating
$x$, $y$ and $z$ satisfy 
$x+y=z$ and
$x*y=z$.
A 
logic
will be said to
{\bf recognize}
successor, 
 addition  and multiplication
as {\bf Total Functions} iff it 
includes
1-3 as axioms.

\vspace*{- 0.8 em}
{\  
\cvl
\beq 
\label{totdefxs}
\forall x ~ \exists z ~~~Add(x,1,z)~~
\enq
\vspace*{- 1.7 em}
\beq 
\label{totdefxa}
\forall x ~\forall y~ \exists z ~~~Add(x,y,z)~~
\enq
\vspace*{- 1.7 em}
\beq 
\label{totdefxm}
\forall x ~\forall y ~\exists z ~~~Mult(x,y,z)~
\enq }

\vspace*{- 0.6 em}

\cvl

We will say
a 
logic
system 
$\alpha$
is 
{\bf Type-M} iff it contains
each of \eq{totdefxs} -- \eq{totdefxm}
as axioms,  
{\bf Type-A} iff it contains only
\eq{totdefxs} and \eq{totdefxa} as axioms,
and {\bf Type-S} iff it contains
only \eq{totdefxs} as an
 axiom. 
A system is called 
{\bf Type-NS} iff it {\it does not} contain
any of these axioms.

Our investigations
 \cite{ww1}--\cite{ww7}
began by observing
some 
Type-A systems can  recognize 
their 
consistency under semantic tableaux deduction,
and 
several
Type-NS systems 
can 
recognize their
 Hilbert consistency.
Many of
these systems were capable of 
proving 
%all 
Peano Arithmetic's
 $\Pi_1$ theorems
in a language
that represents addition and multiplication 
as
the
3-way predicates
of
 Add$(x,y,z)$ and  Mult$(x,y,z)$.

\parskip 1pt

Our self-justifying 
 evasions of the 
Incompleteness
Theorem are difficult to further extend
primarily because the
combined work of Pudl\'{a}k, Solovay, Nelson and Wilkie-Paris
\cite{Ne86,Pu85,So94,WP87} showed
natural 
Type-S systems cannot recognize  their own Hilbert consistency.
Also,  Willard
 \cite{ww2,ww7,ww9}  
%
% wwooo
%
strengthened earlier results of 
Adamowicz-Zbierski 
\cite{Ad2,AZ1} 
to establish that  natural 
Type-M  system cannot recognize their semantic
tableaux consistency.

\cvs

\medskip

A related
 class of evasions of the
 Second Incompleteness
Theorem
was discovered
in \cite{ww9}.
Let us say 
 $~\alpha~$ 
is
a
{\bf  Type-Almost-M} 
axiom system
iff  $~\alpha~$  can prove
statements \eq{totdefsymba} and 
\eq{totdefsymbm}
as theorems
while treating
{\it  
none of  sentences
 \eq{totdefxs} -- 
   \eq{totdefsymbm}
as axioms.} 
(Many axiom systems,
that 
use
 function symbols ``$~+~$'' and
``$~*~$'' for
formalizing addition and multiplication, 
fall technically into the
% obviously  
Type-Almost-M
% (rather than Type-M) 
category.)
\vspace*{- .2 em}
\beq 
\label{totdefsymba}
\forall x ~\forall y~ \exists z ~~~~x+y=z
\enq
\vspace*{- 1.5 em}
\beq 
\label{totdefsymbm}
\forall x ~\forall y ~\exists z 
 ~~~~x*y=z
\enq
\vspace*{-1.4 em}

\noindent
The preceding is of interest because
some surprisingly strong 
(albeit unusual)
Type-Almost-M systems
\cite{ww9}
 have an ability to 
verify their Herbrand but not
also semantic tableaux consistency.

The 
proofs
in our prior        papers
were 
challenging
primarily
 because 
they required
one to separate the
local combinatorial
methods
employed in
\cite{ww93,ww5,wwapal,ww9}'s 
 particular applications 
from
the 
common principles
that underlied behind all
these 
works.
Our 
Theorems \ref{ppp1}, \ref{ppp2} and \ref{pqq4}
will  rectify this problem by
identifying 
common components
that unite these
four paradigms.
(Theorems 
 \ref{pqq3}, \ref{pqq5},
  \ref{ppp6}, 
E.1, G.2 and G.3 will then carry on
in further
directions.)

%%%iii 1

All these theorems will 
contain severe limits on their generality,  so that the
Second Incompleteness Theorem does not contradict them.
It is clearly perplexing to imagine
how humans
are able to 
motivate themselves to cogitate, 
without
 their thought processes possessing
some type of
{\it  at least tentative}
presumption of 
their own consistency.
%It is for this reason that our
Our
research
has thus consisted of an approximately equal effort
in exploring 
both 
\cite{ww2,wwlogos,wwapal,ww7}'s
new
% types  of  
generalizations of the
Second Incompleteness Theorem
and 
%in examining the unusual perspectives of 
\cite{ww93,ww1,ww5,ww6,wwapal,ww9}'s unusual
boundary-case exceptions to it.

It is clear 
every boundary-case exception
to the  Second Incompleteness Theorem
has limited scope because
the  Incompleteness Theorem
is a broadly encompassing result.
This
 paper
 will, thus, 
be addressing
a challenging 
near-paradoxical
 question
about the maximal 
nature of 
self-justification
that
 can
 never be resolved
in a
fully satisfying
manner.
The Second
 Incompleteness Theorem is
clearly sufficiently
central to  
logic
for it to be desirable
to know  
what {\it partial roads of success} a 
self-justifying axiom     system can 
obtain.

%\cvs

%\vspace*{- 0.6 em}

\section{Literature Survey}

\label{survey}
\label{B2-lem}
\label{D2-def}

%\vspace*{- 0.6 em}

Two 5-page surveys of the prior literature about the
Second Incompleteness
Theorem were provided in our
articles \cite{ww5,wwapal}. 
This section will present a 
more abbreviated survey,
focusing
on 
only 
those
developments that are
particularly germane to
 the current article.

The study of incompleteness
began with 
four classic papers  by
G\"{o}del, L\"{o}b, Rosser and Tarski 
\cite{Go31,Lo55,Ro36,Ta36} and
with the
Hilbert-Bernays
 exploration of 
their 
derivability conditions
\cite{HP91,HB39,Ka91}.
Generalizations of these results for weak
axiom systems, such as Q,  began with 
the work of Tarski-Mostowski-Robinson \cite{TMR53}
and 
Bezboruah-Shepherdson \cite{BS76}.

Some 
more
notation is needed to describe more
recent developments.
Let $~x'~$ denote 
the ``successor'' operation that maps 
$x$ onto  $x+1$.
A formula
 $ \varphi(x) $ is called \cite{HP91} a
{\bf Definable Cut} for an
axiom system $~ \alpha~$   iff
$~\alpha~$ can prove:
\begin{equation}
\label{initdefx}
\varphi(0) \mbox{  AND  }
\forall~x~ \{~\varphi(x)\Rightarrow\varphi(~x'~)~ \}
 \mbox{  AND  }
 \forall~x ~\forall~y<x
 ~~ \{ ~\varphi(x)\Rightarrow\varphi(y) ~ \}
\end{equation}
Definable cuts 
and their cousins 
have been studied by
an
extensive literature
\cite{Ad2,AZ1,Be97,Bu86,BI95,HP91,Ka91,Ko6,Kr87,Ne86,Pa71,PD83,PW81,Pa72,Pu84,Pu85,Pu96,Sm85,Sv83,Sv7,VV94,Vi92,Vi93,Vi5,VH73,WP87}.
(They
are
 unrelated to Gentzen's 
notion of
 sequent calculus
deductive ``cut rule'',
which uses the word ``cut'' in a  
different
context).

\smallskip

A Definable Cut 
$~\varphi(x)~$
is called {\bf Non-Trivial}
relative to
an axiom system $~\alpha~$ iff
$\alpha~$ cannot prove  $~\forall ~x~\varphi(x)~,~$
 {\it although it can prove 
\eq{initdefx}.}
Every axiom system $\,\alpha$,
strictly weaker than Peano Arithmetic, 
will contain some non-trivial Definable Cut.
This cut will have the property that 
$\,\alpha\,$
 can   verify $~\varphi(n)$ for each 
fixed integer $~n$,
although it cannot prove
 $~\forall ~x~\varphi(x)~$.

\smallskip

Let
 $ \, \lceil  \,  \Psi  \,  \rceil  \, $
denote
$   \Psi $'s
G\"{o}del number,
and 
Prf$\, _\alpha^d \, \,(t,p)$
denote
that $   p   $ is a proof
of the theorem $   t   $ from the axiom system $   \alpha   $ 
using $\, d \,$'s
deduction method.
An axiom system
$   \alpha   $ will
then
 be said to recognize its own 
{\bf Cut-Localized  d-Consistency} 
relative to a Definable Cut 
 $ \,  \varphi \,   $
iff $   \alpha   $ can prove:
\begin{equation}
\label{dcon}
\forall~p~~~~~~ \{ ~~~\varphi(p)~~\Rightarrow~~
\neg~\mbox{Prf}\,_\alpha^d \, ~(~\lceil 0=1 \rceil~,~p~)~~~ \}
\end{equation}
The recent literature 
has 
sought to identify which
triples $   (   \varphi   ,   d   ,    \alpha   )   $ have this property.
A crucial negative 
 result about
 Cut-Localized  d-Consistency, discovered by
 Pudl\'{a}k
\cite{Pu85},
established a 
significant 
generalization of G\"{o}del's Second Incompleteness Theorem.
It
 showed 
that
an axiom system
$   \alpha   $ must be unable to prove
\eq{dcon}'s  statement 
about 
any of its definable cuts  
 $\,\varphi \,$, $\,$when
$   d   $  represents
 Hilbert deduction and 
$   \alpha   $ is  any consistent extension of $   Q   $.

Solovay \cite{So94} noted how 
Pudl\'{a}k's
result could be combined
  with the techniques of Nelson and Wilkie-Paris
\cite{Ne86,WP87} to establish the following 
theorem
 that will often be
cited in this paper:
\begin{quote}
  \baselineskip = 1.2 \normalbaselineskip 
%notre \small
\small
\baselineskip = 1.25 \normalbaselineskip 
\baselineskip = 1.25 \normalbaselineskip 
{\normalsize \bf  Theorem 2.1 }
{\it 
(Solovay's  1994 Generalization \cite{So94}
of a 1985 theorem of Pudl\'{a}k \cite{Pu85} 
using some of 
Nelson and Wilkie-Paris \cite{Ne86,WP87}'s
methods )} :
Let $ \, \alpha \, $ denote any axiom system which
contains 
\ep{totdefxs}'s
Type-S statement
and which assures 
%that 
the successor operation
always satisfies 
% the axioms of 
 $  \,   x'     \neq 0     $ and
$     x'     =     y' \Leftrightarrow x=y $.
$~$Then $~\alpha~$
will be unable to recognize its
own  Hilbert
consistency,
whenever  it
 treats addition and multiplication
as 3-way relations 
satisfying 
their usual identity,
associative, commutative and 
distributive  properties.
\end{quote}
Solovay never
published any
precise
proof of Theorem 2.1's hybridizing of
the work of Pudl\'{a}k, Nelson and Wilkie-Paris 
\cite{Pu85,Ne86,WP87},  which he 
privately
communicated
 \cite{So94} to us.
A reader can find 
generalizations of the Second Incompleteness Theorem 
that are closely 
related to
Theorem 2.1
in papers by 
Pudl\'{a}k, Buss-Ignjatovic,
\v{S}vejdar and Willard
\cite{BI95,Pu85,Sv7,wwlogos}, as well
as in 
Appendix A of 
\cite{ww1}.

Other interesting observations
that 
preceded our research were that
Wilkie-Paris \cite{WP87} 
demonstrated  that
I$\Sigma_0+Exp$ cannot prove the
Hilbert consistency of even the axiom system $~Q, \,~$
and that Adamowicz-Zbierski \cite{Ad2,AZ1} showed that
I$\Sigma_0+\Omega_1$ 
satisfied the Herbrandized version
of the Second Incompleteness Theorem.
Both 
these results helped  stimulate 
\cite{ww2,ww7}'s
 semantic tableaux generalizations
of  
the Second Incompleteness Theorem for
I$\Sigma_0$.

A fascinating observation by 
L. A.   Ko{\l}odziejczyk
\cite{Ko5,Ko6},
about the difference in lengths between
semantic tableaux and Herbrandized proofs, 
also motivated
our
investigation \cite{ww9}
into some surprising properties  of
unorthodox encodings for I$\Sigma_0$.
  Ko{\l}odziejczyk 
observed
\cite{Ko5,Ko6}
that various generalizations of
the Second Incompleteness Theorem
for  I$\Sigma_0$ and
I$\Sigma_0+\Omega_i$ in
\cite{Ad2,AZ1,Sa11,ww2,ww7,ww9} 
imply
that the proof of 
 the Herbrandized version
of the Second Incompleteness Theorem
can be
more complicated 
than
 its semantic tableaux counterpart.
This is because there can be an exponential difference
between semantic tableaux and Herbrandized proof lengths
under extremal circumstances.
It was due to 
Ko{\l}odziejczyk's
insightful 
 communications
\footnote{ \baselineskip = 1.3 \normalbaselineskip \b55 
 \label{putglue} The Herbrandized and semantic tableaux
definitions of an axiom
 system $\alpha$'s consistency 
 are 
different 
from each other
because the former
requires  skolemizing $\,\alpha \,$'s axioms, $\,$while the
latter permits
\cite{Fi90} 
an existential quantifier
elimination rule
to
replace Skolemization.
This distinction   can
create a
 potential
exponential difference
 between
the lengths of
Herbrand and Semantic Tableaux proofs.
This 
insightful observation,
%distinction,
due to private communications from
L. A.  Ko{\l}odziejczyk \cite{Ko5},
was used in
 \cite{ww9}
to create an axiom system that satisfied the
semantic tableaux but not
also Herbrandized version of the
Second Incompleteness Theorem.}
that \cite{ww9}
developed 
an axiom system that was a boundary-case
exception to the Herbrandized but not also 
the
semantic tableaux version
of the Second Incompleteness Theorem.

The literature on Definable Cuts has
centered its evasions of 
the Second Incompleteness Theorem around 
\ep{dcon}'s localization formalism
(rather than  
employing  analogs of
 Section \ref{secc1}'s
SelfRef$(\alpha,d)$ axiom,
as we did in \cite{ww93,ww1,ww5,ww6,wwapal,ww9}).
Pudl\'{a}k \cite{Pu85}
proved that essentially
every axiom system
$    ~     \alpha   ~      $ 
 of
finite cardinality can be associated with
a definable cut  $  ~      \varphi   ~      $ such that 
$    ~     \alpha   ~      $ 
can prove
sentence \eq{dcon}'s  validity for $    ~     \varphi   ~  $
when $    ~     d   ~      $ is either
the semantic tableaux or Herbrand-styled
deductive method.

Pudl\`{a}k's 
theorem is related to 
Friedman's
observation \cite{Fr79b}
that 
for many finite theories
$S$ and $T$,
the theory $  S  $ has an interpretation in $  T  $
if and only if
I$\Sigma_0+Exp$ can prove that
$T$'s Herbrand consistency implies
$S$'s Herbrand consistency.
Several generalizations of these results 
by 
 Kraj\'{i}cek, 
Pudl\`{a}k,
 Smory\'{n}ski,
\v{S}vejdar 
  and Visser appear  in
\cite{Kr87,Pu85,Pu96,Sm85,Sv7,Vi92,Vi93,Vi5}.
Visser's article \cite{Vi5} 
contains 
an excellent review of this literature,
as well as
many 
additional
new results.
Also,
we will see how 
% some of
some of the reflection machines of
Beklemishev,
Kreisel-Takeuti and
 Verbrugge-Visser
\cite{Be95,Be97,Be3,KT74,VV94,Vi5}
nicely complement
\thx{ppp6}'s
 reflection mechanisms
in alternate types of intended applications.

It was established by 
H\'{a}jek,  
\v{S}vejdar and
Vop\v{e}nka \cite{Sv83,VH73}
that GB Set Theory can
construct
a definable cut $\varphi$ where
it can prove
the statement  \eq{dcon} is  valid when
$~d~$ denotes Hilbert deduction and
$~\alpha~$ 
is ZF Set Theory.
This result was
 surprising because Pudl\'{a}k \cite{Pu85}
showed GB can never verify its own 
Hilbert consistency 
localized
on a definable cut. (Thus, GB will  view its Hilbert consistency as
equivalent to ZF's Hilbert consistency
in a global sense 
 {\it but not} in a  cut-localized
respect.).

%% 
%% H\'{a}jek,  
%% \v{S}vejdar and
%% Vop\v{e}nka \cite{Sv83,VH73}
%% constructed
%% a definable cut $\varphi$ where
%%  GB Set Theory
%% can prove
%% that the statement  \eq{dcon} is  valid when
%% $~d~$ denotes Hilbert deduction and
%% $~\alpha~$ 
%% is ZF Set Theory.
%% This is
%%  surprising because Pudl\'{a}k \cite{Pu85}
%% proved 
%% GB can never verify its 
%% Hilbert consistency 
%% localized
%% on a definable cut.
%% Thus, GB will  view its Hilbert consistency as
%% equivalent to ZF's Hilbert consistency
%% in a global sense 
%%  {\it but not} in a  cut-localized
%% respect.
%% 

\bigskip

In some sense, 
 Kreisel and Takeuti
  \cite{KT74,Ta53}
can be viewed as
the 
first authors to develop a
%%type of   system 
logic
recognizing its own
 consistency
using a
variant of \ep{initdefx}'s
formula.
Their
results
for typed logics
%  have  
formalized a second-order generalization of 
Gentzen's sequent calculus
that can verify
% a declaration of
 its own consistency, 
%in the special case where 
when
no sequent calculus deductive
cuts are performed.
A key aspect of their formalism
can be seen as using
an analog of
\eq{dcon}'s sentence in an implicit 
manner. 
It thus begins
by using a set of objects, which we shall call $~I,~$
that includes all the standard integers 
plus  some
allowed
 {\it non-standard integers}
(that can 
permissibly
represent 
contradictory  proofs).
  %, and
Their second-order logic
 then uses Dedekind's definition of the natural numbers to 
construct a 
subset of $~I~$, $\,$called
 ``$~N~$'', which
includes all the standard integers
and which is disallowed to contain any
contradiction proof.
We will not go into the details here, but this
transition from $\,I\,$ to $\,N\,$ (with an 
accompanying
relativization of the provability predicate onto
$\,N\,$'s more restricted domain) can be viewed
as Kreisel-Takeuti's
analog
 of %    sentence  
\eq{dcon}'s 
 local consistency 
statement
for
\cite{KT74,Ta53}'s
``CFA''
 second-order logic.

It is difficult to compare our   research
(which 
has
relied upon an analog of 
SelfRef$(\alpha,d)$'s
Kleene-like
{\it ``I am consistent''} axiom)
with the
preceding literature
that has
used 
various forms of
Localized  d-Consistency 
statements.
This is because every effort to 
evade the Second Incompleteness
Theorem
employs
some built-in weakness
to evade G\"{o}del's classic paradigm.

\smallskip

Our 
work
in \cite{ww93,ww1,ww5,ww6,ww9}  represented
less than a 
full-scale evasion of the
Second Incompleteness Theorem
mostly
because 
it was
 incompatible with treating 
{\it as  formal
axioms}
the statements in
Equations \eq{totdefxm}
and \eq{totdefsymbm}
that
 multiplication is a total function
\footnote{ \b55 This 
caveat applies
also  to our article \cite{ww9},
although its
Herbandized form of
self-justification
differs from our other papers
by retaining a capacity to treat 
\eq{totdefsymbm}'s
statement about multiplication's totality as a
derived theorem {\it that is
not an} 
%  a formalized  
axiom.
The key point is that  
theorems are
 weaker
than  axioms
under 
Herbrand deduction 
 because only
axioms are used as intermediate steps 
during proofs.
This
explains intuitively
how    \cite{ww9}'s formalism
was able to recognize
its 
%own 
Herbrandized
consistency, $\,$while 
treating
\eq{totdefsymbm}'s statement about
the totality of
   multiplication 
 {\it
as a theorem}.
(We will return to this subject in 
Appendix D.) } .
%(More will be said about this subject at the end of
%Appendix D.) } .
Some 
reasons
why
it was helpful for
\cite{ww93,ww1,ww5,ww6,ww9}
to employ 
 analogs of    
Section \ref{secc1}'s
SelfRef 
axiom are that:
\bed
\topsep -3pt
\itemsep +2pt
\item[A  ]
An axiom's self-referential
{\it ``I am consistent''} declaration
allows a formalism to recognize its consistency
in 
a
global sense, rather than in the 
$\varphi -$localized sense
used by
sentence
\eq{dcon} and 
the analogous 
  Kreisel-Takeuti relativization of their
second-order
proof predicate.  
\item[B  ]
If a logic is employing a deductive method
$d$ that lacks a modus ponens rule, as occurs in
nearly all
self-justifying
systems, then it is
          preferable for  it           to 
view its {\it ``I am consistent''} 
statement as an axiom rather
than as a theorem.
(This is
 because 
weak deductive methods are
 capable of
drawing 
logical inferences
only 
from
axioms
when
modus ponens 
is absent.) 
\item[C  ]
Analogs of
Section \ref{secc1}'s
SelfRef$(\alpha,d)$'s
{\it  ``I am consistent''} axiom
have been shown
by 
\cite{ww93,ww1,ww5,ww6,wwapal,ww9}
to at least partially
formalize the notion of a 
logic  
possessing  {\it  an 
 almost}
instinctive  
form of 
faith in its
own internal consistency.
(This paper will make it apparent that
such an instinctive faith is less than a full-scale
proof. Yet, Theorem \ref{ppp6} and
 Remarks 
\ref{f88},
\ref{remhappy} 
and \ref{recc1}
will make it apparent that such formalizations
of instinctive faith are also
useful.) 
\ennd
We emphasize
that both  virtues and drawbacks of
SelfRef$(\alpha,d)$'s
{\it  ``I am consistent''} 
axiom statements 
have been 
cited
in this paragraph because
every effort to
evade the Second Incompleteness
Theorem can obtain no more than
limited levels of
success.

% \nop

\cvmew

The scope of the challenge we face becomes apparent when one
realizes 
 $\alpha \, + \,$SelfRef$(\alpha,d)$ 
is 
 inconsistent 
for most
$(\alpha,d).~$
  This is because
   $\alpha+$SelfRef$(\alpha,d)$ 
typically
satisfies Part-i
  {\it but not also }
  Part-ii of
Section 
\ref{secc1}'s  definition 
 of a ``self-justifying'' logic.
(Thus, 
a
% G\"{o}del-like 
diagonalization
paradigm will typically 
imply 
$\alpha \, + \,$SelfRef$(\alpha,d)$ 
is inconsistent,
$\,$as a 
consequence of 
it containing $ ~$SelfRef$(\alpha,d)~$ 
as an axiom.)
This is the
reason      Kleene, 
 Rogers and Jeroslow \cite{Kl38,Ro67,Je71}
were 
hesitant
%skeptical
about 
the utility of
SelfRef$(\alpha,d)$'s
mirror-like
axiom sentence.
Our goal 
in \cite{ww93}-\cite{ww9} 
has 
been
 to 
develop generalizations and boundary-case exceptions for
the Second Incompleteness Theorem, so as 
determine exactly
which paradigms
%%\
%%\{\it 
%%\paradigms
%%\are sufficiently
%%\weak}
%%\so that 
%%\they 
%%\
can
support, for example,
Theorem \ref{ppp6}'s limited notion
of self-justification.

The 
%intuitive
 reason 
one would  anticipate some 
{\it limited}
%partial
exceptions to the Second Incompleteness Theorem 
to
exist is 
it is  hard to 
imagine how humans can motivate
themselves to cogitate without 
using 
some
variant 
of 
self-justification.

%%
%%Thus, our %recent 
%%research 
%%project
%%has sought to {\it both} generalize
%%  the Second Incompleteness
%%Theorem and explore when it can be evaded.

\cvl

\section{Generic Configurations}

\label{3uuuu1}

% \vspace*{- 0.2 em}

%Our main results about self-justification will
%appear in Sections 
%\ref{3uuuu2} -- \ref{sect64}.
The
 phrase  {\bf Bounded Quantifier}
will refer
to
expressions of the form
``$\, \exists~v  \leq  T \,$'' or
``$ \,\forall~v  \leq  T\,$''
where $T$ is a term.
A formula  
is
called
{\bf Fully-Bounded$\,$} when all its quantifiers are
   so bounded.
\lem{lex22} will soon explain how Definition \ref{xd+1x1}'s 
formalism can encode conventional arithmetic:

\medskip

\bxbxdd
\label{xd+1x1}
\rm
Let $~\xi~$ denote some 
non-integer
indexing superscript
(whose  properties will 
be 
discussed
 later by
Definition
\ref{def3.3}). 
Then the symbol
$~\Delta_0^\xi~$
will 
denote some 
fixed
special
set of fully-bounded formulae
that
is
closed 
under negation, 
in a 
language
that will be later
 called $~L^\xi~$.
(Thus, if some
formula $~ \Psi~$ is a member of
             $~\Delta_0^\xi~$
then so is 
$~ \neg ~ \Psi~$.)
Items 1-3 
formalize how
$~\Pi_n^\xi~$ 
and 
$~\Sigma_n^\xi$ formulae are
built 
in a straightforward manner
out of
these
 $~\Delta_0^\xi~$
sub-components:
\bee
\topsep -15pt
\item
Every $ \, \Delta_0^\xi \, $ formula is
considered to be also a
$ \, \Pi_0^\xi \, $ 
and 
$ \, \Sigma_0^\xi$ formula.
\item
For $~n \, \geq \, 1~~,~$
a formula 
will be 
called
$ ~~  \Pi_n^\xi~~$
iff it can be
written 
in the canonical form of
$\forall v_1 ~ \forall v_2 ~...~ \forall v_k ~~~ \Phi(v_1,v_2,..v_k \,),~~~$  where
$\Phi$ is  $~\Sigma_{n-1}^\xi.~$
\item
Likewise for $~n \, \geq \, 1~~,~$
 a formula will be
called
$~~ \Sigma_n^\xi~~ $
iff it can be
written in the 
form of
$~\exists v_1 ~ \exists v_2 ~...~ \exists v_k ~~~ \Phi(v_1,v_2,..v_k \,),~~~$
  where
$\Phi$ is  $~\Pi_{n-1}^\xi.~$
\ene
\eedd

\medskip
\noindent
{\bf Notation Convention:}
Our rules for defining $\,\xi$, specified later in this section,
will never have this  superscript 
designate
an integer quantity. This is because 
integer superscripts have a special meaning under 
a typed-based hierarchy,
not intended here.

\nop

\begin{example}
\label{ex3-1}
\rm
Let
$~L~$  denote a 
conventional arithmetic
language that uses
function symbols 
for denoting 
addition and multiplication.
Below are two
examples  
of $\Delta_0-$like formulae that
invoke
Definition \ref{xd+1x1}'s
notation:
\bed
\itemsep 1pt
\item[  a ] The 
symbol ``$~\Delta_0^A~$''
will denote any 
fully bounded
formula that uses the addition, multiplication and maximum
function symbols in an arbitrary manner.
(Thus,  $~\Delta_0^A~$ 
corresponds to what 
many
textbooks
\cite{HP91,Ka91,Kr95}
simply call a ``$~\Delta_0~$'' formula.)
\medskip
\item[  b ] The 
symbol ``$~\Delta_0^{R}~$'' 
will denote a class of 
formulae in $L$'s language
whose bounded quantifiers are  allowed
to use
only 
  the Maximum function symbol. 
Their
bodies,
however, may contain
any combination of 
addition, multiplication
and maximum function symbols.
\ennd
Formulae \eq{ex1}  
and \eq{ex123} 
illustrate the
distinction
between the  $\Delta_0^A$ 
and  $\Delta_0^{R}$ classes. Thus, \eq{ex1} satisfies the first
but not second condition (on account of the presence of the
multiplication symbol used by its bounded quantifiers).
 In contrast,
\eq{ex123}
is an example of a  
 $~\Delta_0^{R}~$ formula.
\beq
\label{ex1}
\exists \, y\leq  x*x ~~~~
\forall  \, z\leq  y*y  ~~~~
\exists \, w\leq  y*z ~~~~
~~~~~: ~~~~~ 
 \{ ~~~~ x  *y=z+w ~~~~ \}
\enq 
\vspace*{-0.5 em}
\beq
\label{ex123}
\exists \, y\leq  x ~~~~
\forall  \, z\leq  y  ~~~~
\exists \, w\leq  \mbox{Max}(y,z) ~~~~
~~~~~: ~~~~~ 
 \{ ~~~~ x  *y=z+w ~~~~ \}
%% x  *y=z+w
\enq 
The distinction between 
 $ \Delta_0^{A}  $
arithmetic formulae and the unconventional
 $ \Delta_0^{R}  $ class may first convey the impression
that these
two 
classes 
have
fundamentally 
different natures. Actually,
\lem{lex22} will
show
that their
relationship is more subtle.
This is because
its formalism
will map  $ \Delta_0^{A}  $ formulae
onto     $ \Delta_0^{R}  $
expressions that are equivalent to it under Definition
\ref{snn}'s Standard-M model
$~$---$~$ in a context where
only the length of these
formulae
is allowed to 
possibly
grow.
This 
% result
% Standard-M 
equivalence
 enabled
 \cite{ww9} to construct
a natural axiomatic formalism
that could recognize its
own
Herbrandized consistency  
{\it but which nevertheless  satisfied} 
the 
idealized form   
\footnote{  \b55 An axiom system $\alpha$ 
is defined to
 satisfy
the 
{\bf ``idealized form''}
of the
semantic tableaux
version 
 Second Incompleteness Theorem
when no  $\beta  \supseteq  \alpha$ can
prove a
 semantic tableaux
proof
 of 0=1
 from 
itself is incapable of existing.
We will summarize
\cite{ww9}'s formalism
and the distinction between
Herbrandized and semantic tableaux 
deduction 
at the 
of 
Appendix D.} 
of
the 
semantic tableaux
version of the 
Second Incompleteness Theorem.
\end{example}

\medskip

\begin{definition}
\label{snn}
{\bf ``Standard-M''} will denote
the
standard model of integers.
\end{definition}

\medskip

The reason for our interest in Standard-M is
that
many 
pairs of formulae
are equivalent
under the Standard-M
model, while
weak axiom systems often
cannot formally prove they are
equivalent.
For instance,  this 
will occur
when  Example \ref{ex3-2}
examines
Definition \ref{def3.3}'s properties.

%\notre-only \medskip

\nop

\bxbxdd
\label{def3.3}
A  {\bf Generic Configuration},
often identified by 
 the superscript symbol of
 $~\xi~$, is defined to be a 5-tuple
$~(\, L^\xi \, , \, \Delta_0^\xi \,  , \, B^\xi \, , \, d \, , \, G \, )~$
where:
\bee
\topsep -12pt
\rm
\item
$ L^\xi $  
is a 
language that includes 
logical symbols
for ``$0$'',  ``$1$'',  ``$2$'',
``$=$'', ``$\leq$'' 
 and
for the operation of ``Maximum(x,y)''.
$ L^\xi $   also includes a 
 sufficient number of 
function  and constant symbols 
so that every integer $~k~$ can be encoded by
some term 
$~T_k~$ 
specifying
$~k\,$'s value.
\item
$~ \Delta_0^\xi ~$ 
corresponds to 
any variation
of
Definition \ref{xd+1x1}'s 
class of 
``fully-bounded'' formulae
that is rich enough to assure that there exists
 two  $~\Delta_0^\xi~$ 
formulae, henceforth
called ``Add$(x,y,z)$''
and ``Mult$(x,y,z)$'', $\,$for formalizing 
the graphs of
addition and multiplication.
(It will  generate 
$\xi$'s set of $\Pi_n^\xi$
and $\Sigma_n^\xi$ sentences,
using Definition \ref{xd+1x1}'s 3-part formalism.)
\item 
$ \, B^\xi \,$  
denotes 
a$~$ {\normalsize {\bf ``Base Axiom System''}}, 
$~ \,$whose axiom-sentences are true under
the 
Standard-M model
and which is
$\Sigma_1^\xi$ complete.
(Thus,
$  B^\xi $  can  prove every true $\Sigma_1^\xi$ sentence, and
it can
likewise refute all false $\Pi_1^\xi$ sentences.)
\item
$~d~$   denotes  $ \, \xi \,  \, $'s  
 method of deduction.
It is required to be sufficiently conventional to
satisfy the usual indirect-implication 
property \footnote{  \b55   
$~~$This is that irregardless of whether or not
$~d~$ contains a built-in
modus
ponens rule, it does support some
form of
a (possibly
quite 
lengthy)
 proof of a
theorem $Z$,
when it is able to prove
$X$, $Y$ and $~(X \wedge Y)~\rightarrow ~Z~$ as
theorems. }
associated with
 G\"{o}del's Completeness Theorem. 
\item
$~g~$   denotes  a method for
encoding the  G\"{o}del numbers
of proofs.
\ene
\eedd

\medskip

\begin{example}
\label{ex3-2}
\rm
Let us recall
Example \ref{ex3-1} defined
``$~\Delta_0^A~$'' 
as essentially the conventional
textbook notion
\cite{HP91,Kr95}
 of an arithmetic 
 ``$~\Delta_0~$'' formula.
This example  will  outline how
well-known techniques can map
every $\Delta_0^A$  formula 
onto a semantically 
equivalent 
 $\Delta_0^\xi$ 
formula under the Standard-M model.

Our discussion will have
Seq$(x)$   denote a 
function that
maps 
non-negative integers onto 
binary strings
in
lexicographic order.
Thus Seq$(x)$  maps
0 onto the empty
string, 
$\,$the integers
1 and 2 onto the strings of ``0'' and ``1'', $\,$the integers 3--6
onto$\,$  ``00'',   ``01'',   ``10'',   ``11'', etc.   
(Formally,
Seq$(x)$
is
an operation
that maps integer $~x~$ onto the bit-string 
that occurs to the immediate right of the leftmost ``1'' bit in
the binary encoding of $~x+1.~~)$

Given any $\,k-$tuple  $~(~x_1\, , \, x_2\, , ~ ... ~ x_k~),~$
let  $~$STRING$(   x_1   ,   x_2   ,\, ... \,   x_k   )$
denote the
 concatenation of
Seq$(x_1), \, ... ~$Seq$(x_k)$.
For any 
integers $v$ and $w$ satisfying $v  \leq w^2$, 
it is clear that 
there  exists $(x_1,x_2,x_3)$ where
STRING$(x_1,x_2,x_3)$ 
represents $v  $'s
binary encoding and
each $x_i  \leq   \mbox{Max}(w,4)$.

\smallskip

An example will now illustrate
the approximate structure of an inductive methodology
for mapping
$\Delta_0^A$ formulae 
onto their equivalent   $\Delta_0^\xi$  counterparts
in the Standard-M model.
Let SQUARE$(x_1,x_2,x_3,w)$ be a 
$\Delta_0^\xi$ formula which 
specifies 
that
STRING$(x_1,x_2,x_3)$ represents  an integer 
$ \leq \, w^2~.~ \,$ Also,  
$\,$let $\phi^*(x_1,x_2,x_3)$
and $\phi(v)$ represent 
a pair of  $\Delta_0^\xi$ and
$\Delta_0^A$ formulae 
that
are equivalent
under the Standard-M model
$\,$when STRING$(x_1,x_2,x_3)$ is an encoding for  $~v~$. 
Then one
possible method
for mapping $\Delta_0^A$ formulae onto their 
equivalent $\Delta_0^\xi$ 
counterparts 
(in the Standard-M model)
could 
map the
formula \eq{exam21} 
onto 
\eq{exam22}'s alternate form:
\beq
\label{exam21}
\forall ~ v\, \leq \, w^2~~~~~\phi(v)
\enq
\smallskip
\begin{center}
$\forall ~ x_1 \leq~ \mbox{Max}(w,2)~~~~~\forall ~ x_2 \leq ~ \mbox{Max}(w,2)
  ~~~~~\forall ~ x_3 \leq ~ \mbox{Max}(w,2)$
\end{center}

\vspace*{-0.5 em}

\beq
\label{exam22}
  \{~~ \mbox{SQUARE}(x_1,x_2,x_3,w)
  ~~ \Rightarrow~~
\phi^*(x_1,x_2,x_3)~~\}
\enq
\lem{lex22}
indicates
Example \ref{ex3-2}'s
translational methodology generalizes easily to
       all combinations of   $\Delta_0^A$  inputs
and
generic configurations $\, \xi \,$,
via
an 
approximate
inductive generalization of the transition from
sentence \eq{exam21} to \eq{exam22}. (Its 
procedure essentially
performs
iteratively a finite number of such transitions, so as to
translate   all 
the
clauses of
 an
initial
$\Delta_0^A$ 
% format
formulae
into their 
$\Delta_0^\xi$ 
counterparts via an inductive methodology.
The intuition behind these transitions is 
they will repeatedly replace a single variable,
such as $~v~$ in 
sentence \eq{exam21},
with a multiplicity of variables,
such as $(x_1,x_2,x_3)$ in 
\eq{exam22}.)
\end{example}

\smallskip

\begin{lemma}
\label{lex22}
{\rm (Paris-Dimitracopoulos
\cite{PD82} )}
For every generic configuration $\xi$, each
$\Delta_0^A$ formula
can be translated
into an  equivalent 
$\Delta_0^\xi$ 
formula
in the
Standard-M
model
via a generalization of Example \ref{ex3-2}'s
process.
{\rm (This 
 clearly
 implies 
$\Pi_j^A$ 
formulae
can
also be
translated
into
$\Pi_j^\xi$ 
expressions.)}
\end{lemma}

\smallskip

Paris-Dimitracopoulos \cite{PD82}
sketched   
an 
%approximate
analog of Lemma \ref{lex22}'s
translation 
algorithm,
using
only
 slightly different notation,
 that
is applicable to 
any formalism that satisfies Parts (1) and (2)
of 
 Definition \ref{def3.3}.
Their formalism
thus uses an 
inductively-iterated
analog of the prior example's
replacement of a single variable $~v~$ in formula
\eq{exam21} with \eq{exam22}'s multiplicity of variables 
% in the tuple 
$~(x_1,x_2,x_3)~,~$ 
so as to perform
 Lemma \ref{lex22}'s translation 
task.
It will  be unnecessary for a reader to 
consider the
details behind 
% either 
\cite{PD82}'s Theorem 1 or
Lemma \ref{lex22}'s 
similar translation 
mechanism
%%% 
%%% 
because
 the remainder of this article will never
% explicitly
 use
them again.
%its mechanism.
Instead, their sole purpose 
% of  Lemma \ref{lex22}
has been
% merely
 to provide an {\it implicit backdrop}
for our results by illustrating how
the study of the $\Pi_1^\xi$ sentences of
Definition \ref{def3.3}'s generic configurations
% also
provides information
about   $\Pi_1^A$ sentences 
(after the needed translating is done).

\smallskip

Four  examples of 
self-justifying
systems that employ  
Definition \ref{def3.3}'s 
 $\Pi_1^\xi$ sentences 
 will be
illustrated in 
\appD.
These examples   are 
too
complicated
to be 
examined
 before Sections \ref{3uuuu1} -- \ref{sect64}
are 
read.
However, the next example 
should convey some useful  intuitions:

\medskip

\nop

\cvwew

\begin{example}
\label{ex3-3}
\rm
Let
$   x_0,   x_1,   x_2,     ...    $ 
and  $   y_0,   y_1,   y_2,     ...  $
denote sequences defined by:
\vspace*{- 0.7 em}
\beq
\label{zs}
x_0~~~=~~~2~~~=~~~y_0
\enq
\beq
x_{i}~~~=~~~x_{i-1}~+~x_{i-1}
\label{as}
\enq
\beq
y_{i}~~~=~~~y_{i-1}~*~y_{i-1}
\label{bs}
\enq
\end{example}
For $\, i \, > \,0 \,$, 
$\,$let $ \, \phi_{i} \, $ 
and $ \, \psi_{i} \, $ 
denote the  
sentences in 
\eq{as} and \eq{bs}
respectively.
Also, 
 let
  $ \, \phi_{0} \, $ and
$ \, \psi_{0} \, $ 
denote \eq{zs}'s
sentence.
Then
 $ \, \phi_0, \, \phi_1, \, ... \, \phi_n \, $
imply
 $ \, x_n \, = \, 2^{ n+1} \, , \, $ and 
 $ \, \psi_0, \, \psi_1, \, ... \, \psi_n \, $
 imply $ \, y_n \, = \, 2^{2^n} \, $.
Thus, the  latter sequence
grows at a 
faster
rate than 
the former.
Much of our research has used the difference between the
growth rates of 
$   x_0,   x_1,   x_2,     ...    $ 
and  $   y_0,   y_1,   y_2,     ...  $
as a motivating example explaining why 
\ep{totdefxa}'s Type-A axiom systems can
support a stronger form of boundary-case exception
to the semantic tableaux version of the Second Incompleteness
theorem than
can 
Type-M systems.

\parskip 2pt

\smallskip

Let
Log$(\, y_n \,) \, = \, 2^{n} \, $ 
and
Log$(\, x_n \,) \, = \, {n+1} \, $ 
thus
designate the lengths of the binary codings for
$ \, y_n \, $
 and 
$ \, x_n \, $.
Then $ \, y_n\,$'s coding 
has a length
$\,  2^{n} \,  $, which is
 {\it much larger} than
the $  n+1  $
steps that  $ \, \psi_0, \, \psi_1, \, ... \, \psi_n \, $
use to
define its existence.
However,
 $ \, x_n\,$'s 
length  has a 
 smaller
 size of
$ \, {n+1} \, $.
These observations are useful because every proof 
of the 
Incompleteness Theorem
involves 
a 
G\"{o}del
 number $ \, z \,$ 
coding a sentence
that has a  capacity
to self-reference its own definition.
The faster
growing series $y_0,\,y_1,\,,\,...\,y_n$
should  
be intuitively 
anticipated
to have
 this 
self-referencing
capacity because 
 $~y_n\,$'s binary encoding 
has a 
$~2^{n+1}~$  length that 
dwarfs the
size of the $O(n)$
steps 
used
 to define its
value. Leaving aside 
\cite{ww2,ww7}'s
%
% wwoooo
%
many details,
this
fast growth 
explains
roughly
 why many Type-M 
logics
satisfy the semantic tableaux version of 
the Second Incompleteness Theorem.

\smallskip

This paradigm also 
illustrates intuitively
why some 
 Type-A  systems, employing
\cite{ww93,ww1,ww5}'s
semantic tableaux formalism,
can 
represent
boundary-case exceptions to the
 Second
Incompleteness
Theorem.
 This is because such formalisms 
lack access to 
\ep{totdefxm}'s axiom that multiplication is a total function.
(They are
unable,
thus,
 to
easily
 construct numbers $ \, z \, $ that can
self-reference their own definitions
because they have access  only
to the slower growing
              addition primitive.)
In particular
assuming only that each sentence in
the axiom-sequence
  $   \phi_0,   \phi_1,   ...   \phi_n   $
(from  \ep{as} )
requires  a mere 
two bits 
for its encoding, 
the length    $     n+1     $ of 
 $   x_n   $'s
binary encoding 
will be 
smaller
than the
length of its
defining 
 sequence.

\smallskip

This short length for  $   x_n   $
had
motivated 
\cite{ww93,ww1,ww5,ww6}'s
evasion of  
the semantic tableaux 
version
of the
Second Incompleteness Theorem. 
It
suggested that the self-referencing
needed in a
G\"{o}del-like diagonalization argument would stop being
feasible
when \ep{as}'s slow-growing 
$x_1,\,x_2,\,x_3,\,...$
sequence represents the fastest growth that is possible.

One of the several goals in this article 
will be
 to formalize 
a generalizations of
\cite{ww93,ww1,ww5,ww6}'s
self-justifying methodologies by using
Definition \ref{def3.3}'s
generic configurations.
The  proofs
of  our main theorems
will, of course, be
more subtle than
the hand-waving intuitions appearing in this example.
For instance, the combined work of
Pudl\'{a}k, Solovay, Nelson and Wilkie-Paris \cite{Ne86,Pu85,So94,WP87}
(summarized by Theorem 2.1)
raised the subtle issue that
  no Type-S   system
can
prove a theorem affirming its
own Hilbert consistency. 
Another complication 
is that the
\ep{bs}'s implication for proofs that use
the multiplication operative has
different side effects for
Herbrandized and  
semantic tableaux
deduction
(on account of
 Ko{\l}odziejczyk
\cite{Ko5,Ko6}'s
previously mentioned
observations about the potential
exponential difference between the lengths of these
proofs
under extremal
circumstances).

\smallskip

Our
main theorems
will show
that
self-justifying systems, 
using four
deduction methods,
are
capable of  proving 
all of Peano Arithmetic's  $\Pi_1^\xi$ 
theorems.
Interestingly,
self-justification
will
be compatible with \cite{ww5}'s modification
of semantic tableaux deduction, 
that includes a
 modus ponens  rule 
 for 
 $\Pi_1^\xi$  and  $\Sigma_1^\xi$ type sentences.
However,
 \cite{wwlogos}
has shown 
an analogous
  modus ponens rule for
 $\Pi_2^\xi$  and  $\Sigma_2^\xi$  sentences
is incompatible
with self justification.
(Thus, the contrast between our
main 
results and the
Second Incompleteness
Theorem's generalizations will be quite tight.)

\cvl

\section{Five Helpful Definitions and An Informative Lemma}

\label{3uuuu2}

\label{D.4-theorx}

This section will  introduce 
five definitions and 
 prove 
a 
Lemma \ref{lemex4} 
about
self justification.
This lemma 
will be 
% much 
weaker
than
Sections \ref{3uuuu3}
and \ref{sect64}'s
main results.
Its
main
purpose
 will be to provide
a
useful
starting example.

\smallskip

\bxbxdr
\label{xd+1x3}
The symbol
``E$(n)$'' 
will denote some
term 
in Definition \ref{def3.3}'s
language  $ \, L^\xi \, $
that represents the
value  
 $ \, 2^n \, . \, $ 
In using 
this symbol, we  do
not presume that
$ \, L^\xi \, $ possesses a function
symbol for 
the exponent operation.
Thus if  $ \, L^\xi \, $
has only a function symbol for
multiplication,
then  
 E$(n)$ 
could 
designate
the term of $ \, $``$ \, 2*2* \, ... \, *2 \, $''$ \, $ with
$ \, n \, $ repetitions of ``2'' $. \, $ (Alternatively,
 E$(n)$ 
can be
 defined via
applying
 $\, 2^n \,$ iterations
of the successor function to zero,
or by having a
special constant symbol  designating
 $ \, 2^n \, $'s value.
Essentially, any reasonable method can
be used to define  E$(n)$'s value)

\eedd

\cvrew

\nop

\smallskip

\bxbxdr
\label{xd+1x4}
Let
$\,\Uxp \,$ denote
a prenex normal sentence.
Then 
{\bf Scope$_E$($\Uxp,$N)$~$} will  denote
a sentence identical to $~\Uxp~$ except that
every unbounded universal quantifier  ``$~\forall~v~$''
is changed to  ``$ \, \forall \, v \, < \, E(N) \, $'',
and every
 unbounded existential quantifier  ``$ \, \exists \, v \, $''
is changed to ``$ \, \exists \, v \, < \, E(N) \, $''.
(No change is made among the bounded quantifiers within
the
$~\Delta_0^\xi$ 
part of the sentence
$\,\Uxp \,$.) 
For example, if 
$\,\Uxp \,$ denotes
the $\Pi_1^\xi$ 
sentence of
$~\forall \, v_1~\forall \, v_2~...~\forall \, v_k~~~\phi(v_1,v_2,...v_k)~$
then
\eq{scopede} 
illustrates
Scope$_E$($\Uxp,$N)'s form.
Likewise if
$~\Uxp ~$  is
the $~\Sigma_1^\xi~$ 
sentence of
$~\exists \, v_1~\exists \, v_2~...~\exists \, v_k~~~\phi(v_1,v_2,...v_k)~$
then
\eq{scopedx} illustrates
Scope$_E$($\Uxp,$N)'s form.

\vspace*{- 0.4 em}

{ 
\small  
\cvl
\beq
\label{scopede}
\forall ~ v_1~ < ~E(N)~~\forall ~ v_2~ < ~E(N)~~ ...
\forall ~ v_k~ < ~E(N) ~~~~~ : ~~~~~
\phi(v_1,v_2,...v_k)~
\enq
\beq
\label{scopedx}
\exists ~ v_1~ < ~E(N)~~\exists ~ v_2~ < ~E(N)~~ ...
\exists ~ v_k~ < ~E(N) ~~~~~ : ~~~~~
\phi(v_1,v_2,...v_k)~
\enq}
\eedd

{\bf Special Note about Definition \ref{xd+1x4}'s Meaning.} 
If
$\Uxp$ is a 
$\Delta_0^\xi~$ 
sentence then
Scope$_E$($\Uxp,N)$ 
will 
be equivalent to  $~\Uxp~$
for every $N \geq 0$
by definition.
(This is  because 
$\Delta_0^\xi~$  formulae
contain
no unbounded  quantifiers
that undergo change when 
  $\Uxp$ is mapped onto  
Scope$_E$($\Uxp,N).~~)$

\medskip

{\bf More About
this Notation:}
The
potentially lengthy
syntactic object of
``$\,$Scope$_E$($\Uxp,$N)$\,$'' {\it
 will actually 
 not}
be used 
in our physical encodings
of proofs.
Instead, these encodings will
 use the 
more 
desirably
compressed
object of
``$\,\Uxp\,$''
(which has no
possibly bulky
 $E(N)$ term). The {\it sole function} of 
$\,$Scope$_E$($\Uxp,$N)$\,$ 
will be for us to speculate about what 
Boolean value
$\,\Uxp\,$
 {\it would theoretically assume}
(under the Standard-M
 model)
if  $\,\Uxp\,$'s quantifiers
were  modified
 so that their ranges
were changed to be
bounded by $E(N).$
$~\,$(It  turns out that  
Scope$_E$($\Uxp,$N)$\,$'s
finitized quantifier-range
will 
help
%greatly
simplify our
analysis.)

%proofs.)

\medskip

\bxbxdr
\label{xd+1x5}
A
$\Pi_1^\xi$
or $\Sigma_1^\xi$
sentence
 $ \Uxp $ will be called
{\bf Good(N)} 
when the entity
Scope$_E$($\Uxp,N)$ is 
true
under the
Standard-M model
\footnote{\f55  \label{fgood}  A
quite
unusual aspect of Definition \ref{xd+1x5} is that its
Good$(N)$ condition 
has opposite properties when it is applied to
$\Pi_1^\xi$ 
and $\Sigma_1^\xi$ 
sentences
in one particular respect.
This is because for each  $~N,~$ the
 Good$(~N~)$ condition  is weaker than
the
Good$(~\infty~)$ condition 
for $\Pi_1^\xi$ 
sentences, while it is stronger than it for 
 $\Sigma_1^\xi$ 
sentences. (For instance,
$~\forall \, x~\phi(x)~$ is
stronger than
$~\forall \, x <E(N)~\phi(x)~$,
but 
$~\exists \, x~\phi(x)~$ is weaker than
$~\exists  \, x<E(N)~\phi(x)~$.)}.
Also, a set of 
$\Pi_1^\xi$
sentences, 
denoted as $\theta$, is called  Good(N)
iff 
all of its
sentences are Good(N).
 
\eedd

\bxbxdr
\label{chg}
If $\Uxp $ 
is a $\Pi_1^\xi$ sentence then
 $\zzz \Uxp \,)$ 
will denote the largest integer $J$ such
that $\Uxp $ satisfies the 
Good$(J)$ condition.
(It will equal $\infty$ if 
$\Uxp $ satisfies 
Good$(J)$ for all $J$.)
Also,
if $\theta$ is a set of 
$\Pi_1^\xi$ sentences, then
 $\zzz \theta \,)$ 
will denote the largest $J$ where each sentence in $\theta$
is
Good$(J)$.
\eedd

\medskip

{\bf A Very Helpful Start :}
Several 
more
definitions will 
be needed
before Section \ref{sect64}
  can present our strongest results.
The remainder of this section will  illustrate
how the current formalism   is already
 sufficient for introducing a
useful
starting
lemma.

\smallskip

\begin{definition}
\label{deftight} \rm
Let
      \gggen  
again   denote 
a generic configuration called $~\xi~$,
and let us
presume that
its 
 base axiom system
 $~B^\xi ~$ 
is comprised exclusively of 
 $ \, \Pi_1^\xi \, $ sentences.
Also,
let
 $ \, \beta \,  \supset  \,  B^\xi  \, $
      denote a second axiom system,  comprised
also of 
 $ \, \Pi_1^\xi \, $ sentences,
that 
(unlike $ \, B^\xi \,$)
can
possibly  be
inconsistent.
(If  $ ~  \beta ~ $
 is inconsistent then 
 let $ \,  q_\beta  \, $ denote the shortest proof of
$ \, 0=1  \,$
from $~\beta~$.$~$)
Then the
generic
configuration
 $ \, \xi \, $ 
will be called
 {\bf Tight} if 
iff {\it every 
inconsistent} set of 
 $ \, \Pi_1^\xi \, $ sentences
 $~\beta \, \supset \, B^\xi  ~ $
satisfies the following
constraint:
\beq
\label{piss1}
\mbox{ Log(}\,q_\beta \, )  ~ ~ \geq ~~ \zzz \beta ~)
 ~~+~~2
\enq
\end{definition}

Lemma \ref{lemex4} will 
prove 
$B^\xi+$SelfRef$(B^\xi,d)$
satisfies
Section \ref{secc1}'s
 self-justification
criteria
 whenever
$ \, \xi \,$
is tight.  
This ``tightness''
will clearly fail to be satisfied by many
generic configurations.
This is
 because 
the Second Incompleteness Theorem is 
a widely encompassing result,
which imposes severe restrictions on its
allowed
 exceptions.
Lemma \ref{lemex4}'s 
mini-result 
will be of interest 
primarily
because it will be generalized
substantially in
Sections \ref{3uuuu3} 
and \ref{sect64}.

\medskip

\begin{lemma}
\label{lemex4}
If a 
generic configuration
  \gggen  
is tight then
 $B^\xi+$SelfRef$(B^\xi,d)$
will be a consistent
self-justifying axiom system.
\end{lemma}

\medskip

{\bf Proof Sketch:} Our justification of \lem{lemex4}'s
mini-result  will be
simpler than 
 the
next section's 
proof of 
\thx{ppp1}'s stronger result.
 The current proof 
will also be kept
brief and informal because the 
same topic will be visited  more 
rigorously during
Section  \ref{3uuuu3}'s
discourse.

Let $~\Psi~$ 
denote 
Section \ref{secc1}'s
SelfRef$(B^\xi,d)$  sentence.
We will  
omit formalizing
 $~\Psi\,$'s exact 
 $\Pi_1^\xi$ 
encoding here
because Appendix A will 
provide a more general
fixed-point
 construction,
using
Definition \ref{xd+1x7}'s stronger paradigm.
%% 
%% Let $~\Psi~$ 
%% denote 
%% Section \ref{secc1}'s
%% SelfRef$(B^\xi,d)$  sentence.
%% We will  
%% omit formalizing
%%  $~\Psi\,$'s exact 
%%  $\Pi_1^\xi$ 
%% encoding 
%% here because Appendix A will later
%% provide an encoding for a more general fixed-point
%% sentence in the context of
%% Definition \ref{xd+1x7}'s stronger paradigm.
%% 
The
current
proof
will 
simply 
presume 
 $~\Psi\,$'s fixed point statement 
can receive a
 $\Pi_1^\xi$ 
encoding
under
sentence
\eq{fixnew},
where 
``$~\mbox{Prf}^d_{~B^\xi+{\rm SelfRef}(B^\xi,d)}(~\lceil \, 0=1 \, \rceil~,~p~)~$''
is a $~\Delta_0^\xi$  formula, indicating that
$~p~$ is a proof
of 0=1 from 
$~B^\xi~+~$SelfRef$(B^\xi,d)~$
under  $~d \, $'s
deduction method.
\beq
\label{fixnew}
\forall ~p~~~\neg ~ 
~\mbox{Prf}^d_{~B^\xi+{\rm SelfRef}(B^\xi,d)}(~\lceil \, 0=1 \, \rceil~,~p~)
\enq

%\cvs

The
Definition \ref{chg}'s symbol 
 $~\sharp~$ 
will be
% become 
helpful
at this juncture. The application of
$\xi$'s tightness to 
% sentence
\eq{fixnew}'s
$\Pi_1^\xi$ styled encoding  
will 
%thus 
imply
 \footnote{  \b55 \g55 \ep{piss2} is 
easy to justify when
one  presumes
there is an
available 
$\Pi_1^\xi$ encoding of
\eq{fixnew}'s statement,
which we call  $ \, \Psi \, $.
(This 
 $\Pi_1^\xi$ 
presumption
is reasonable because an analog of
$ \, \Psi\,$'s
exact
  $\Pi_1^\xi$ encoding 
will be discussed later by  Definition \ref{xd+1x7}
and 
in
Appendix A.)
The invalidity  
of $ \, \Psi \, $
% in this context,  
will 
thus assure the existence of 
a
proof of $0=1$ from
 $ \, B^\xi \, + \, $SelfRef$(B^\xi,d)$. Moreover, 
Definition \ref{chg}'s notation  implies
$\zzz \Psi  \, )$ will equal Log$(q) \, - \, 1 \, $ when 
$q$ denotes the  
 shortest proof of $0=1$ from
 $ \, B^\xi \, + \, $SelfRef$(B^\xi,d)$. 
The latter shows \eq{piss2} is
valid.} that
\eq{piss2} must be true 
when
 $~\Psi~$ 
is 
false
under the Standard-M model
and
when
the shortest proof of $0=1$ from
 $~B^\xi~+~$SelfRef$(B^\xi,d)$
is denoted as $~q$.
\beq
\label{piss2}
\mbox{ Log(}\,q~)  ~~ = ~~ \zzz \Psi ~) ~~+~~1
\enq
Also
 Definition \ref{chg}
trivially implies 
$~\zzz B^\xi+\Psi ~) 
~ = ~ 
\zzz \Psi ~)~$ 
(because  all of
 $~B^\xi \,$'s axioms are
true  under the Standard-M model).  Thus, 
\eq{piss2}
yields
\eq{piss3}.
\beq
\label{piss3}
\mbox{ Log(}\,q~) ~~ = ~~ 
\zzz B^\xi+\Psi ~) 
 ~~ +~~1
\enq
But the point is that
the 
Tightness constraint's  \ep{piss1},
$~$used in the 
%%nope particular 
context where $~\beta~$
is the axiom system
of $~ B^\xi+\Psi ~,~$
implies
%%nope  that
$~\mbox{ Log(}\,q \, ) ~ \geq ~ \zzz B^\xi+\Psi ~)~+~2~$.
This  directly contradicts \ep{piss3}'s equality, 
{\it whenever} the proof  ``$~q~$'' of $~0=1~$ cited 
in our discussion 
{\it does formally exist.}
This 
contradiction
%latter observation
is precisely what 
is needed to corroborate 
Lemma \ref{lemex4}'s claim that 
the axiom system $~B^\xi~+~$SelfRef$(B^\xi,d)$
must be
   consistent. (Thus,  
 $~B^\xi~+~$SelfRef$(B^\xi,d)$
must be  consistent 
because otherwise
a proof $~q~$ of  $~0=1~$ would exist and have its
 $~\mbox{ Log(}\,q \, ) ~ \geq ~ \zzz B^\xi+\Psi ~)~+~2~$ inequality
contradict \ep{piss3}.)  $~~\Box$

\medskip

%\cvs
%\cvl

{\bf How  Lemma \ref{lemex4}
May Be Interpreted :}
We 
%firstly 
remind the reader that many
(but not all)
generic configurations will
fail to satisfy 
Lemma \ref{lemex4}'s 
tightness hypothesis. This
is because any configuration
satisfying this hypothesis  represents one of those unusual
boundary-case
exceptions to the Second Incompleteness Theorem
that are feasible.

Lemma \ref{lemex4}
was intended to capture
the simplest variant of a self-justifying
phenomena (that employs 
Definition \ref{chg}'s machinery).
Its proof was
 kept
informal
because 
more sophisticated
 self-justifying formalisms will be explored
in Sections \ref{3uuuu3} and \ref{sect64}.
They will apply to
 four different types
of generic configurations, each of whose base axiom systems
$~B^\xi~$  
can be made capable of proving
all
Peano Arithmetic's
$~\Pi_1^\xi~$ theorems --- 
in a context where these  systems use a
broader variant
of {\it ``I am consistent''} axiom-statement
than does
Lemma \ref{lemex4}'s
``SelfRef$(B^\xi,d)$''
sentence.

\section{The First Two Meta-Theorems about Self-Justification}

\label{3uuuu3}

The core theorems in this section will employ
the following notation:
\bee
\item
The symbol $~\theta$ will 
denote 
any
recursively enumerable (r.e.)
set of
$\Pi_1^\xi$ sentences, 
henceforth called a {\bf R-View}.
(An
R-View
does not need to be valid under the
Standard-M model. It 
only needs to be r.e.)

\medskip

\item
 {\bf RE-Class$(\xi)$} 
shall denote the
set of all possible ``R-Views''
$~\theta~$ that 
can be built out of $~\xi \,$'s language
$~L^\xi~$. (This 
 permits
 both valid and invalid
R-Views  to appear in
RE-Class$(\xi)$. 
We choose this unrestricted definition because
no recursive 
decision
procedure can
identify
all  the true  $\Pi_1^\xi$ sentences
in the Standard-M model.)
\ene

\nop

\bxbxdr
\label{astab}
Let   $~\xi~$
denote 
the 5-tuple
  \gggen 
representing
one of
Definition \ref{def3.3}'s 
generic configurations,
and let RE-Class$(\xi)$
and its R-Views ``$\theta$''
be defined as in the previous paragraph.
Then
$ \,\xi \, $ is called
{\bf A-Stable\sss } iff 
each
$ \, \theta \, \in \, $RE-Class$(\xi)$
satisfies the 
following 
invariant:

\begin{description}
\item[  *$~$   ]
If $ \, \Uxp \, $ is a
$\Pi_1^\xi$ theorem 
of
axiom system $~\theta \cup B^\xi~$
via a proof $~p~$ 
whose length satisfies
 Log$(p)~\leq \zzz \theta)~+1 $
then  $ \, \Uxp \, $ 
will 
satisfy  Good$\{~ \, \tftt  \zzz \theta) \, ~\}~$.
\end{description}
\eedd

\begin{remark}
\label{re3-1}
\rm
The invariant
 $ \, * \,$ 
states
short proofs
(with lengths 
 $ \leq \,  \zzz \theta)~+1~~) $
  will 
produce at least 
{\it partially useful} deductions,
in that
their $\Pi_1^\xi$
theorems 
will always
 satisfy
  Good$\{ \,  \, \tftt  \zzz \theta) \,  \, \} \, $,
{\it irregardless of whether or not}
$\theta$'s axioms are {\it technically} true.
(This makes the study of  A-stability
very
interesting
unto itself,
 apart from its
applications
        in the
current article).
\end{remark}

\phx{ppp2} will show
the presence of
 A-stability,
alone,
is sufficient
for 
constructing 
self-justifying systems. 
This  will imply 
every A-stable configuration
must contain 
some
embedded weakness
(as every evasion
of the Second Incompleteness Theorem 
always does).
On the other hand,
Appendix F will explain how
A-stability 
and its 
``EA-stable'' cousin 
(defined later)
are both
epistemologically interesting.
Thus, A-stability has 
redeeming features.

\medskip

\bxbxdr
\label{estab}
A Generic Configuration 
$\,\xi \,$
 will be called
 {\bf E-Stable}  iff 
all of the
$\theta  \in \,$RE-Class$(\xi)~$
satisfy
  $~**~$.  $~$(This 
construct is the counterpart
for $\Sigma_1^\xi$ sentences of
the Item  $ *$ in
Definition \ref{astab}.)

\begin{description}
\item[ **   ]
If $ \, \Uxp \, $ is a
$\Sigma_1^\xi$ theorem 
 derived
from the
axiom system $~\theta \cup B^\xi~$
via a proof $~p~$ 
whose length satisfies
 Log$(p)~\leq \zzz \theta)~+1 $
$~$then  $ \, \Uxp \, $ 
will 
automatically satisfy 
Good$\{ ~ \frac{1}{2}~\lfloor ~$Log$(p)~\rfloor~-~1     ~\}~$.
{ (This invariant
{\it further
implies} 
\footnote{ \g55 
This point is easy to confirm when one
remembers that $\Sigma_1^\xi$ sentences $~\Upsilon~$
have the property that $~A<B~$ implies 
 Scope$_E$($\Uxp,A$)
is stronger than  Scope$_E$($\Uxp,B$).
It is then obvious
that the
 Good$\{ ~ \frac{1}{2}~\lfloor ~$Log$(p)~\rfloor~-~1     ~\}~$
criteria implies the validity of
 Good$\{~ \, \tftt  \zzz \theta) \, ~\}~$
for $\, \Sigma_1^\xi \,$ sentences
because the invariant $\,** \, $ presumes 
its  $\, \Sigma_1^\xi \,$ theorems have proofs $\, p \,$
satisfying 
`` Log$(p)~\leq \zzz \theta)~+1 $ ''. }  
that  $ \Uxp  $  
will also
satisfy 
 the Good$\{ \, \tftt  \zzz \theta) \,~\} $
criteria.)} 
\end{description}
\eedd

\smallskip

\begin{remark}
\label{re3-2}
\rm
The invariants $ \, * \, $ and $ \, ** \,  \,  $ 
are partially analogous 
to each other
because both imply that
if $ \, p \, $ is a proof short enough to
satisfy  Log$(p) \, \leq \zzz \theta) \, +1 \,  $ 
then their resulting theorem will satisfy
Good$\{ \,  \, \tftt  \zzz \theta) \,  \, \} \, $.
However, there is 
a
% an important 
distinction between 
Definitions
\ref{astab} and \ref{estab},
as well.
This is because the
prior section's
 Footnote \ref{fgood}
observed that $ \, \Sigma_1^\xi \, $ sentences are stronger
when they meet a 
Good$( \, N \, )$ 
rather than a
Good$( \, \infty \, )$ threshold, $ \, ${\it while the reverse is true for 
  $ \, \Pi_1^\xi \, $ sentences.}
Thus $ \, ** \,$'s short proofs
$ \, p \, $  
(satisfying  Log$(p) \,\leq \zzz \theta)~+1~ $)
will have the
special
 property that
their theorems  $ \, \Uxp~$ 
will satisfy a 
``Good$\{ ~ \frac{1}{2}~\lfloor ~$Log$(p)~\rfloor~-~1     ~\}~$''
constraint 
that is actually stronger than their formal
$~\Sigma_1^\xi~$ 
statements.
\end{remark}

\smallskip

The \appD will provide four examples
of generic configurations that are either E-stable
or A-stable (or often 
% very often 
both).
Its most prominent example will be a 
configuration
$~\xi^*~$ 
that uses  semantic tableaux deduction and
recognizes addition as a total function.
Theorem D.\ref{D.4-theorx} will 
imply such systems can be made
% both 
self-justifying
and 
able to
 prove 
% all
 Peano Arithmetic's
 $\Pi_1^*$ theorems.

\njp
%\smallskip

Three more  definitions
are needed 
to help introduce 
our first theorem

\smallskip

\bxbxdr
\label{eastab}
A generic configuration $~\xi~$
will be called
{\bf EA-stable} iff it 
is both 
E-stable
and A-stable.
(It will thus
satisfy both 
$~*~$ and $~**~$.)
\eedd

\smallskip

Our next definition is related to the fact
 that  
many definitions of consistency
are logically 
 equivalent 
from the perspective
of strong enough logics,
but they 
are {\it  often not provably equivalent}
from the perspectives of 
weak logics.

\medskip

\bxbxd
\label{lev}
\rm
Let $   \xi   $
denote the generic configuration
  \gggcp, 
and $   \alpha   $ be an axiom system satisfying
$   \alpha  \supseteq  B^\xi   $.
Then
$   \alpha   $ is called
{\bf Level$(k^\xi)$ 
Consistent} when there
exists  no proofs from $\alpha$ via
$ d   $'s deduction method  
of both a $\Pi^\xi_k$ sentence and 
of  the
$\Sigma^\xi_k$ sentence that is its  negation
\eedd

\medskip

Most of this article will focus on 
self-justifying systems that recognize their
Level$(k^\xi)$ consistency 
when 
 $k$ 
equals 0 or 1.
Our next definition will
be applied mostly to these two cases.

\smallskip

\bxbxdr
\label{xd+1x7}
Given any  
$k \geq 0$, 
a
generic configuration
  \gggen     
and an
axiom system
 $\beta  \supset  B^\xi$,
the symbol
SelfCons$^k(\beta,d)$
will
denote a 
self-referencing
$\Pi_1^\xi$ sentence declaring
$\,  \beta   +  $SelfCons$^k(\beta,d)      $'s
formal
Level(k$^\xi)$ consistency,
as is 
illustrated
below
by  
the 
statement $ ~ +~ $.  (An  encoding for
SelfCons$^k(\beta,d)$ 
will be provided
by Appendix A.)
\begin{quote}
$+~~~$ 
There exists  no two proofs 
(using deduction method $   d   $)
of {\it both}
some $\Pi_k^\xi$ sentence and
of the 
$\Sigma_k^\xi$ sentence, 
$\,$that 
represents its negation,
$\,$from the  union of
the axiom system $~\beta~ $
with $~this~$ added
sentence  ``SelfCons$^k(\beta,d) \,$''
{\it (looking at itself).}
\end{quote}
\eedd

\begin{remark}
\label{re3-3}
\rm
We will focus on
Definition \ref{xd+1x7}'s 
SelfCons$^k(\beta,d) \,$ axiom 
mostly in
 the settings
where $\,k=0$ or 1.
This is because 
SelfCons$^k(\beta,d) \,$ 
will 
typically 
be too strong for it to generate boundary-case exceptions
to the Second Incompleteness Theorem when
 $~k \geq 2~$.
It turns out that even when
 $~k=1~,~$  Definition \ref{xd+1x7}'s
SelfCons$^k(\beta,d) \,$ statement
will be significantly stronger than the
axiomatic
declaration $\, \bullet \,$
used by 
Section
\ref{secc1}'s
 SelfRef$(\beta,d)$
axiom.
This is 
because 
SelfCons$^1(\beta,d) \,$ asserts 
 non-existence of 
simultaneous proofs
for a  $~\Pi_1^\xi~$ sentence and its negation,
$\,$while  SelfRef$(\beta,d)$
establishes
merely
the
 non-existence of a proof of $0=1$.
\end{remark}

\smallskip

\begin{theorem}
\label{ppp1}
Let $~\xi~$
denote a
generic configuration
  \gggen   
that is 
EA-stable.  Then 
the corresponding  axiom system of 
$B^\xi+$SelfCons$^1(B^\xi,d)$
must satisfy Section
\ref{secc1}'s  definition of
self-justification.
\end{theorem} 

\parskip 3pt

\smallskip

{\bf Proof.}
The justification 
of \phx{ppp1} will
be 
a more elaborate  version
 of 
\lem{lemex4}'s
mini-proof.
It will
replace
Definition \ref{deftight}'s
Tightness constraint with
an EA-stability requirement.
It will also
 replace
SelfRef$(\beta,d)$'s 
``I am consistent''
axiom
with a 
{\it stronger}
SelfCons$^1(\beta,d)$
statement.

\cvb
\parskip 5pt

Our proof will focus on
showing 
  $~B^\xi+$SelfCons$^1(B^\xi,d)~$ 
is consistent (and thus 
% it
satisfies 
the subtle
Part-ii component
of Section \ref{secc1}'s
definition of
``Self-Justification'').
It will be awkward 
during our discussion
to 
write repeatedly  the
expression ``$~B^\xi+$SelfCons$^1(B^\xi,d)~$''.
Therefore,
 ``$S$'' will be the abbreviated
name for this axiom system.
We will also employ the following notation:
\bee
\itemsep +3pt
\item
$~\mbox{Prf}_S(t,p)~$ will denote 
that $~p~$ 
is
 a proof of the theorem $~t~$ from
the above mentioned formalism ``$~S~$''
(using 
$~\xi \,$'s
deduction method of $d$ ).
\item 
$~\mbox{Neg}^1(x,y)~$ will denote
that $~x~$  is
the G\"{o}del encoding of a
 $\Pi_1^\xi$ sentence and
that $~y~$ is a $\Sigma_1^\xi$ 
sentence which represents
$~x \,$'s formal negation.
\ene 
Appendix A explains how to combine the theory of LinH
functions \cite{HP91,Kr95,Wr78} with \cite{ww1}'s
fixed point methods to provide 
$\mbox{Prf}_S(t,p)$ 
and 
$\mbox{Neg}^1(x,y)$ with
$\Delta_0^\xi$ encodings.
(A reader can
 omit examining Appendix A,
if he 
just
accepts this fact.)
Thus, \eq{pencode} 
can
be 
viewed \footnote{ \s55 
Our notation convention has the 
abbreviated
formula of
``$~\mbox{Prf}_S(t,p)~$''
in \ep{pencode} 
corresponding 
to
 Appendix A's 
``$~\mbox{SubstPrf} \, _{\beta}^d  \, 
( \bar{n}    ,    t    ,    p   ) ~$'' formula,
in a context where $~n~$ specifies the G\"{o}del number of
the expression
$~\Gamma^k(g)~$  
in
Appendix A's 
 \ep{encode12} and where
the superscript $~k~$ 
within this
 formula
``$~\Gamma^k(g)~$''
is set equal to 1.  
This 
implies that sentence
\eq{pencode} does
assert the Level$(1^\xi)$ 
consistency
of S.  \fend }
in this context
as being   SelfCons$^1(B^\xi,d)\,$'s 
formalized
$\Pi_1^\xi$ 
statement,
declaring 
the Level$(1^\xi)$ 
consistency
of S:
\begin{equation}
\label{pencode}
\forall   \, x \, \forall   \, y \, \forall   \, p \, \forall   \, q 
~~~  \neg  ~~
 \{ \,\mbox{Neg}^1(x,y) ~ \wedge~ 
\mbox{Prf}_S  \, (     x    ,    p   ) ~
\wedge  ~
\mbox{Prf}_S \,    (       y    ,    q   )  \, \}
\end{equation}

Let   $ \Phi $   denote 
\eq{pencode}'s
 sentence. 
We will
use it to
prove
\phx{ppp1}'s claim that  $ S $
is consistent. Our proof 
will be a proof by contradiction.
It will 
begin with
the 
assumption that $ S $ is inconsistent.
This
implies 
$\Phi$ is
false under the Standard-M model.
Hence, 
Definition \ref{chg}
implies
\begin{equation}
\label{punk}
 \zzz \Phi~) 
   ~~ < ~~ \infty
\end{equation}

\noindent
\ep{punk} 
thus
indicates there  exists
$(\bar{p},\bar{q},\bar{x},\bar{y})$
satisfying \eq{pencodepunk}
(because such a tuple
will be a counter-example to
\eq{pencode}'s assertion).  
The particular
$(p,q,x,y)$
 satisfying \eq{pencodepunk}
with  minimum value for 
$ \mbox{Log}\{~\mbox{Max}[p,q,x,y] ~\}$ 
can 
then
be easily shown$\,$
\footnote{ \sm55 \label{footp1} 
Let $~L~$ denote the 
minimum value for 
$ \mbox{Log}\{~\mbox{Max}[p,q,x,y] ~\}$  for a tuple
$(p,q,x,y)$
 satisfying \ep{pencodepunk}. Then by definition,
$~\Phi~$ satisfies Good$(L-1)$ but not
Good$(L)$.
From Definition \ref{chg}, this 
establishes
 the validity of
\ep{punkless}  (because 
the minimal
$(\bar{p},\bar{q},\bar{x},\bar{y})$
satisfying sentence \eq{pencodepunk}
has
$ \mbox{Log}\{\mbox{Max}[  \bar{p}  ,  \bar{q}  ,  \bar{x}  ,  \bar{y}  ] \}
 \, = \, L.~$) }
  $\,$to also
satisfy
\ep{punkless}.
\begin{equation}
\label{pencodepunk}
  \,\mbox{Neg}^1(~ \bar{x}~, ~ \bar{y}~) ~ \wedge~ 
\mbox{Prf}_S  \, (     ~ \bar{x}~    ,    ~ \bar{p}~   ) ~
\wedge  ~
\mbox{Prf}_S \,    (       ~ \bar{y}~    ,    ~ \bar{q}~   )  
\end{equation}
\begin{equation}
\label{punkless}
 \mbox{Log}~\{~\mbox{Max}[~\bar{p}~,~\bar{q}~,~\bar{x}~,~\bar{y}~] ~~\}
   ~~~ =~~~\zzz \Phi~) ~~+~~1
\end{equation}

$~~~$ We will now use 
\eq{pencodepunk} and \eq{punkless} 
 to
bring our proof-by-contradiction to its conclusion.
Let $ \, \Upsilon \, $ denote the $\Pi_1^\xi$ sentence specified
by $ \, \bar{x} \, $. Then
$ \, \neg  \,  \Upsilon \, $ 
will
correspond to the $\Sigma_1^\xi$ sentence denoted
by $ \, \bar{y} \, $.
Also, 
\eq{pencodepunk} and \eq{punkless} imply that 
both
 $ \, \Upsilon \, $ and
$ \, \neg  \,  \Upsilon \, $
have  proofs 
such that the logarithms of
their G\"{o}del numbers are bounded
by
$\zzz \Phi \,)  \,+ \,1 \,$.
These facts
imply
 that
 $ \, \Upsilon \, $ and
$ \, \neg  \,  \Upsilon \, $
{\it both} 
satisfy 
 Good$\{ \,  \, \tftt  \zzz \Phi) \,  \, \} \, $
under our formalism.
(This is because
if we take $ \, \theta \, $ in Definitions
\ref{astab} and \ref{estab} to be simply $ \, \Phi \,$'s
1-sentence statement, $\,$then the invariants of
 $ \, *  \,$ and $ \, **  \,$ from these two definitions
both impose the
same
Good$\{ \,  \, \tftt  \zzz \Phi) \,  \, \} \, $
 constraint
on 
 $ \, \Upsilon \, $ and
$ \, \neg  \,  \Upsilon \, $.)

It is infeasible, 
however,
for a sentence and its
negation to both satisfy the same goodness
constraint.
 This completes \phx{ppp1}'s
proof-by-contradiction
because 
the initial
 assumption that 
$S$ was inconsistent has led to an
infeasible conclusion.
$~~\Box$

\bigskip

\cvl

\parskip 3pt

Our next definition will help formalize
a useful cousin of
\phx{ppp1}.

\smallskip

\bxbxdr
\label{ostab}
A Generic Configuration of 
$\xi $
will be called
 {\bf $~$0-Stable$~$} when every
particular
$\theta \, \in $ RE-Class$(\xi)$
satisfies the 
invariant of  $***$. 
(This invariant 
is strictly weaker than 
its counterparts   $*$  and $**$
in Definitions 
   \ref{astab} and \ref{estab}.) 
\begin{description}
\item[ ***   ]
If $ \, \Uxp \, $ is a
$\Delta_0^\xi$ theorem 
 derived
from the
axiom system $~\theta \cup B^\xi~$
via a proof $~p~$ 
whose length satisfies
 Log$(p)~\leq \zzz \theta)~+1~,~ $
then  $ \, \Uxp \, $
is 
true under the Standard-M model.
\end{description}
\eedd

\begin{theorem}
\label{ppp2}
If 
the configuration
$\xi$ is 
0-stable 
then
$B^\xi+$SelfCons$^0(B^\xi,d)$
is      a self-justifying 
formalism. {\rm (Appendix C 
shows this result  applies
also to
E-stable and  A-stable  configurations.)}
\end{theorem} 

\medskip

\phx{ppp2}'s proof
is similar
to  \phx{ppp1}'s
 proof. 
The difference between these two
propositions is that
\phx{ppp2}
has reduced
       SelfCons's
 superscript  
 from 1 to 0, so that its
 hypothesis
can encompass a 
theoretically
 broader 
set of applications.
(The Appendix C summarizes how
 \phx{ppp1}'s proof can be
easily modified to
also prove  \phx{ppp2}.)

\medskip

\begin{remark}
\label{newrem}
{\rm
Theorems 
\ref{ppp1} and \ref{ppp2}
should clarify 
the  nature of 
\cite{ww93}--\cite{ww9}'s
formalisms.
This is because proofs-by-contradictions 
are notorious 
in the mathematical literature
for being
confusing.
They should be simplified whenever possible.
This has been done
mainly through
\phx{ppp1}'s
short proof.
(It 
 applies to
three of Appendix D's
four examples of 
generic 
configurations,
and
\phx{ppp2} applies to 
Appendix D's
fourth example.)
Furthermore, Section \ref{sect64} 
will
show how more elaborate
self-justification systems can verify all
Peano Arithmetic's $\Pi_1^\xi$ theorems.}
\end{remark}

\section{Four Further Meta-Theorems}
\label{sect64}

We need one preliminary lemma
before  exploring 
how strong self-viewing
                logics may become
before 
they   cross
  the inevitable
boundary between 
self-justification
and 
 inconsistency,
$\,$implied by 
G\"{o}del's
Theorem.

\begin{lemma}
\label{llpp3}
Let  $  \xi  $ 
denote a generic configuration  
   \gggcp  
$\,$, and $  \mheta  $ denote an r.e. set of 
$\Pi_1^\xi$ sentences,
each of which holds true in the Standard-M model.
Let $  \mxi  $ denote a 
5-tuple that differs from 
 $  \xi  $  
in that its base axiom system
is  $  B^\xi \, \cup \, \mheta  $
(rather than
$  B^\xi  $).
These conditions imply that
$  \mxi  $ is a
generic configuration,
and 
it
will satisfy
the following
four invariants:
\bed
\item[     i   ]
If $~\xi~$ is  0-stable\sss 
then 
$~\mxi~$  will also be
0-stable\sss.
\item[      ii    ]
If $~\xi~$ is  A-stable\sss 
then 
$~\mxi~$ will also be
A-stable\sss.
\item[     iii   ]
If $~\xi~$ is  E-stable\sss 
then 
$~\mxi~$  will also be
E-stable\sss.
\item[     iv   ]
If $~\xi~$ is  EA-stable\sss 
then 
$~\mxi~$  will also be
EA-stable\sss.
\ennd
\end{lemma}

Lemma  \ref{llpp3}'s proof
is fairly straightforward. It
has been placed in
the Appendix B. 
This section will
use  Lemma  \ref{llpp3}
to prove
four
meta-theorems
 that are consequences of
its
formalism.

\smallskip

\begin{definition}
\label{dap4-1}
\rm
Let  $\xi$  again
denote
a generic configuration  
  \gggcp,  
and
$\theta$ denote some r.e. set of 
$\Pi_1^\xi$ sentences
(which are not
required to be  true  under the Standard-M model).
For the cases where $\,k \,$ is either
0 or 1,
the symbol 
$G^\xi_k(\, \theta \,)$
will
denote the following axiom system:
\vspace*{- 0.2 em}
\begin{equation}
\label{4gedef}
G^\xi_k(~ \theta ~)~~= ~~
\theta~\cup~B^\xi~\cup~\mbox{SelfCons}^k\{~[~\theta \,\cup \,B^\xi~]~,d~\}
\vspace*{- 0.2 em}
\end{equation}
Also when
$ \, k \, = \, 0 \,$
or 1,   the function 
$ \, G^\xi_k ~ $
(which maps $ \, \theta \, $ onto 
$G^\xi_k(\, \theta \,)~~~)$ is 
called 
{\bf Consistency Preserving} iff 
$ \, G^\xi_k(\, \theta \,) \, $is 
assured to be consistent whenever 
all the sentences in $~\theta ~$
are  
true  under the Standard-M model.
\end{definition}

\smallskip

We emphasize 
 consistency preservation is 
unusual in 
logic. This is because
$ \, G^\xi_k(\, \theta \,) \, $
comes from adding a self-justifying axiom to an
initially consistent
formalism $~B^\xi+\theta~$, and the 
Second Incompleteness Theorem
demonstrates that sufficiently powerful
formalisms are simply incompatible with
such an axiom.
 However, 
there will   be
four
specialized 
paradigms, defined in Appendix D,
that are exceptions to this rule.
They will be related to our next result:

\smallskip

\begin{theorem}
\label{pqq3} 
The function 
$~G^\xi_1~$ 
shall satisfy Definition \ref{dap4-1}'s
consistency preservation
property
when $~\xi~$ is
EA-stable. Likewise 
the function 
$~G^\xi_0~$ 
will be consistency preserving when $~\xi~$ is
any one of A-stable, E-stable or 0-stable.
{\rm 
(Thus 
in each case,
$~G^\xi_k(\, \theta \,)~$
will  
be consistent
when
all the sentences in $~\theta ~$
are  
true in the Standard-M model.)}
\end{theorem}

\medskip

{\bf Proof}.
It will be
 convenient
for our proof
 to use a
 dot-style
notation,
 analogous to
Lemma \ref{llpp3}'s
terminology. 
 Thus, 
\bee
\item
 $~\mheta~$ denotes any r.e. set of 
$\Pi_1^\xi$ 
sentences {\it that are 
each
true} under the Standard-M model.
\item
 $\mxi$
is the tuple
$ \, ( \, L^\xi  \, ,  \, \Delta_0^\xi  \,  ,
 \,  B^\xi \cup \mheta  \, ,  \, d   \, ,  \,  G  \,  )$.
(It
differs from 
 $\xi$
by replacing
 $\xi\,$'s  base axiom system
of $\,B^\xi \,$ with $\,B^\xi  \cup  \mheta$.)
\ene
 Part-iv of Lemma \ref{llpp3}
indicates
the EA-stability of 
$\xi$ 
 implies
the EA-stability of 
$  \mxi $.
 Moreover  \phx{ppp1} 
applies to all EA-stable
configurations, 
including  
$\mxi $.
Thus,
$G^\xi_1(\, \mheta \,)$ 
is  consistent
because $\xi$ is EA-stable
and all of $~\mheta \,$'s sentences are true in
the Standard-M model.
 This
proves
\phx{pqq3}'s 
first 
claim.

An 
almost
identical proof,
where 
\phx{ppp2} 
simply
 replaces \phx{ppp1}
as the central self-justifying
engine,
will corroborate 
\phx{pqq3}'s second claim .
Thus,
$~G^\xi_0~$ 
is consistency preserving when $~\xi~$ is
one of A-stable, E-stable or 0-stable. $~~~\Box$

\medskip

\begin{remark}
\label{re4-1}
\rm
Most
 generic configurations
$~\xi~$ 
will not satisfy \phx{pqq3}'s hypothesis
because
$~G^\xi_k(\theta)~$ will typically be
inconsistent,
irregardless 
\footnote{\sm55 Conventional generic configurations $\xi$
will satisfy the Hilbert-Bernays derivability conditions
\cite{HB39,HP91}.
Their
$~G^\xi_k(\theta)~$ will thus
be automatically inconsistent
because of
a G\"{o}del-like diagonalization argument.}
of whether or not
all of $~\theta\,$'s axioms are valid
under the Standard-M model.
The significance of \phx{pqq3} 
is that it shows that some
outlying
 exceptions to this
general rule do prevail when
$~\xi~$ satisfies one of 
\phx{pqq3}'s 
four stability conditions.
These
exceptions
are related to the
3-page abbreviated philosophical discussion
that will appear later in Appendix F.
They  
will 
assure 
that
Theorem \ref{pqq3}'s 
formalisms
of
 $G^\xi_1(\theta)$ 
and $G^\xi_0(\theta)$
are {\it always
consistent}
whenever 
$~\theta\,$'s axioms are valid
in the Standard-M model.
\end{remark}    

\smallskip

Our next definition will enable
self-justifying formalisms
to prove
the
$\Pi_1^\xi$ 
theorems 
of  any
consistent 
r.e. axiom system
that uses  $~L^\xi \,$'s language.

\smallskip

\begin{definition}
\label{dap4-2}
\rm
Let  $\,\xi\,$  denote
a generic configuration  
  \gggcp.   
Let $   \nal{B}  $ 
denote any recursive axiom system whose language
is an extension of 
$\, L^\xi \,$.
For an arbitrary 
deduction method $ \,  \nal{D} \,  $
(which may be 
possibly
different from $\,\xi \,$'s deduction method
$ \,d \,$ ),
 let
$\mbox{Prf}_{ \nal{B}}^D \,( \, \lceil \, \Psi \, \rceil \,  , \,  q \, )\,  $
denote a $\Delta_0^\xi$ formula indicating 
$  q  $ is a proof of the theorem $  \Psi  $ from 
axiom system
$ \,  \nal{B}\,$, 
using
deduction method
 $ \,  \nal{D} \,  $.
Then the  {\bf Group-2 Schema} for
$  (   \nal{B}  ,  D  )  $ is
defined as an infinite set of axioms that includes
one instance of \eq{group2}'s axiom 
for each 
$\Pi_1^\xi$ sentence $\, \Psi$. 
\beq
\label{group2}
\forall ~q~~~ \{~
\mbox{Prf}_{ \nal{B}}^D \,( \,  \lceil \, \Psi \, \rceil \,  , \,  q \, )
 ~~~\longrightarrow ~~~ \Psi ~~ \}
\enq
\end{definition}
{\bf  $~~~~$Comment about this Notation:}
The Definition \ref{dap4-2}
had called 
\eq{group2} a 
 ``Group-2 Schema''  so as to keep our terminology consistent
with
\cite{ww93,ww1,ww5,wwapal}'s notation.

\bigskip

%% Several of our articles
%% \cite{ww93,ww1,ww5,wwapal}
%% have  used the term
%%  ``Group-2 Schema'' to refer to
%% Definition \ref{dap4-2}'s
%% formalism. The same notation is
%% thus 
%% used here.

\smallskip

\begin{theorem}
\label{pqq4} 
Let  $~\xi~$ denote any
arbitrary generic configuration,
and 
$  (   \nal{B}  ,  D  )  $ 
denote any pair consisting of
an axiom system and a deduction method
(which, once  
again, are allowed to be 
different from
 $~\xi \,$'s deduction method
and axiom system).
Then if all of $  (   \nal{B}  ,  D  )  $'s
$\Pi_1^\xi$ theorems are 
true in the Standard-M model,
the following two invariants will hold:
\bed
\topsep -4pt
\itemsep -2pt
\item[   i ] 
If $~\xi~$ is EA-stable then
 there will exist an
r.e.$~$self-justifying system that can prove
all of 
$ (   \nal{B}  ,  D  )  $'s
$\Pi_1^\xi$ theorems 
and 
recognize its own Level($~1^\xi~)~$ consistency.
\item[  ii ] 
Likewise,
if
 $~\xi~$ is
one of A-stable, E-stable or 0-stable,
 there
will  exist an
r.e.$~$self-justifying system that can confirm
all of 
$ (   \nal{B}  ,  D  )  $'s
$\Pi_1^\xi$ theorems 
and 
which can
recognize its own Level($~0^\xi~)~$ consistency.
\ennd
\end{theorem}

\smallskip

{\bf Proof.}
\phx{pqq4} 
follows from
\phx{pqq3}.
Thus let
 $\theta$
denote
 the set of 
all
$\Pi_1^\xi$ sentences 
that are members of
Definition \ref{dap4-2}'s Group-2 schema.
Then every one
of $\theta\,$'s 
Group-2 axioms 
must be
true  under the Standard-M model
(because 
the hypothesis of  
\phx{pqq4} 
 indicated
all of $  (   \nal{B}  ,  D  )  $'s
$\Pi_1^\xi$ theorems are
true in this model).
Hence,  \phx{pqq3}  implies:
\bee
\item
$G^\xi_1(\, \theta \,)$
 is consistent when $\xi$ is EA-stable.
\item
$G^\xi_0(\, \theta \,)$ is consistent when $\xi$ is one of
 A-stable,  E-stable or  0-stable.
\ene
Since $G^\xi_0(\, \theta \,)$ and $G^\xi_1(\, \theta \,)$
are self-justifying systems that prove all of 
$ (   \nal{B}  ,  D  )  $'s
$\Pi_1^\xi$ theorems, 
Items 1 and 2
will substantiate
\phx{pqq4}'s
two claims. $\Box$

\bigskip

An
 awkward aspect of Definition \ref{dap4-2}'s
``Group-2'' schema is that 
it employs an
infinite number of instances
of \eq{group2}'s Group-2-like axiom sentences.
It turns out 
that this
Group-2 scheme
can be
compressed
 into a {\it single axiom sentence,}
if one is willing to 
settle for a
slightly 
diluted 
variant of
 $ (   \nal{B}  ,  D  )  $'s
$\Pi_1^\xi$ knowledge.
To 
formalize this  concept,
the following  notation shall be used:
\bee
\item
Check$^\xi(t)$ will denote a $\Delta_0^\xi$ formula
that checks to see whether $t$ represents the G\"{o}del
number of a $\Pi_1^\xi$ sentence. 
\item
Test$^\xi(t,x)$ 
will denote any $\Delta_0^\xi$ formula
where
\eq{testsim}'s invariant 
is 
true   under
the Standard-M model
for
every $  \Pi_1^\xi  $ sentence 
$   \Psi   $
simultaneously.
There are infinitely many 
different  $\Delta_0^\xi$ 
formulae that can
serve as 
Test$^\xi(t,x)$ predicates
satisfying 
this condition.
(Example \ref{new-exam} will illustrate
one such encoding of a
Test$^\xi(t,x)$ predicate.)
\beq
\label{testsim}
\Psi ~~~ \longleftrightarrow~~~ \forall ~x~~
\mbox{Test}^\xi(~\lceil~\Psi~\rceil~,~x~)
\enq
\ene
The expression \eq{globsim}
will be called 
a
{\bf Global Simulation Sentence}
for representing
$( \nal{B},D)$
via $\xi \,$.
Its $\mbox{Test}^\xi(t,x)$ clause essentially allows
$\xi$ to simulate the $\Pi_1^\xi$ knowledge of
$( \nal{B},D)$'s 
set of
theorems.
\beq
\label{globsim}
\forall ~t~~
\forall ~q~~
\forall ~x~~\{~~
[~~\mbox{Prf}_{ \nal{B}}^D \,(   t   ,   q   )~~ \wedge ~~
\mbox{Check}^\xi(t)~~]~~~
\longrightarrow ~~~ 
\mbox{Test}^\xi(t,x)~~~ \}
\enq

\bigskip

\begin{example}
\label{new-exam}
\rm
For any generic configuration $\,\xi~=\,$ \gggcp,
$\,$let NegPrf$^\xi(t,x)$ denote
a  $\Delta_0^\xi$ formula
specifying
that
$ \, t \, $ is a $\Pi_1^\xi$ sentence and
that
$ \, x \, $ is a proof under
 $\,d\,$'s deduction
method  of the 
 $\Sigma_1^\xi$ sentence that represents $ \, t\,$'s negation.
Also, let   Test$^\xi_0(t,x)$ 
be defined as follow:
\beq
\label{new1-exam}
\mbox{Test}^\xi_0(t,x)~~~~ =_{\mbox{def}}~~~~ \neg~ \mbox{NegPrf}^\xi(t,x)
\enq
For each $\Pi_1^\xi$ sentence $\Psi$,
%%  vvvv  
it is easy to verify 
\footnote{  \b55 Part-3 of Definition \ref{def3.3}
indicated that $~\xi\,$'s base axiom system is
``$~\Sigma_1^\xi~$ complete''. $~$(It is
thus able to prove
all $~\Sigma_1^\xi~$ sentences that are true in the true in
Standard-M model,
and it will likewise refute all 
$~\Pi_1^\xi~$ sentences that are false.)
The statement \eq{new1-exam}
then
immediately implies that
\eq{newtestsim}
must  be true
 under the Standard-M model
for
 every $~\Pi_1^\xi~$ sentence $~\Psi~$.}
that
the statement 
\eq{newtestsim}
is true
under
the 
Standard-M model.
\beq
\label{newtestsim}
\Psi ~~~ \longleftrightarrow~~~ \forall ~x~~
\mbox{Test}^\xi_0(~\lceil~\Psi~\rceil~,~x~)
\enq
\end{example}
\vspace*{- 1.3 em}
$~~~$ 
%Thus, 

The latter
% immediately
 implies 
that
\eq{new2-exam} 
is 
% an example of
a global simulation sentence
for  $  (   \nal{B}  ,  D  )  $.
\beq
\label{new2-exam} 
\forall ~t~~
\forall ~q~~
\forall ~x~~\{~~
[~~\mbox{Prf}_{ \nal{B}}^D \,(   t   ,   q   )~~ \wedge ~~
\mbox{Check}^\xi(t)~~]~~~
\longrightarrow ~~~ 
\mbox{Test}_0^\xi(t,x)~~~ \}
\enq
We emphasize that there are countably infinite different examples
of  $\mbox{Test}_i^\xi(t,x)~$ predicates that generate
global simulation sentences and
that statement \eq{new2-exam} 
illustrates only one such example.

%\bigskip

\begin{definition}
\label{gsim}
\rm
Let $( \nal{B},D)$ denote any ordered pair whose set
of $~\Pi_1^\xi$ theorems are true under the Standard-M model.
Let
Test$^\xi_1 \, , \, $ 
Test$^\xi_2 \, , \, $ Test$^\xi_3  \,~....~ $
denote the set of 
$\Delta_0^\xi$ formulae 
where
statement
\eq{testsim}
is true
under the Standard-M model
for every
$\Pi_1^\xi$ sentences
$\Psi$.
Then TestList$^\xi$ will denote
a
list of all  these 
Test$^\xi_i \,$ predicates.
Also for each
Test$^\xi_j $ 
formula
in  TestList$^\xi\,$,
$\,$the symbol
{\bf GlobSim$^D_{ \nal{B}} \,(\xi,j)$} 
will denote 
the special version  
of \eq{globsim}'s global simulation formalism
that employs
Test$^\xi_j \,  $'s machinery.
 \end{definition}

%\smallskip 

\begin{remark}
\label{re4-2}
\rm
A comparison between 
Definition \ref{dap4-2}'s
Group-2 schema with  \ref{gsim}'s
global simulation sentences
will
reveal
neither 
is 
strictly
 better than the other. 
Both         have their own separate advantages.
Thus,
the attractive aspect about Definition \ref{gsim}'s
GlobSim$^D_{ \nal{B}} \,(\xi,j)$  sentence is that it is a
finite-sized object that can simulate
the
infinite set of  axioms
associated with
Definition \ref{dap4-2}'s Group-2 schema. 
The accompanying
drawback
\footnote{ \f55  The
chief
 difficulty arises
essentially
        because
a  $\Pi_1^\xi$
theorem of
 $ (   \nal{B}  ,  D  )  $ 
may contain an arbitrarily long combination
of bounded universal and bounded existential quantifiers.
Thus,
some generic configurations will have 
base axiom systems $B^\xi$ 
that are so weak that
their combination with \eq{globsim}'s
Global Simulation Sentence 
is insufficient 
{\it to prove} the validity of \eq{testsim}'s
equivalence statement
{\it  for all}
 $\Pi_1^\xi$ sentences $\Psi$ simultaneously.
In particular, such   proofs
will often be infeasible
when  $\Psi \,$'s
sequence of bounded  quantifies
has a length greatly exceeding the length of
the G\"{o}del encoding for  \eq{globsim}'s
global simulation
statement.}
$\,$of 
a global simulation sentence
 is that the union of it
with the base-axiom system  $B^\xi$ 
will typically be inadequate  to prove every $\Pi_1^\xi$
sentence  that is a theorem of
 $ (   \nal{B}  ,  D  )  $.
$\,$Instead, in a context where
  $\Psi \,$ is a $\Pi_1^\xi$
theorem of
 $ (   \nal{B}  ,  D  )  $,
the sentence 
 GlobSim$^D_{ \nal{B}} \,(\xi,j)$ 
will usually
provide only
enough {\it  fragmented information} to prove
the statement \eq{quasisim} (which is equivalent to
$\Psi$ 
{\it under the Standard-M model).}
\beq
\label{quasisim}
 \forall ~x~~
\mbox{Test}_j^\xi(~\lceil~\Psi~\rceil~,~x~)
\vspace*{- 0.4 em}
\enq
While \eq{quasisim} may 
be insufficient
%not 
% provide
%enough information 
to prove $ \, \Psi \, $ from  
$ \, B^\xi \, $, it 
still 
(according to the 
statement \eq{testsim} )
has  the desired property of being
equivalent to $\Psi$
under the Standard-M model.
(This means that
the knowledge of \eq{quasisim}'s
truth is helpful, {\it even if it is unknown
from $ \, B^\xi\,$'s perspective} to be
equivalent to
$ \, \Psi \, $ . )
\end{remark}

%\smallskip

Our prior 
articles
\cite{ww93,ww1,ww5,wwapal} 
did not use
Definition \ref{gsim}'s
global simulation formalism.
They
employed,
instead,
Definition \ref{dap4-2}'s
Group-2 axiom schema.
\phx{pqq5} 
will be
   the 
 analog of  \phx{pqq4} 
for
global simulation.
It will be
 useful 
when
%in applications where 
one desires
to compress all the information held
by a Group-2 schema into a
single 
finite-sized object.

\smallskip

\begin{theorem}
\label{pqq5}
Let  $ \,\xi \,$ 
denote the generic configuration  
  \gggcp,   
 $\, \nal{B}\,$ denote a
 recursively enumerable
axiom system
and $ \,  \nal{D} \,$ 
denote any
deduction method
(which can be different than $\xi \,$'s deduction method    $~d~~)~$.
Suppose that all the 
$\Pi_1^\xi$
theorems generated by
 $ \, ( \,  \nal{B} \, , \, D \, ) \, $ are
true  under the Standard-M model.
Then the following  invariants
do hold:
\bed
\item[   i ]
If $ \, \xi \, $ is EA-stable 
then 
for each  $\, j \,$
there 
exists a 
{\bf finitized
extension}
$~\beta_j~$ of $~B^\xi$ that 
recognizes its 
Level$(1^\xi$)
self-consistency
 and which
contains
the sentence
GlobSim$^D_{ \nal{B}} \,(\xi,j)$. 
\item[  ii ] 
Likewise
for each $ j $,
if $\xi$ is 
  E-stable, A-stable  or 0-stable
then there 
exists a 
{finitized
extension}
$~\beta_j~$ of $~B^\xi$ that 
recognizes its 
Level$(0^\xi$)
self-consistency
 and which
contains
the sentence
GlobSim$^D_{ \nal{B}} \,(\xi,j)$. 
\ennd
\end{theorem}

%\medskip

{\bf Proof:}
Let
$\theta$ denote the 1-sentence R-View of 
``GlobSim$^D_{ \nal{B}} \,(\xi,j)$''
formalized by  Definition
\ref{gsim},
$~$and let 
$~\beta_j^1~$  and 
$~\beta_j^0~$ denote 
$ ~ \, B^\xi \, \cup \, \theta  \, +  \, $SelfCons$^1(
\,B^\xi \, \cup \, \theta  \, )  \,~ $
and
$  \, B^\xi \, \cup \, \theta  \, +  \, $SelfCons$^0(
\,B^\xi \, \cup \, \theta  \, )  \, $,
respectively.
These axiom systems
correspond 
to
the objects 
that Definition 
\ref{dap4-1} 
had called $G^\xi_1(\theta)$ 
and $G^\xi_0(\theta)$. 

\medskip

\phx{pqq5}'s  hypothesis 
indicates
%had specified that
all the 
$\Pi_1^\xi$
theorems of
 $( \nal{B},D)$  are true
 under the Standard-M model.
Thus,
it follows that $\, \theta \,$ is 
also
true in the Standard-M model.
Hence \phx{pqq3}
implies that 
$ \, \beta_j^1 \, = \, G^\xi_1(\theta)$ 
is a consistent  system 
satisfying \phx{pqq5}'s 
claim (i).
Likewise, \phx{pqq3}
implies 
$ \, \beta_j^0 \, = \, G^\xi_0(\theta)$ 
satisfies
\phx{pqq5}'s 
 second claim.
$~~\Box$

\bigskip

\begin{remark}
\label{re4-n}
\rm
Theorems  \ref{pqq4} 
and \ref{pqq5} raise a fascinating 
%open 
question:
{\it Is the trade-off between these 
% two 
formalisms needed?}
That is, 
%Thus,
can self-justifying systems use only a 
 a finite number of added axioms
beyond those lying in $~\xi\,$'s base 
%axiom 
system of $B^\xi$
and 
also
%{\it simultaneously} 
duplicate all 
 $( \nal{B},D)$'s $\Pi_1^\xi$ theorems in a pure 
sense (i.e. without simulation) ?
We will return to this topic in Appendix G.
\end{remark}

%%nov150ut 
%%nov150ut \begin{remark}
%%nov150ut \label{old-re4-n}
%%nov150ut \rm
%%nov150ut Theorems  \ref{pqq4} 
%%nov150ut and \ref{pqq5} 
%%nov150ut are difficult to compare
%%nov150ut  because
%%nov150ut neither result is uniformly better than the other.
%%nov150ut The advantage to Theorem  \ref{pqq4} 
%%nov150ut is that it duplicates all of
%%nov150ut  $( \nal{B},D)$'s $\Pi_1^\xi$ theorems in a pure sense,
%%nov150ut rather than using Theorem  \ref{pqq5}'s
%%nov150ut weaker simulated construct.
%%nov150ut The virtue of Theorem  \ref{pqq5}, 
%%nov150ut in contrast,
%%nov150ut  is that its self-justifying systems 
%%nov150ut include only
%%nov150ut Our suspicion is that Theorem  
%%nov150ut \ref{pqq5} is of greater  pragmatic interest,
%%nov150ut while
%%nov150ut \thx{pqq4} is the more mathematically surprising result.
%%nov150ut \end{remark}

%% \begin{remark}
%% \label{recc1}
%% \rm
%% This article deliberately
%% used the phrase ``boundary-case exception''
%% to characterize our limited  evasions
%% the Second Incompleteness Theorem. 
%% This
%% cautious terminology was employed
%% because there are two fundamental barriers
%% limiting the generality 
%% of all such
%% results:  

\bigskip
\bigskip
\bigskip
\bigskip

%{\bf $~~~~$ A BROADER PERSPECTIVE:}

% {\bf $~~~~$ A Yet Further Issue$~$:}

% {\bf $~~~~$ A Yet Further Issue$~$:}
{\bf $~~~~$ REFLECTION PARADIGMS$~$:}
$~$Our  
last
%final 
goal
is to show how
self-justifying systems 
%can 
support
%% 
%% will be to formalize
%% how Theorems \ref{ppp1},
%%  \ref{pqq3}, \ref{pqq4} 
%% and \ref{pqq5} support 
%% 
{\it unusually}
strong 
reflection principles.
% for their $\Pi_1^\xi$ theorems. 
Let $\mbox{Reflect}_{\alpha,d}(~\Psi~)~$ 
%thus
denote  \eq{reflect}'s statement 
when
%in a context where
$~\Psi~$ is a sentence with G\"{o}del number
 $~\lceil \, \Psi \, \rceil~$, and 
$~(\alpha,d)~$ denotes an axiom system
and deduction method.
\beq
\label{reflect}
\forall ~p~~~[~~ \mbox{Prf}_{\alpha,d}(~\lceil \, \Psi \, \rceil~,~p~)
 ~~~ \Rightarrow ~~~ \Psi~~]
\enq
L\"{o}b's Theorem
\cite{HP91,Lo55,So76b}
 implies that conventional
systems $\, (\alpha,d) \,$, 
possessing at least
 Peano Arithmetic's strength, are unable to
prove 
% the sentence 
 $\mbox{Reflect}_{\alpha,d}(~\Psi~)~$ except for in the degenerate cases
where they can prove $~\Psi$.

% itself.

%\cvt

\medskip

Moreover, it is easy to generalize
L\"{o}b's Theorem
(via say  \cite{ww1}'s
Theorem 7.2) 
so that
a wide class of formalisms $~\alpha~$, weaker than
%PA,
Peano Arithmetic,
are also
 unable to prove
% the validity of 
 $\mbox{Reflect}_{\alpha,d}(~\Psi~)~$ 
for
% essentially 
all $\Pi_1^\xi$
sentences $~\Psi~$ simultaneously. 

\medskip

The intuition
behind this generalization is quite simple.
Let
$~\mho~$
denote 
a $\Pi_1^\xi$ encoding for
the classic
 G\"{o}del sentence
declaring: 
{ \it  `` There is no proof of me 
%{\it this sentence}
from the axiom system $~\alpha~$ using $~d\,$'s deduction method''.}
Then \cite{ww1}'s
Theorem 7.2  
uses a very routine 
diagonalization argument
to show 
most formalisms
$~\alpha~$ will be
 inconsistent 
%(and thus useless)
if they  prove  $\mbox{Reflect}_{\alpha,d}(~\mho~)$'s
statement. 

\medskip

%% 
%% The intuition
%% behind this result is quite simple.
%% Let
%% $~\mho~$
%% denote 
%% a $\Pi_1^\xi$ encoding for
%% the classic
%%  G\"{o}del sentence
%% declaring: 
%% { \it  `` There is no proof of me 
%% %{\it this sentence}
%% from the axiom system $~\alpha~$ using $~d\,$'s deduction method''.}
%% Then 
%% Theorem 7.2  
%% used an essentially 
%% conventional
%% diagonalization argument
%% to show 
%% its formalisms
%% $~\alpha~$ would be 
%%  inconsistent (and thus useless)
%% if they proved  $\mbox{Reflect}_{\alpha,d}(~\mho~)$'s
%% statement. 

Our next theorem will show, surprisingly, that
the preceding limitation is much less 
stringent 
% severe
than 
it may
initially  appear to be.
This is
because 
%EA-stable 
Level$(~1^\xi~)$ self-justifying 
axiom systems
are capable of proving  
very close analogs to
\eq{reflect}'s
impermissible

\nop

\noindent
 reflection principle,
using a ``translational'' methodology.

\smallskip

Thus, let $~T~$  denote an
algorithm that maps a
% initial 
$~\Pi_1^\xi~$ sentence $~\Psi~$ onto a 
``translated''
sentence  $~\Psi^T~$ that is equivalent
to  $~\Psi~$ 
under the Standard-M model
{\it and which is also written in a  $~\Pi_1^\xi~$ 
format.} 
(See footnote\footnote{\label{imper} Part-3 of Definition
\ref{def3.3} indicated that generic configurations are
$\Sigma_1^\xi$ complete. 
%In this context, 
Our requirement that
 $~\Psi^T~$ 
must have a       $\Pi_1^\xi$ format
thus
causes 
$~T~$ 
to
gain  much added
meaning.
This is 
because 
the
axiom system $~B^\xi~$ will then  
automatically  disprove
$~\Psi^T~$ 
 whenever
it is false under the Standard-M model.
(Thus, $T$'s mapping of 
$~\Psi~$  onto $~\Psi^T~$ 
gains much  significance
when 
$~\Psi~$ and
$~\Psi^T~$ do rest on the same  $\Pi_1^\xi$ 
level of the arithmetic
hierarchy.)}
for why 
it is absolutely imperative that
{\bf both these} requirements be included in
$~T\,$'s definition.)
Also, let
$\mbox{Reflect}^T_{\alpha,d}(~\Psi~)~$ denote the
%  corresponding
translational
modification of \eq{reflect}'s reflection principle that
replaces $\Psi$ with $\Psi^T$. 
%in its
%modified {\bf ``Translation'' Reflection Principle.}
\beq
\label{T-reflect}
\forall ~p~~~[~~ \mbox{Prf}_{\alpha,d}(~\lceil \, \Psi \, \rceil~,~p~)
  ~~~ \Rightarrow ~~~ \Psi^T~~]
\enq

\begin{theorem}
\label{ppp6}
Let  $ \,\xi \,$ 
denote the EA-stable configuration of  
  \gggcp,   
and let
$~\alpha~=~   B^\xi \, +  \, $SelfCons$^1(\,B^\xi \,)$
denote  $ \,\xi \,$'s corresponding Level(1) self-justifying
axiom system. Then
there will exist 
a translation methodology $~T~$ where
 $~\alpha~$ can prove the validity
of \eq{T-reflect} for all its $\Pi_1^\xi$ sentences simultaneously.
%% 
%% {\rm (In other words, this means that
%% $~\alpha~$ will  know that
%% its 
%% proof of $~\Psi~$ automatically
%% implies a counterpart
%% of $~\Psi~$ is true in
%% the Standard-M model ! )}
%% 
\end{theorem}

% {\bf Abbreviated Proof:}

\bigskip

{\bf Proof:} Let us use
Example  \ref{new-exam}'s notation.
It observed that  $~\Psi\,$'s
$\Pi_1^\xi$ statement
was equivalent 
under the Standard-M model
to
 ``$ \,  \forall  \, x~ \mbox{Test}^\xi_0( \, \lceil \, \Psi \, \rceil \, , \, x \, )$''.
Thus, 
let us view $\,T \,$  as 
being
a mapping of
the first sentence onto the second. 

Our proof of Theorem \ref{ppp6} will next  use
the following observations:
\bee
\item The non-existence of a proof of $~ \neg ~ \Psi~$ from
$~    B^\xi \, +  \, $SelfCons$^1(\,B^\xi \,)$
trivially implies the non-existence of a proof of the same theorem
from $~B^\xi~$
 (because the latter axiom system is 
simply
a subset of
the former).
\item 
Moreover, Example  \ref{new-exam}'s notation treats
``$ \,  \forall  \, x~ \mbox{Test}^\xi_0( \, \lceil \, \Psi \, \rceil \, , \, x \, )$'' as being equivalent to the declaration that no proof of
$~\neg~\Psi~$ from $~B^\xi~$ exists. 
\ene
Hence, $~\alpha~$ can prove 
\eq{T-reflect}'s  statement by noting
that 
$~p\,$'s proof of 
a $\Pi_1^\xi$ sentence
$~\Psi~$ implies 
(via $\, \alpha \,$'s $~$SelfCons$^1~$  axiom) 
the non-existence
of a proof of $~\neg ~\Psi~$, which 
(via Items 1 and 2)
implies 
 ``$ \,  \forall  \, x~ \mbox{Test}^\xi_0( \, \lceil \, \Psi \, \rceil \, , \, x \, )$''
$~~\Box$.

\cvs

\bigskip
\bigskip

\begin{remark}
\rm \label{f88}:    
$~$\phx{ppp6} and the statement
\eq{T-reflect}'s 
Translational Reflection Principle
may possibly be useful devices in unraveling some of the
mystery that has enshrouded G\"{o}del's Second Incompleteness
Theorem, since its inception.
This is
partly
 because G\"{o}del was 
explicitly uncertain
about the 
generality of the
Second Incompleteness Theorem 
in his initial 
 1931 
seminal
paper \cite{Go31}
about this subject.
% to appreciate the preceding 
%oint.
His centennial paper
about Incompleteness
 thus included 
the following 
quite
poignant
caveat:
\begin{quote}
%\small
%\baselineskip = 1.0 \normalbaselineskip 
$~\bullet~~~$ : 
``It must be expressly noted that
Theorem XI
(i.e the Second Incompleteness Theorem) 
represents no contradiction of the formalistic
standpoint of Hilbert. For this standpoint
presupposes only the existence of a consistency
proof by finite means, and {\it there might
conceivably be finite proofs} which cannot
be stated in ... ''
%
%%%%%%%%%%%%%%%%%%%%%%%% 
%
\end{quote}
Some of the issues
that troubled 
G\"{o}del 
in the 
statement 
$\bullet$
can 
perhaps 
be 
partially
resolved if one
compares the reflection principles of sentences
\eq{reflect} 
and \eq{T-reflect}.
This is because 
\eq{reflect} is probably
unnecessary to explain how thinking beings can
appreciate their
 $~\Pi_1^\xi~$
 theorems 
%$~$---$~$ 
when 
\phx{ppp6}'s specialized
logics can,
instead,
 use the fact that
its $\Pi_1^\xi$ sentences  satisfy at least
\eq{T-reflect}'s
modified
 reflection principle. 
%% 
%% (with
%% $~\Psi^T~$ being required to be both 
%% $~\Pi_1^\xi~$ and equivalent to $~\Psi~$ in the
%% Standard-M model). 
%% 
%%
%
%In other words, 
Thus,
some of the mystery surrounding the
Second Incompleteness
Effect can be
clarified
 when one notices that 
%the validity of
\eq{T-reflect}'s translational reflection princible 
is a useful precept, that was
shown by \phx{ppp6} to be
technically
%essentially
 unrelated to G\"{o}del's observation
that no
reasonable formalism  
can prove
%proves 
\eq{reflect}'s 
{\it purist}
%a reflection 
principle
 for all 
 $~\Pi_1^\xi~$
sentences
simultaneously.
\end{remark}

%%% 
%%% 
%%% Some of the issues 
%%% G\"{o}del 
%%% had raised in statement $~\bullet~$
%%% can be
%%% partially
%%%  resolved when one
%%% compares the reflection principles of sentences
%%% \eq{reflect} 
%%% and \eq{T-reflect}.
%%% This is because 
%%% \eq{reflect}'s reflection principle is probably
%%% unnecessary to explain how thinking beings can
%%% appreciate their
%%%  $~\Pi_1^\xi~$
%%%  theorems $~$---$~$ 
%%% when some logics can,
%%% instead,
%%%  use the fact that
%%% every $~\Pi_1^\xi~$ statement  satisfies at least
%%% \eq{T-reflect}'s
%%% modified
%%%  reflection principle. 
%%% %% 
%%% %% (with
%%% %% $~\Psi^T~$ being required to be both 
%%% %% $~\Pi_1^\xi~$ and equivalent to $~\Psi~$ in the
%%% %% Standard-M model). 
%%% %% 
%%% %%
%%% In other words, some of the mystery surrounding the
%%% Second Incompleteness
%%% Effect can be
%%% resolved
%%%  when one notices that the validity of
%%% \eq{T-reflect}'s translational reflection princible 
%%% (for some {\it specialized} 
%%% logics)
%%% is unrelated to G\"{o}del's observation
%%% that no
%%% reasonable formalism can feasibly prove 
%%% \eq{reflect}'s 
%%% {\it purist form} of
%%% a
%%% reflection principle
%%%  for all 
%%%  $~\Pi_1^\xi~$
%%% sentences
%%% simultaneously.
%%% %% \end{remark}

%% 
%% ) $~$---$~$
%% {\it while at the same time} \phx{ppp6}
%% permits a 
%%  formalism
%% to,
%% nevertheless, 
%% support  \eq{T-reflect}'s modified reflection principle and
%% simultaneously
%%  prove all Peano Arithmetic's  $~\Pi_1^\xi~$
%% theorems.

\bigskip

\medskip

\begin{remark}
\label{remhappy}
\rm
Theorem \ref{ppp6} is
also
significant
because it explains how
its  specialized logics 
 can grapple with a 
$\, \Pi_1^\xi \, $ encoded
G\"{o}del sentence
$~\mho~$ 
which asserts {\it ``There is no proof of me''}.
This issue is challenging because
routine constructions, such as
\cite{ww1}'s Theorem 7.2,  
demonstrate that no natural logic
can 
verify
statement
  \eq{reflect}'s validity 
for all $~\Pi_1^\xi~$ 
sentences $~\Psi~$
(on account of the
well-known
syllogism
posed by 
$~\mho$'s  G\"{o}del sentence).
%% 
%% 
%% on account of the challenge posed by 
%% $~\mho\,$'s available  
%%  $~\Pi_1^\xi~$ encoding. 
%% 
%% 
Theorem \ref{ppp6} 
constructs,
%demonstrates,
however,
a 
%method to 
%serious
reply to this challenge. 
%% 
%%% 
%%%  that this 
%%% challenge
%%% %difficulty
%%% is not quite as disturbing
%%% as  might first 
%%% appear
%%% %be anticipated.
%%% 
This is
because 
%%%Theorem \ref{ppp6}'s 
its self-justifying systems do
%%% can
surprisingly
%{\it at least,}
prove,
without difficulty \footnote{Since 
$\Psi$ and $\Psi^T~$ are equivalent under
the Standard-M Model {\it but not also
equivalent}  
from the perspective of the 
%axiom 
system 
$~\alpha~=~   B^\xi \, +  \, $SelfCons$^1(\,B^\xi \,)$,
the conventional contradictions, produced by
% sentence 
\eq{reflect}'s reflection principle, disappear when it
is replaced by 
\eq{T-reflect}.
%'s
% theoretical
% alternative principle. 
(Thus, there is no danger that
$\alpha$ could
use \eq{T-reflect}'s reflection principle to 
prove 
an analog of $\,\mho \,$'s 
forbidden G\"{o}del sentence.) \fend
%the forbidden statement
%{\it ``There is no proof of me''} 
%using
%%
%% 
%% The 
%% ``difficulty'' that conventionally arises is that an
%% axiom system 
%% $~\alpha~$
%% becomes inconsistent when it proves
%% $\mbox{Reflect}_{\alpha,d}(~\mho~)~$'s sentence because this causes
%% it to prove $~\mho\,$'s sentence, which asserts {\it ``There is no
%% proof of me''.}  
%% This difficulty vanishes,
%% however, 
%%  when $~\alpha~$
%% proves the translational reflection sentence of 
%% $\mbox{Reflect}^T_{\alpha,d}(~\mho~)$.
%% This is because the latter replaces
%% $~\mho~$ in  \ep{reflect}
%% with \eq{T-reflect}'s sentence $~\mho^T~$.
%% This second sentence  is equivalent to $~\mho~$
%% under the Standard-M
%% model.
%% However,
%% $~\alpha~$ is unable to recognize
%% this equivalence (and  its
%% thus escapes the syllogism of
%% a G\"{o}del-like diagonalization argument).
%% 
%% 
}, the validity of
\eq{T-reflect}'s 
{\it translated modification}
of  \eq{reflect}'s
% formally 
unobtainable
$\Pi_1^\xi$ 
styled
variant of
a reflection principle.
\end{remark}

\bigskip

\begin{remark}
\label{new14}
\rm
We encourage the reader to examine the work of
Beklemishev, Kreisel-Takeuti
 and Verbrugge-Visser
\cite{Be95,Be97,Be3,KT74,VV94,Vi5}
to see
alternative
reflection principles
and 
their uses.
(The  constraint on proof-length,
by Verbrugge-Visser,
 is certainly
one 
alternative to
\thx{ppp6}'s machinery.
Likewise, 
%\cite{KT74}'s 
Kreisel-Takeuti's
Second Order
Logic {\it CFA}
reflection is another alternative,
although it will not \footnote{Kreisel-Takeuti indicate on
page 25 of \cite{KT74} that CFA's reflection principle for
a first order formula ``A'' infers the validity of
{\it only the relativized formula} ``$~A^N~$'' from $A$'s proof.
Also, their proof predicate is similarily relativized. 
Thus CFA's second-order logic
reflection principle, while fascinating, does
not generalize to first-order logic environments. }
generalize to first-order 
logics.) 
One 
complicating aspect of
%drawback to
our  \thx{ppp6}'s
%translational 
reflection
%principle
method
 is
%  quite 
that
Appendix E 
% formally 
proves it becomes
%%% unfortunately
% fully
 inoperable
when an axiom system is 
% simply
sufficiently conventional
to satisfy G\"{o}del's Second Incompleteness Theorem.
In essence, \thx{ppp6}'s translational reflection
principle is a
specialized
methodology,
% an attractive type of
%machinery,
%reflection methodology, 
% that is
intended for logics using
Definition \ref{xd+1x7}'s
Level(1) style
of self-justification. 
% ``SelfCons''
% self-justifying axiom.
\end{remark}

%%%% cccccccccc
%%%% Some readers may wish to examine the work of
%%%% Beklemishev, Verbrugge and Visser
%%%% \cite{Be95,Be97,Be3,VV94,Vi5}
%%%% to see
%%%% some  different 
%%%% perspectives about 
%%%% reflection principles.
%%%% (The Verbrugge-Visser mechanism,
%%%% with its  constraints  \cite{VV94,Vi5} on proof-length,
%%%%  is 
%%%% a
%%%% % an especially
%%%% fascinating alternative to
%%%% \thx{ppp6}'s machinery.) 
%%%% One 
%%%% complicating aspect of
%%%% %drawback to
%%%% our  \thx{ppp6}'s
%%%% %translational 
%%%% reflection
%%%% principle is
%%%% %  quite 
%%%% that
%%%% Appendix E 
%%%% % formally 
%%%% proves it becomes
%%%% %%% unfortunately
%%%% % fully
%%%%  inoperable
%%%% when an axiom system is 
%%%% % simply
%%%% sufficiently conventional
%%%% to satisfy G\"{o}del's Second Incompleteness Theorem.
%%%% Thus, self-justifying systems support a type
%%%% of 
%%%% %surprising 
%%%% reflection principle that is 
%%%% uniquely intended 
%%%%  for their particular
%%%% % special types of
%%%% applications.
%%%% %inapplicable in conventional settings.

\bigskip

%\smallskip

\begin{remark}
\label{recc1}
\rm
The preceding discussion clearly shows 
%that 
self-justifying logics are 
% very 
% quite 
tempting.
At the same time,
it is necessary 
to be 
very
cautious
because there are also two 
fundamental 
barriers
limiting 
such results:
\bed
\topsep -4pt
\itemsep +2pt
\item[   a  ] 
The
first is the
 Theorem 2.1 arising from the joint work of
Pudl\'{a}k, Solovay, Nelson and Wilkie-Paris \cite{Ne86,Pu85,So94,WP87}.
It
showed
no 
reasonable 
system recognizing successor
as a total function can verify its own Hilbert consistency.
%Moreover,
Also,
Willard \cite{ww2,ww7} established 
%% ww0
analogous results 
under semantic tableaux consistency 
for 
% most
systems recognizing multiplication
as a total function. 
%These results show that every effort to evade
Thus, each effort to evade
the  Second Incompleteness Theorem 
must
encounter 
%inevitable
 robust
barriers.
\item[   b  ]
A second issue is that Definition \ref{xd+1x7}'s
``SelfCons''
{\it ``I am consistent''} 
axiom
sentence is less than ideal because 
it causes axiom systems to produce
essentially
a
1-line proof of their
own 
 consistency.
Such an excessively compressed proof corresponds more
closely to an axiom system formulating
 an
{\it instinctive faith} in its own consistency
(rather than it  supporting 
a 
% more elaborate 
full-length proof-justification
of this fact).
\ennd
Part of the reason 
self-justifying systems are of
interest, despite these limitations, is that they illustrate
how some formalisms {\it are compatible} with at least an 
{\it instinctive faith} in their own self-consistency.
(This compatibility issue is non-trivial
because Item (a) 
implies
there are
many circumstances where a generalization of the Second Incompleteness
Theorem will make it infeasible for a formalism to satisfy
{\it both} 
Parts (i) and (ii) of Section 1's definition
of Self-Justification.)
Moreover, three of Appendix D's four sample self-justifying
configurations,
called $~\xi^*~$, $~\xi^{**}~$ and $~\xi^R~$,
 will be Type-A systems that
 recognize
addition as a total function. These configurations will
thus
possess
the following three 
significant
finitized 
 features:
\bee
\topsep -4pt
\itemsep +2pt
\item
They will be able to construct the entire infinite set of integers
{\it by finite means} because they recognize addition as a total function.
\item
For any r.e. logical
configuration $  (   \nal{B}  ,  D  )  $, 
it will be possible to develop a 1-sentence
 {\it finitized} 
extension for the base axiom systems of any of  
the configurations of  $~\xi^*~$, $~\xi^{**}~$ and $~\xi^R~,~$
which deploy
 \eq{globsim}'s
 Global Simulation Sentence
to simulate the 
 $~\Pi^\xi_1~$ knowledge of  $  (   \nal{B}  ,  D  )  $.
This means 
that some fully
finite-sized
extensions of the base-formalisms of 
  $~\xi^*~$, $~\xi^{**}~$ and $~\xi^R~$ will contain a
non-trivial amount of 
 $~\Pi^\xi_1~$ styled
knowledge,
since  $  (   \nal{B}  ,  D  )  $ can 
correspond to, say,
%represent for example
Peano Arithmetic. 
\item
The 
key point is
 that a 1-sentence extension of an
axiom system  containing  features (1) and (2)
can formalize 
how a
logic
 can possess
an instinctive faith in its own
consistency
via \thx{pqq5}'s explicitly
{\it finitized} structure.
(Moreover, \thx{ppp6}'s 
Translational  Reflection Principle
is applicable to 
Appendix D's
 generic configurations of
 $~\xi^*~$ and $~\xi^{**}~$.
It will 
thus imply that  their single 
{\it finitized}
 Level-1 self-justifying
axioms enable 
% each of
them to prove
an  {\it infinite number} of incarnations
of \eq{T-reflect}'s
translational 
reflection principle, where  each
$\Pi_1^*$ sentence  
$\Psi$ 
is mapped onto one such 
%well-definded unique 
unique instance.)
%%% 
%%% 
%%% Appendix G's
%%%  more 
%%% sophisticated
%%%  Theorems G.2 and G.3,
%%% %%% 
%%% %%% accompanied by 
%%% %%% Appendix F's
%%% %%% epistemological  
%%% %%% discussion,
%%% %%% in{re
%%% will
%%% further strengthens   
%%% these results.)
%%% 
%%% 
\ene
The contrast between Items 1-3's positive remarks
about ``finitized'' cogitation with
Items (a) and (b)'s 
opposing comments 
is 
obviously
formidable. 
It is clearly preferable
to
view these
positive results 
cautiously and treat them
as being no more than
boundary-case  exceptions  to
the Second Incompleteness Theorem. 
The essential reason why 
these exceptions are of
interest
is that G\"{o}del's famous centennial paper has
implicitly raised
the following 
%unresolved 
puzzling issue:
\begin{quote}
 \#  How is it that Human Beings
 manage 
to muster 
the physical 
drive
to think (and prove theorems) when the many
generalizations of 
G\"{o}del's Second 
Incompleteness Theorem assert 
conventional logics
lack knowledge of
their own consistency?  
\end{quote}
There will, of course, never be any perfect answer to the 
puzzle
posed by
$~\# ~$ 
because philosophical paradoxes and ironical dilemmas never yield
perfect answers. 
However, 
 part of an
imperfect
 answer to 
$~\# ~$ is that 
Items 1-3 reply to
 Challenges (a) and (b) by formalizing how
a thinking being can muster an
% at least
approximate partial
{\it instinctive faith} in its own self-consistency.
(Moreover, the
tight
 contrast 
between various 
generalizations of the Second
Incompleteness Theorem
\cite{Ad2,AZ1,BS76,BI95,HP91,HB39,Ko6,Lo55,Pu85,Sa11,So94,Sv7,WP87,ww2,wwlogos,ww7}
with the
self-justifying systems appearing
in Appendixes D and G
suggests 
that these
come close to
being
maximal forms
of
feasible   
results.)
\end{remark}

%\smallskip

Our remaining discussion 
will consist of four optional sections, called
 Appendixes D, E,  F and G,
% consists of three 
% optional
% sections,
which
can be skimmed,
omitted or  examined
in 
any
% whichever 
order 
 the reader 
%so 
prefers.
A 
%brief 
summary of 
their contents 
is given 
%provided
below:
\bed
\itemsep +4pt
\item[   I  ] 
The Appendix D provides four examples of
generic configurations
 that
utilize
Theorems \ref{ppp1}, \ref{ppp2},
\ref{pqq3}, \ref{pqq4} 
and \ref{pqq5}.
Its most prominent
examples involve 
\ep{totdefxa}'s
Type-A axiom systems
where the deduction method is
either semantic tableaux
or a modified version of tableaux that permits
a modus ponens rule for $\Pi_1^\xi$ and $\Sigma_1^\xi$ sentences.
\item[ II  ]
The Appendix E introduces a generalization of
the  Second Incompleteness
Theorem which shows that \thx{ppp6}'s
Translational Reflection Principle
applies {\it only to} 
self-justifying logics. (It is thus  fully inoperative
for conventional logics.
This may explain why 
\thx{ppp6}'s
self-justifying systems are an interesting topic.)
\item[ III  ]
The Appendix F differs from the rest of
this paper
by having a 
philosophical
slant.
It will 
offer a 3-page summary about why we suspect 
\thx{pqq3}'s 
self-justification
formalism and Remark \ref{recc1}'s
notion of ``instinctive faith'' are useful.
%and related to Theorem \ref{ppp}'s reflection principle.
\item[  IV  ]
The Appendix G 
introduces 
a ``Braced$^\xi(  \Phi  ,j  )$''
construct
 and two new theorems that hybridize
the methodologies of
 Theorems \ref{pqq4} 
and \ref{pqq5}.
These  results will
improve upon
 Theorem \ref{pqq4} because 
their self-justifying  systems 
contain only a finite 
number of axiom-sentences beyond those lying in
$~\xi \,$'s base formalism
of  $B^\xi $.
They 
will
improve upon 
 Theorem \ref{pqq5} because
they can prove the important  
 Braced$^\xi(  \Phi  ,j  )$  
subset of
 $ (   \nal{B} , D )  $'s $\Pi_1^\xi$ theorems
in a full sense
 (rather than in
Remark  \ref{re4-2}'s 
weaker 
simulated respect).
Appendix G's results are useful
because {\it for arbitrary $k$}
and 
for
any of Appendix D's four sample configurations,
every
$\Pi_1^\xi$ theorem
of  $  (   \nal{B}  ,  D  )  $
 containing $~k~$ or
fewer bounded and unbounded quantifiers will 
be proven 
by its Theorem G.3 to be self-justifying
in an undiluted pure sense. (This is because each
such $\Pi_1^\xi$ sentence, with fewer than $~k~$
quantifiers, will 
lie in
 some
 fixed  Braced$^\xi(  \Phi  ,j  )$  
set $~$---$~$  where solely 
the value of
$~k~$ determines the values for
 $~\Phi~$ and $~j ~$. )
% \nop
\ennd 

%\njp

\cvmew
\cvs

It is 
probably
desirable
 to 
 concentrate 
primarily
 on
 Theorems \ref{ppp1}, \ref{ppp2},
\ref{pqq3}, \ref{pqq4}, \ref{pqq5}
and  \ref{ppp6}
during 
one's first reading of this paper.
This is because
Appendixes A-G
are
less central than these 
% five 
core theorems,
although their material
does
add 
several
useful
further
perspectives to this subject.

%%notre
  
\medskip

\cvb
\parskip 3pt

%%notre-only  
%%notre-only    {\bf Reminder to the Referee and Editor:} The cover
%%notre-only    page of this article indicated any (or all) of 
%%notre-only    Appendixes A-F can be removed from this article if it
%%notre-only    is deemed too long. In that case, the manuscript could
%%notre-only    cite the Cornell Archive paper 
%%notre-only    \cite{ww11}
%%notre-only    containing these appendixes.

%% 
%%   {\bf Reminder to the Referee:}
%% % and Editor:} 
%% The cover
%%   page of this article indicated any (or all) of 
%%   Appendixes A-F can be removed from this article if it
%%   is deemed too long. 
%% (This is because 
%% these results
%% are stored in \cite{ww11}'s archives.)

%% In that case, the manuscript could
%%   cite the Cornell Archive paper 
%%   \cite{ww11}
%%   containing these appendixes.

\section*{7.  Concluding Remarks}

The research in 
this article 
has been
 a continuation of our prior
research \cite{ww93}-\cite{ww9} that simultaneously
has simplified, unified
and extended the prior results. 
It
has explored self-justification
with a 3-part approach
where:
\bee
% \itemsep +3pt
\item
Sections \ref{3uuuu2} and \ref{3uuuu3}
introduced
three different
  stem components
that can be used to generate
 self-justifying
 systems. (These are
     the relatively simple Lemma  \ref{lemex4}
and the mathematically more sophisticated
Theorems  \ref{ppp1}  and \ref{ppp2}.)
\item
Section \ref{sect64} and Appendix G then
generalized our initial stem-like theorems
in the
six 
% five 
different directions 
formalized
by
Theorem  \ref{pqq3},  \ref{pqq4},
  \ref{pqq5}, 
  \ref{ppp6}, 
G.2 and G.3
\item
Appendix D 
subsequently
provided four examples of generic configurations
that are applications of 
Section \ref{sect64}'s results.
\ene
This 3-part approach is
very
different from the methods
used
in our  prior articles
\cite{ww93,ww5,wwapal,ww9}. The latter
examined
 particular isolated applications
in 
thorough
detail (rather than
compartmentalize
and separate
 the analysis
into three stages). The virtue of this 3-stage analysis
is 
it leads to many new theorems, in addition to
unifying our prior results.

%%cc \bigskip

It is desirable to 
categorize
the maximal
generality and 
strongest form of
boundary-case exceptions
for the Second Incompleteness Theorem
that are feasible because
G\"{o}del's centennial discovery 
 beckons the 
scholarly community
to sharpen their understanding of his 1931 
landmark discovery,
that
has
fundamentally
 reshaped mathematics.

%%cc \bigskip

It should be emphasized that
our over-all 
research
in \cite{ww93} -- \cite{ww9}
 has spent an approximately
equal 
effort
in
 exploring generalizations of
the Second Incompleteness Theorem
\cite{ww2,wwlogos,wwapal,ww7} 
and 
in
examining
its 
viable 
boundary-case exceptions
\cite{ww93,ww1,ww5,ww6,wwapal,ww9}
(although the current article
focused on
the latter topic).
 This is because
the
Second 
Incompleteness Theorem
is 
a starkly 
robust result that
 imposes 
sharp
limits
on how strong self-justifying systems may become.

%%cc \bigskip

Finally, we encourage the reader to take another brief glance
at Remarks \ref{f88} -- \ref{recc1}.
They
offer a brief
summary
of both the
strengths and limitations
of our
chief
 results.
They also explain how Theorem \ref{ppp6}'s reflection
principle for
$\Pi_1^\xi$ styled theorems is a  
very  unexpected
 result.

%quite unexpected
%and pleasing result.

% implications for  
% Remark \ref{remhappy} additionally explains
% what types of

%reflection principles. 

%%%%% that are associated with 
%%%%boundary-case

%%%%%  \bigskip
%  

%%cccorn-ff  
%%cccorn-cc  
%%cccorn-cc  

  \bigskip
  \bigskip

    {\bf ACKNOWLEDGMENTS: }
    I thank
    Bradley Armour-Garb,
    Seth Chaiken and
    Kenneth W. Regan
    for many useful suggestions about how to
    % significantly
     improve
    the presentation 
    style
    of some earlier drafts of this article.
    As was noted in the text of this article, some
    of my research
    during the last 17 years was influenced by some private 
    communications
    that I had with Robert M. Solovay 
    in 1994 and
     Leszek 
    Ko{\l}odziejczy  in 2005.

\newpage

\parskip 2pt

\cvt

\section*{Appendix A: The
$~\Pi_1^\xi~$ encoding for 
SelfCons$^k(\beta,d)$}

This appendix will  
summarize how to 
formalize
a $\Pi_1^\xi$ encoding for
Definition \ref{xd+1x7}'s 
SelfCons$^k(\beta,d)$ predicate.
It
will
use
the following notation:
\bee
\itemsep +1pt
\baselineskip = 1.0 \normalbaselineskip
\topsep -3pt
\itemsep 4pt
\item
Neg$^k(x,y)$,
will denote 
a
$\Delta_0^{\xi}$
formula 
indicating that $~x~$ is the G\"{o}del
number of a 
$\Pi_k^\xi$ sentence and
that $~y~$
represents the  
$\Sigma_k^\xi$ sentence
which is its logical negation.
\item
$\mbox{Prf} \, _\beta^d~( \, t \, , \, p \,)$ 
will denote a formula  designating
that
 $~p~$ is  a proof of theorem $~t~$ from the axiom
system  $~\beta~$ using the deduction method $~d.~$
\item
$\mbox{ExPrf} \, _\beta^d\,( \, h \, , \, t \, , \, p \,)$  
will denote that
$ \, p \, $ is a proof
(using $d$'s deduction method)
 of
a theorem $t$ 
from the union 
of the axiom system
$\beta$ with the added
sentence whose G\"{o}del number equals
$\, h \,$. 
\item
 $\mbox{Subst} \, ( \, g \, , \, h \, )$ will denote
G\"{o}del's
substitution formula --- which yields TRUE when $\, g \,$
is an encoding of a formula
and $\, h \,$ encodes a sentence
that replaces all occurrence of free variables in $g \,$ with
a term of $\bar{g}$ (that specifies $g$'s
G\"{o}del number). 
\item
$\mbox{SubstPrf} \, _\beta^d~( \, g \, , \, t \, , \, p \,)$  
will denote the 
hybridization of Items 3 and 4
that yields a Boolean value of TRUE 
when there
exists an integer $~h~$ 
satisfying 
  $\mbox{Subst} \, ( \, g \, , \, h \, )$ and
$\mbox{ExPrf} \, _\beta^d~( \, h \, , \, t \, , \, p \,)$.
\ene

It is easy to apply \cite{ww1}'s methodologies
to confirm 
Items 1-5 
can be encoded as $\Delta_0^\xi$ formulae.
Thus,  
Appendixes C and D of  \cite{ww1}
explained how  the theory of
LinH functions 
\cite{HP91,Kr95,Wr78}
implied  there existed 
 $\Delta_0$ encodings 
for  formulae 1-4,
and
these 
 $\Delta_0$ encodings can be
easily rewritten
 \footnote{  \b55 This rewriting
of conventional $\Delta_0$ formulae into
a  $\Delta_0^\xi$ format
 is possible because  
Part 2  of Definition \ref{def3.3}
indicated that
two 3-way predicates 
of
Add$(x,y,z)$
and Mult$(x,y,z)$ do 
encode  addition and multiplication
in a $\Delta_0^\xi$ styled form.}
as
 $ ~ \Delta_0^\xi ~$ expressions.
\ep{encode} 
 uses this information to formulate a
$\Delta_0^{\xi}$ encoding for 
$\mbox{SubstPrf} \, _{~\beta}^d \,(   g   ,   t   ,   p  )$'s 
graph. It is
equivalent to 
$~$``$~\exists ~h~ \{ ~\mbox{Subst}    (    g    ,    h    )~\wedge~
\mbox{ExPrf} \, _{\beta}^d(    h    ,    t    ,    p   )\,  \}  \, \,$''$,~$
  but \ep{encode} is written in
 a $\Delta_0^{\xi}$ format --- {\it unlike} the quoted
expression.
\begin{equation}
\label{encode}
\mbox{Prf} \, _{\beta}^d~( \, t \, , \, p \,)~~~\vee~~~\exists ~h\leq p
~~ \{ ~ \mbox{Subst} \, ( \, g \, , \, h \, )~\wedge~
\mbox{ExPrf} \, _{\beta}^d~( \, h \, , \, t \, , \, p \,)~ \}   
\end{equation}

Using \eq{encode}'s
  $\Delta_0^{\xi}$ encoding for 
$\mbox{SubstPrf} \, _{\beta}^d(   g   ,   t   ,   p  )$, it is
easy to encode 
SelfCons$^k(\beta,d)$
 as  a $\Pi_1^{\xi}$ 
axiom-sentence.
Thus,  let
$ \, \Gamma^k(g) \, $ 
denote  \eq{encode12}'s formula, and
let  $ \, \bar{n} \, $ denote $ \, \Gamma(g)$'s
G\"{o}del number.
\begin{equation}
\label{encode12}
\forall   \, x \, \forall   \, y \, \forall   \, p \, \forall   \, q 
~~~  \neg  ~~
 \{ \,\mbox{Neg}^k(x,y) ~ \wedge~ 
\mbox{SubstPrf} \, _{\beta}^d  \, (  g    ,    x    ,    p   ) ~
\wedge  ~
\mbox{SubstPrf} \, _{\beta}^d  \,    (    g    ,    y    ,    q   )  \, \}
\end{equation}
Then  $~$``$~\Gamma^k(~ \bar{n}~)~$''$~$ 
is a $\Pi_1^{\xi}$ 
encoding 
for SelfCons$^k(\beta,d)$'s formalization of
the statement $+$ from Definition
\ref{xd+1x7}.
Thus,  $~\Gamma^k(~ \bar{n}~)~$
is 
encoded is as follows:
\begin{equation}
\label{encode12n}
\forall   \, x \, \forall   \, y \, \forall   \, p \, \forall   \, q 
~~~  \neg  ~~
 \{ \,\mbox{Neg}^k(x,y) ~ \wedge~ 
\mbox{SubstPrf} \, _{\beta}^d  \, (  \bar{n}     ,    x    ,    p   ) ~
\wedge  ~
\mbox{SubstPrf} \, _{\beta}^d  \,    (     \bar{n}   ,    y    ,    q   )  \, \}
\end{equation}

{\bf Reminder about \ep{encode12n} :}
This sentence's
definition for 
SelfCons$^k(\beta,d)$ 
does not assure 
\eq{encode12n}
is
true under the Standard-M model. Indeed for nearly
all $(\beta,d)$, it will be false
when $k \geq 2$. This
is the reason that the study of
SelfCons$^k(\beta,d)$,
under Theorems \ref{ppp1} and \ref{ppp2},
has focused on the cases where $~k~$ equals 0
or 1.
Moreover,
the preceding
construction did assure that 
SelfCons$^k(\beta,d)$
had a
$~\Pi_1^\xi~$ encoding
because 
such
an
{\it ``I am consistent''} axiom 
carries 
 more meaning
than  a 
$\Pi_2^{\xi}$ encoded axiom.

\parskip 0 pt

\section*{Appendix B: The Proof of Lemma \ref{llpp3}}

\cvt

Lemma \ref{llpp3} is
a crucial interim step used to verify
each of
Theorems \ref{pqq3}, \ref{pqq4} and \ref{pqq5}.
Its proof will employ
the following three 
straightforward observations:
\bed
\item[Fact B.1 ]
Lemma  \ref{llpp3}'s 
hypothesis
implies that
 $\mxi$ 
is a
generic configuration. (This is because 
it specified that
 $\xi$ 
was a generic configuration  
and that
all the $\Pi^\xi_1$ sentences of
$~\mheta~$ were
true
in the
 Standard-M model.
Thus,
 $\mxi$ 
must also be a
generic configuration.)
\item[Fact B.2 ]
The
associative identity of 
 $ \, \theta \, \cup \, (\mheta \,  \cup B^\xi) \,  \, = \,  \, 
( \, \theta \, \cup \, \mheta \, ) \cup B^\xi \, $ 
obviously holds. It 
implies
a sentence
$ \, \Uxp \, $
is a theorem of
 $ \, \theta \, \cup \, (\mheta \,  \cup B^\xi) \, $ 
if and only if 
it is a theorem of 
 $ \, ( \, \theta \, \cup \, \mheta \, ) \cup B^\xi \, $.

\item[Fact B.3 ]
Lemma \ref{llpp3}'s hypothesis 
directly
\footnote{\s55 The identity of 
$~ \zzz \mheta ~)~=~\infty~$
must be true
because the 
hypothesis 
of Lemma  \ref{llpp3}
 indicated that 
all the $\Pi^\xi_1$ sentences in
$~\mheta~$ are 
true  under the Standard-M model.
By  Definition \ref{chg}, 
this
 implies $~ \zzz \theta ~)~= ~ \zzz \mheta \cup \theta ~)~$.}
implies $ \zzz \theta )=  \zzz \mheta \cup \theta )$.
\ennd

\smallskip

The justification of
claims (i)-(iv)
are
 consequences of 
Facts B.1 through B.3.
We
will
 provide a detailed proof of only Claim (i) here
 because all
four claims have 
similar proofs.

\medskip

{\bf Proof of Claim (i) :}
The hypothesis of Claim (i) indicated that $~\xi~$ 
was 0-stable. Therefore, it
satisfies
 Definition \ref{ostab}'s
invariant of $~***~$. 
In a context where
$~\phi ~$ is a variable
designating a r.e. set of $\Pi_1^\xi$ sentences 
and  $ \, \Uxp \, $ is a variable corresponding to a 
$\Delta_0^\xi$ sentence, 
 the invariant  $~***~$
can be rewritten in a  quasi-rigorous form as :
\beq
\label{whatshit1}
\forall ~ \phi ~~ \forall ~ \Uxp 
 ~~~~~  \mbox{\it the below
statement, called $~ \Psi_1(\phi,\Uxp)$,  is true  }
\end{equation}
\begin{quote}
$\forall ~p~~$ If $ \, \Uxp \, $ is a
$\Delta_0^\xi$ theorem 
 derived
from the
axiom system $ \, \phi \,  \cup B^\xi \, $ 
via a proof $ \, p \, , \, $ 
whose length satisfies
 Log$(p) \, \leq \zzz \phi \, ) \, +1 \, , \,  $
 then  $ \, \Uxp \, $ 
is true under the Standard-M model.
\end{quote}
Since \eq{whatshit1}'s
 universally quantified variable 
$~\phi~$
can designate any r.e. set of $\Pi_1^\xi$ sentences, it 
may
designate the object  ``$~\theta ~\cup~\mheta~$'', where
$~\mheta~$ is the r.e. set of  $\Pi_1^\xi$ sentences
defined by Lemma \ref{llpp3}'s hypothesis
(and $~\theta~$ is any second r.e. set of sentences).
Thus,  \eq{whatshit1} 
directly
implies :
\beq
\label{whatshit2}
\forall ~ \theta ~~ \forall  \Uxp
 ~~~~~  \mbox{\it the below
statement, called $~ \Psi_2(\theta,\Uxp)$,  is true  }
\end{equation}
\begin{quote}
$~\forall ~p~~$ If $ \, \Uxp \, $ is a
$\Delta_0^\xi$ theorem 
 derived
from the
axiom system $~(~\theta~\cup~\mheta~)~ \cup B^\xi~$ 
via a proof $~p~,~$ 
whose length satisfies
 Log$(p) \, \leq \zzz  \, \theta\cup\mheta)~   +1, $
 then  $ \, \Uxp \, $ 
is true  under  Standard-M.
\end{quote}
Facts B.2 and B.3 enable one to 
simplify 
\eq{whatshit2}'s
 terms of  $( \, \theta \, \cup \, \mheta \, ) \,  \cup B^\xi$ 
and  $\zzz ( \, \theta \, \cup \, \mheta \, ) \, ~    )$  
and thus to derive \eq{whatshit3} as a consequence.
\beq
\label{whatshit3}
\forall ~ \theta ~~ \forall  \Uxp ~~~~
   \mbox{\it the below
statement, called $~ \Psi_3(\theta,\Uxp)$,  is true  }
\end{equation}
\begin{quote}
$~\forall ~p~~$ If $ \, \Uxp \, $ is a
$\Delta_0^\xi$ theorem 
 derived
from 
axiom system $\theta~\cup~(~\mheta~ \cup B^\xi~)$ 
via a proof $~p$, 
whose length satisfies
 Log$(p) \, \leq \zzz \theta) \, + \,1 $,
 then  $  \Uxp  $ 
is true under the Standard-M model.
\end{quote}
We will now use Fact B.1's observation that
$\mxi$ 
is a generic configuration.
The sentence
\eq{whatshit3} 
indicates
this configuration
satisfies  
Definition \ref{ostab}'s invariant of
$***$.
Hence, 
$\mxi $
 is
0-stable.
$\Box$

%%cccorn-ff

  \medskip
  
  \cvs
  
  {\bf Brief Comments about the Justifications of Claims (ii)-(iv):}
  This appendix has  omitted 
  proving 
  (ii)-(iv) 
  for the sake of brevity.
  Their proofs are   similar
   to Claim (i)'s proof.
  For instance,   Claim (ii)'s proof differs
  from  Claim (i)'s proof  by having
  the 0-stability invariant in $***$
   replaced
  by the 
  A-stability invariant of $~*~$.
  This will  cause the analogs of
  \eq{whatshit1} -- \eq{whatshit3} 
  to undergo the 
  following two 
  simple
  changes under Claim (ii)'s proof :
  \bee
  \topsep -10pt
  \itemsep +2pt
  \item 
  $ \, \Uxp \, $ 
  will  represent a 
  $\Pi_1^\xi$ 
  (rather than a
  $\Delta_0^\xi$ ) theorem-statement under the 
  revised versions of $\Psi_1$, 
  $\Psi_2$ and $\Psi_3$  
   used in Claim (ii)'s proof
  \item 
  The requirement
  (in 
  sentences \eq{whatshit1}-\eq{whatshit3} of Claim i's proof)
   that
  $  \Uxp  $  be
  true in the Standard model
  is
   changed 
  to the stipulation that 
  $  \Uxp  $ 
     satisfies 
   the Good$\{ \, \tftt  \zzz \phi) \, \}$
  and
  Good$\{  \tftt  \zzz \theta)  \}$
  conditions
   under the revised forms of
   $\Psi_1$, 
  $\Psi_2$ and $\Psi_3$  
  used   to prove Claim (ii).
  \ene
  Other
  minor adjustments in Claim (i)'s proof
  shall verify Claims iii and iv.

\nop

\section*{Appendix C: The Proof of \phx{ppp2} }

Our  proof of  \phx{ppp2}
is a straightforward modification of
\phx{ppp1}'s 
 proof. 
It
will be divided into
two lemmas.

\medskip

{\bf Lemma C.\ref{B1-lem}.}
{\it 
Every generic configuration that is either E-stable or
A-stable will automatically satisfy
Definition \ref{ostab}'s 0-stability condition.}

\medskip

{\bf Proof.} $~$
Lemma C.\ref{B1-lem} is 
a 
consequence of
the 
``Special Note''  appearing at the end of 
Definition \ref{xd+1x4}.
For every  $N\geq 0$,
it
 indicated that
if
$\Uxp$ is a 
$\Delta_0^\xi~$ 
sentence then
Scope$_E$($\Uxp,N)$
is  equivalent to  $~\Uxp~$.
This implies 
(via Definition \ref{xd+1x5})
that if 
$\Uxp$ is a  
Good(N)
$\Delta_0^\xi~$ formula
then Scope$_E$($\Uxp,N)$ is automatically
true under the Standard-M model.

The latter observation makes it easy to confirm
Lemma C.\ref{B1-lem}. This is 
because  every
$\Delta_0^\xi$ sentence
is a
$\Pi_1^\xi$ and
$\Sigma_1^\xi$ statement.
Hence, the application of the invariants
$~*~$  and $~**~$
from Definitions 
\ref{astab} and \ref{estab}, in the degenerate case
where
$\Uxp$ is a
 $\Delta_0^\xi$ theorem,
corroborates 
Lemma C.\ref{B1-lem}'s claim
(by showing that
 Definition \ref{ostab}'s 
 invariant of
 $\,***\,$ does hold).
$~~\Box$

\medskip

The remainder of this appendix will  focus on
Definition \ref{ostab}'s 0-stability condition.
(This is sufficient to justify \thx{ppp2}
because Lemma C.\ref{B1-lem}
showed all 
 E-stable and A-stable
configurations are
 0-stable.)

\smallskip

{\bf Lemma C.\ref{B2-lem}.}
{\it 
Let $~\xi~$
denote a
generic configuration
 \gggen 
that is
0-Stable.  Then 
the axiom system of
$B^\xi+$SelfCons$^0(B^\xi,d)$
will be consistent (and 
hence self-justifying).}

\smallskip

\smallskip

{\bf Proof:}
It will be awkward 
to write repeatedly
the expression   ``$~B^\xi+$SelfCons$^0(B^\xi,d)~$''
during our proof.
Therefore,  $\zhz $ will
be an abbreviated
name for this 
system. 
Our justification
 of Lemma C.\ref{B2-lem},
is similar to \phx{ppp1}'s proof,
$\,$except  it replaces a 
Level$(1^\xi)$  form of self-justification
with  
a 
Level($0^\xi)$.
It
will thus
  be
 abbreviated and 
use the following notation:
\bee
\itemsep +4pt
\item
$~\mbox{Prf}_\zhz (t,p)~$  is a
$\Delta_0^\xi$ formula 
specifying
$~p~$ 
is
 a proof of the theorem $~t~$ from
the  axiom system $~\zhz ~$
(using 
$~\xi \,$'s
deduction method of $d$ ).
\item
$~\mbox{Neg}^0(x,y)~$ 
 is a
$\Delta_0^\xi$ formula 
 indicating
$~x~$ is
the G\"{o}del encoding of a
 $\Delta_0^\xi$ 
sentence and
$~y~$ is a   $\Delta_0^\xi$ 
sentence  representing
$~x \,$'s negation.
\ene 
Expression \eq{qpencode} 
denotes
 $~\zhz \,$'s 
 Level$(0^\xi)$
self-justification axiom.
It is encoded using
Appendix A's
methodology,
similar to 
its counterpart 
used in \phx{ppp1}'s proof (i.e.  
\ep{pencode} ).
\begin{equation}
\label{qpencode}
\forall   \, x \, \forall   \, y \, \forall   \, p \, \forall   \, q 
~~~  \neg  ~~
 \{ \,\mbox{Neg}^0(x,y) ~ \wedge~ 
\mbox{Prf}_\zhz   \, (     x    ,    p   ) ~
\wedge  ~
\mbox{Prf}_\zhz  \,    (       y    ,    q   )  \, \}
\end{equation}

Our proof of Lemma C.\ref{B2-lem}
will be a proof by contradiction.
It will thus
 begin with
the 
contrary
assumption that $~\zhz ~$ is inconsistent
and have  $~\Phi~$ denote
\eq{qpencode}'s sentence.
The inconsistency of $~\zhz ~$ 
 implies that 
$~\Phi~$ is
false under the Standard-M model.
Hence via
Definition \ref{chg}, we get:
\begin{equation}
\label{qpunk}
 \zzz \Phi~) 
   ~~ < ~~ \infty
\end{equation}
\noindent
\ep{qpunk} 
implies there  exists
a tuple $(\bar{p},\bar{q},\bar{x},\bar{y})$
satisfying \eq{qpencodepunk}.
(This is because such a 
$(\bar{p},\bar{q},\bar{x},\bar{y})$ 
corroborates \eq{qpunk}'s implication that a 
counter-example to
\eq{qpencode}'s sentence 
does
exist.)
\begin{equation}
\label{qpencodepunk}
  \,\mbox{Neg}^0(~ \bar{x}~, ~ \bar{y}~) ~ \wedge~ 
\mbox{Prf}_\zhz   \, (     ~ \bar{x}~    ,    ~ \bar{p}~   ) ~
\wedge  ~
\mbox{Prf}_\zhz  \,    (       ~ \bar{y}~    ,    ~ \bar{q}~   )  
\end{equation}
% Furthermore, the  
The
$(p,q,x,y)$
 satisfying \eq{qpencodepunk}
with  minimum value for 
$ \mbox{Log}\{\, \mbox{Max}[p,q,x,y] \, \}$ 
will
additionally
%also
satisfy
\eq{qpunkless}.
(This  observation follows
from the analog of the
Footnote \ref{footp1}
appearing 
in \phx{ppp1}'s proof. 
Thus, Section \ref{3uuuu3}'s Equations
of
\eq{pencodepunk}
and
\eq{punkless}
are 
the
analogs of the current
 Equations
\eq{qpencodepunk}
and
\eq{qpunkless}.
The 
Footnote \ref{footp1} 
showed
the particular
$(p,q,x,y)$
 satisfying \eq{pencodepunk}
with  minimum value for 
$ \mbox{Log}\{~\mbox{Max}[p,q,x,y] ~\}$ 
satisfied
\eq{punkless}.
By the same reasoning,  
the minimal 
$(\bar{p},\bar{q},\bar{x},\bar{y})$
satisfying
\eq{qpencodepunk} will     satisfy
\eq{qpunkless}.)
\begin{equation}
\label{qpunkless}
 \mbox{Log}~\{~\mbox{Max}[~\bar{p}~,~\bar{q}~,~\bar{x}~,~\bar{y}~] ~~\}
   ~~~ =~~~\zzz \Phi~) ~~+~~1
\end{equation}

Equations 
\eq{qpencodepunk} and \eq{qpunkless} 
shall
bring 
Lemma C.\ref{B2-lem}'s
 proof-by-contradiction to its 
sought-after
 end. Thus,
let $~\Upsilon~$ denote the
 $\Delta_0^\xi$ 
sentence specified
by $~\bar{x}~$. Then
$~\neg ~ \Upsilon~$ 
corresponds to the 
 $\Delta_0^\xi$ 
 sentence denoted
by $~\bar{y}~$. 
\ep{qpunkless}
indicates
that
both
 $~\Upsilon~$ and
$~\neg ~ \Upsilon~$
have  proofs 
such that the logarithms of
their G\"{o}del numbers are bounded
by
$~\zzz \Phi~) ~+~1~$.
Using
Definition \ref{ostab}'s
invariant of $~***~$,
these facts establish  
\footnote{\s55 To  
apply Definition \ref{ostab}'s
 invariant $~***~$ in the present setting,
one simply sets $~\theta\,$'s R-View 
equal to $~\Phi \,$'s
1-sentence statement. 
Then $~***~$ implies that both 
 $~\Upsilon~$ and
$~\neg ~ \Upsilon~$
must
be 
true under the Standard-M model
because both  their proofs had lengths 
$\, \leq ~\zzz \Phi~) ~+~1~$.
} that both 
 $~\Upsilon~$ and
$~\neg ~ \Upsilon~$
are
true under the Standard-M model.

But it is impossible for
  a sentence and its
negation to be both
true.
This finishes
Lemma C.\ref{B2-lem}'s
proof
because the
temporary
 assumption that 
$\zhz $ was inconsistent
has
led to a
contradiction.
$~~\Box$

\medskip

\phx{ppp2} is a
 consequence of 
the  Lemmas C.\ref{B1-lem} and
 C.\ref{B2-lem} because
Lemma
 C.\ref{B2-lem}'s formalism 
generalizes
to all
E-stable and A-stable configurations
via Lemma C.\ref{B1-lem}'s reduction methodology.

\nop

\section*{Appendix D: Applications and Examples}

This appendix will illustrate
four examples of generic configurations 
that satisfy
Theorems \ref{ppp1} and \ref{ppp2} 
(and 
which
therefore  are self-justifying).
It will be divided into three parts.
Section D-1 will define our first example of
self-justifying configuration, called $~\xi^*~$.
It will use
semantic tableaux
 deduction.
Section D-2 will 
prove that          $~\xi^*~$ is  EA-stable. 
Section D-3 will 
briefly sketch
        three additional examples of
   stable generic      configurations

It is likely preferable to examine
Sections \ref{3uuuu1} -- 
\ref{sect64} before  this appendix.
However, Section
 D-1's 
 short 2-page discussion 
can be
 read
quite easily
either before or
 after 
Section \ref{sect64}.

\subsection*{D-1. $~$Definition of
the EA-Stable Configuration  
 $~\xi^*~$ }

Our first example of an 
EA-Stable configuration,
 called \xxi,
will be defined in
this section.
Its deduction method will be
semantic tableaux.
Its 
 base axiom 
system $\, B^* \,$ will be a
 Type-A  formalism,
which 
treats addition but not
multiplication as a total function
(i.e. see  \ep{totdefxa} $~).$

The closest analog of  $B^*$
and  \xxi 
 in our prior work
% had
appeared in Section 5 of \cite{ww5}.
%(Its formalism was
(It 
differed
from  \xxi 
partly because
it
%\cite{ww5}
did not 
use Definition  \ref{def3.3}'s
unifying
notation.)
In 
\cite{ww5},
% \cite{ww5}'s  discussion,
a function 
 $\, F \, $ 
was     called
 {\bf Non-Growth} 
iff  
$ F(a_1,a_2,...a_j) 
\leq   Maximum(a_1,a_2,...a_j)$
for all  $a_1,a_2,...a_j$. Six  examples of  
non-growth functions are:
\begin{enumerate}
\item
{\it Integer Subtraction} 
where ``$~x-y~$'' is defined to equal zero 
in {\it the special case} where
 $~x \leq y,$
\item
{\it Integer 
Division}
where ``$~x \div y~$'' equals
$~x~$ when $~y=0,~$ and
it equals $~\lfloor ~x/y ~\rfloor~$ otherwise,
\item
$Root(x,y)~$ which equals $~ \lceil ~x^{1/y}~ \rceil$ when $~y\geq 1,~$
and
it equals $~x~$ when $~y=0.$
\item
$Maximum(x,y),~~$
\item
$ Logarithm(x)~=~\lfloor ~$Log$_2(x)~ \rfloor~$ when $~x \geq 2,~$
and  zero otherwise. 
\item
$Count(x,j)~~=~~$the number of ``1'' bits
among $~x$'s rightmost $~j~$ bits.
\end{enumerate}
These operations were called
{\bf Grounding Functions}
in \cite{ww5}.
The term
{\bf U-Grounding Function}  
referred 
to a
set of
functions that included the
Grounding  operations
 plus 
the further 
primitives  of addition and
{\it Double$(x)=x+x$.} 
(The Double operation  
is helpful because it significantly
 enhances
\cite{ww5}'s linguistic efficiency \footnote{\s55 The symbol
{\it Double$(x)$}
was technically unnecessary in \cite{ww5}'s formalism
because
$x+x$ can encode
 Double$(x)$.
However,
its notation
adds expressive power to
\cite{ww5}'s language because, for example,
Double(Double(Double(Double(x))) 
requires less memory space to encode than
than
$~x~$ added to itself 16 times.
} .)

\medskip

The symbol $~\Delta_0^*~$ 
will be the
analog of Definition \ref{xd+1x1}'s 
 $\Delta_0^\xi$
construct under
the U-Grounding function language.
(It will be defined to be any formula in a U-Grounding 
language where all its quantifiers are bounded.)
It is easy in this context to encode
a $\Delta_0^{    *}$ 
formula $Mult(x,y,z)$ 
for 
representing
  multiplication's graph.
For instance,
\ep{neweq1} 
is one such $\, \Delta_0^{    *} \, $ formula
(which actually does not employ any bounded quantifiers):
\begin{equation}
\label{neweq1}
[~(x=0    \vee    y=0 ) \Rightarrow z=0~ ]~ ~\wedge ~~ 
[~(x \neq 0 \wedge y \neq 0~) ~ \Rightarrow ~
(~ \frac{z}{x}=y  ~\wedge \, ~  \frac{z-1}{x}<y~~)~]
\end{equation}
Expression
\eq{neweq1} 
is significant because Part 2
of Definition \ref{def3.3}
indicated every generic configuration must have 
 available some 
method to represent the graphs of
addition and multiplication 
in a
 $\Delta_0^\xi$ 
styled
format, similar to     
\eq{neweq1}'s paradigm.
 (Addition can be
treated trivially
because the 
U-grounding language
           possesses 
an addition function symbol.)
The
footnote
 \footnote{\s55  Section \ref{3uuuu1} 
explained
during
its discussion
of  
Lemma \ref{lex22}
 that
every generic configuration $  \xi  $
must have a means to encode 
the graphs of
addition and multiplication
as $  \Delta_0^{\xi}  $ formulae
(visavis Parts 1 and 2 
of 
 Definition \ref{def3.3}).
This enabled Lemma \ref{lex22}'s
procedure to translate
all of conventional arithmetic's 
$\Sigma_j$ and $\Pi_j$
formulae 
into
equivalent 
$\Sigma_j^{\xi}$ and $\Pi_j^\xi$
expressions. \fend }
% in $ \, \xi \,$'s language. 
serves as a reminder about why
these $\Delta_0^*$ encodings are 
needed. 
Our next goal is to define 
the generic configuration 
 \xxi that
 Section D-2 will prove  
is EA-stable.

%\hgskip

\bigskip

{\bf Definition
D.\ref{D1-def}.}
The 
language
 $~L^{*}~$ of the generic configuration 
 \xxi 
will be built in a natural manner
out of 
the eight 
U-grounding  function operations,
the usual
atomic predicate
symbols of ``$~=~$'' and 
 ``$~\leq~$'',
and the
three constant symbols 
 $~K_0,~ K_1\,$ and $\,K_2\,$ (that define the
integers of 0, 1 and 2).
The
other components
of  \xxi 's configuration
  are defined
below:
\bed
\cvh
\itemsep 4pt
\item[    i  ]
As previously noted,
  $~\Delta_0^{*}~$
is defined to
represent the set of all
formulae in $ \, L^* \,$'s 
language, 
whose  quantifiers are bounded 
in an
arbitrary
 manner by terms employing 
the U-Grounding function
symbols.
(It will thus generate via Definition \ref{xd+1x1}
the  $\Pi_n^*$ and $\Sigma_n^*$ sentences of
$ \, L^* \,$.)
\item[   ii  ]
The base axiom system 
 $~B^*~$ 
for \xxi 
will be allowed to be any
consistent
set of   $\Pi_1^{\xi^*}~$ sentences
that is capable of proving every
 $\Delta_0^*$ sentence that is valid
in Standard-M.
It will 
also
include   sentence
 \eq{d1a}'s 
{\it very precise} \footnote{\g55 \label{fodii}
The two appearances of the term ``$~x+y~$'' in sentence
\eq{d1a} 
may at first appear to be redundant. (This
statement is 
equivalent to sentence  \eq{totdefsymba}'s
declaration that addition is a total function,
which had avoided such redundancy.) The virtue of
\eq{d1a}'s format is that it is a    
 $\Pi_1^*$ styled 
statement, unlike \eq{totdefsymba}'s
 $\Pi_2^*$ styled format. 
This sharpened
 $\Pi_1^*$ perspective
will help simplify some of 
our proofs. }
$\Pi_1^*$  styled
declaration that addition is a total function.
\beq 
\label{d1a}
\small
\forall x ~~~~~~\forall y~~~~~~ \exists~ z \leq  x+y ~~:~~~ ~~~
  \{~~~z=x+y~~~\}
\enq
\item[  iii  ]
\xxi's deduction method will be the semantic tableaux
method.
\item[   iv  ]
\xxi' s 
G\"{o}delized
method $~g~$
for encoding
a semantic tableaux proof
 can be essentially
any natural method 
 that satisfies 
the 
minor stipulation
that
 at least $~5 \, J~$ bits
are required
 to encode a semantic tableaux
proof that has $~J~$ function symbols. 
This stipulation is called the 
{\bf ``Conventional Tableaux Encoding Requirement''}.
It is trivial 
\footnote{\g55 The Conventional Tableaux Encoding Criteria
requires that the 
 G\"{o}del 
number of a semantic tableaux proof,
with  $~J~$ function symbols,
 must be 
least as large as $~32^J~$.
It is clear that
all
the usual methods for 
generating the G\"{o}del codes
satisfy this criteria.
This is because 
any proof that has  $~J~$ function symbols will contain at least
 $~2  \, J~$ logical symbols and thus employ 
at least
 $~5  \, J~$ bits. \fend
} to corroborate   
 that all the
usual methods for encoding semantic tableaux proofs
satisfy this 
criteria.
(The  Appendix A of \cite{ww5} provides
one 
example of 
a possible tableaux encoding method.
Any other natural mechanism for encoding 
tableaux proofs is equally suitable.)
\ennd

\smallskip

Section D-2
will,
interestingly,
prove  that   $\xi^*$ is EA-stable.
This 
will imply (via \phx{ppp1})
that
$~B^*~\cup~\mbox{SelfCons}^1\{~B^*~,d~\}$ 
is 
self-justifying.
Theorems \ref{pqq4},  \ref{pqq5}, G.2 and G.3 will,
formalize, in this context,
four different methods in which 
$~B^*~$ can be extended to 
construct self-justifying formalisms that are able to prove
Peano Arithmetic's 
$~\Pi_1^*$ theorems.

\nskip 

\medskip

 Thus while  self-justifying axiom systems
contain unavoidable weaknesses, they  also possess
the nice feature that they are able to prove many
of the useful theorems of mathematics.

\nskip

\subsection*{D-2. $~$Proof of the  EA-Stability of
 $~\xi^*~$ }

\nskip

This section
 will prove \xxi 
is EA-stable and thus
satisfies the paradigms of 
Theorems \ref{ppp1}, \ref{pqq3}, \ref{pqq4}, \ref{pqq5}
and  \ref{ppp6}.
Our proof will be based on modifying some of the methodologies
from \cite{ww5}, so that they
become applicable to \xxi .
Many readers may prefer to omit examining
both this part of Appendix D and Section 
 D-3 because they 
are unnecessary for understanding
the 
% further
 material in
Section \ref{sect64} and
Appendixes E and F. (Our recommendation is that the latter
material be read first.)

\nskip
\smallskip

Our notation for defining a
  semantic tableaux
proof
        in the next paragraph
will be similar to the 
conventional
      definitions
appearing in 
Fitting's and Smullyan's textbooks  \cite{Fi90,Sm68}.
It will employ
\cite{ww5}'s
notation
so that we can employ
two of its
lemmas during our analysis of semantic tableaux proofs.

\nskip

%\cvmew

{
%notre \cvh
In our discussion,
$~\Phi~$
will be called
a {\bf  Prenex-Level($\,$m$^* \,$)}
  sentence
 iff it is a 
 $\Pi_m^*$ or $\Sigma_m^*$ 
expression
that satisfies the usual prenex requirement (that
all its unbounded 
quantifiers   lie  
in its leftmost part).
If $\Phi$ is  Prenex-Level$(m^* )$
then {\bf Reverse($\Phi$)} 
shall 
denote a second
 Prenex-Level$( \, m^* \,)$
sentence that is equivalent to
$~\neg~\Phi~$  rewritten 
\footnote{\s55  For example,
if $\Phi$ denotes ``$~\forall \, x ~ \exists \,y ~~\psi(x,y)~$''
then Reverse$(\Phi)$ 
would be written  as
 ``$~ \exists \, x ~ \forall \,y ~~\neg~\psi(x,y)~$''.}
%% \footnote{\s55 Thus  
%% % For example,
%% if $\Phi$ denotes ``$~\forall \, x ~ \exists \,y ~~\psi(x,y)~$''
%% then Reverse$(\Phi)$ is 
%% %would be written  as
%%  ``$~ \exists \, x ~ \forall \,y ~~\neg~\psi(x,y)~$''.}
in a
 Prenex Level$( \, m^* \, )$ form.
For a fixed axiom system $ \, \alpha ,  \, $ its
{\bf $\Phi$-Based Candidate Tree}
will  be defined
to be a  tree structure 
whose root is
the sentence 
Reverse$(\Phi)$ 
 and whose all other nodes are
either axioms of $ \, \alpha \, $ or deductions from higher
nodes of the tree,
via the rules 1--8 
given
below. (The symbol
``{\bf $ \,  $A$  ~  \Longrightarrow  ~  $B$  \, $}''
in  rules 1-8
will mean
that {\bf $ \,   $B$  \, $} is a valid deduction
from its ancestor {\bf $ \,  $A$  \, $}
in the germane 
 deduction tree.)
\begin{enumerate}
\itemsep 5pt
%notre \small
%%corbl \baselineskip = 1.05 \normalbaselineskip
\item $~ \Upsilon \wedge \Gamma \, ~ \Longrightarrow ~ \, \Upsilon
~$ 
and 
$~ \Upsilon \wedge \Gamma \, ~ \Longrightarrow ~ \, \Gamma ~$ .
$~~~$
\item $ \,  \neg  \,\neg \, \Upsilon  \,  \Longrightarrow  \,  \Upsilon. \, $  
Other rules for
the ``$ \, \neg \,$'' symbol  are:
$ \, \neg ( \Upsilon \vee \Gamma )  \,  \Longrightarrow  \,  \neg \Upsilon
\wedge \neg \Gamma$,
%\newline
$ \, \neg ( \Upsilon \rightarrow \Gamma )  \,  \Longrightarrow  \,   \Upsilon
\wedge \neg \Gamma \, $,
$ ~~~~\, \neg ( \Upsilon \wedge \Gamma )  \,  \Longrightarrow  \,  \neg
\Upsilon \vee \neg \Gamma \, $,
 $~ \,   \neg \, \exists v \, \Upsilon  (v)  \,  \Longrightarrow  \,  
\forall v   \neg \, \Upsilon  (v)  \, $
 and $ ~~\,   \neg \, \forall v \, \Upsilon  (v)  \,  \Longrightarrow  \,  
\exists v \,  \neg  \Upsilon  (v)$
\item A pair of sibling nodes $~ \Upsilon ~$ and $~ \Gamma ~$ is
allowed
when their ancestor is
$~\Upsilon \, \vee \, \Gamma.~$
\item A pair of sibling nodes $ \,  \neg \Upsilon  \, $ and $ \,  \Gamma  \, $ is
allowed
when their ancestor is
$ \, \Upsilon \, \rightarrow \, \Gamma$.
\item $~ \exists v \, \Upsilon  (v) ~ \Longrightarrow ~ \, \Upsilon(u) ~$
where $\,u \,$ is a newly introduced ``Parameter Symbol''. 
\item 
 $~ \exists v \leq s ~ \, \Upsilon  (v) ~~ \Longrightarrow ~ ~
u \leq s ~ \wedge~ \Upsilon(u) ~$
is the
 variation of Rule 5 for  bounded existential quantifiers
of the form $~$``$~ \exists v \leq s ~$''.
\item $\forall v \, \Upsilon  (v)  \,  \Longrightarrow    \, \Upsilon(t)  \, $
where $t$ denotes a U-Grounded 
term. These terms may be any one of
a constant symbol,
a  parameter symbol (defined  by a 
prior application of  Rules 5 or 6 to some
 some ancestor of the current node),
or a U-Grounding function-symbol with  
recursively defined
inputs.
\item 
 $\forall v \leq s  \, \Upsilon  (v) ~~ \Longrightarrow ~~
t \leq s \, \rightarrow \, \Upsilon(t) $ is 
the variation of Rule 7 
for a
 bounded  quantifier such as
 ``$~ \forall v \leq s ~ $''
\end{enumerate} }

%%corbl  \baselineskip = 1.1 \normalbaselineskip 
\noindent
Let us say a leaf-to-root branch (in a candidate
tree) is {\bf Closed} iff it contains both some sentence
$ \,  \Upsilon  \, $ and its negation ``$ \,  \neg \, \Upsilon  \, $''.  
Then a
 {\bf  Semantic
Tableaux Proof} of $  \Phi  $,
from the axiom system  $  \alpha $,
is defined to be a 
$\Phi$-Based Candidate Tree
whose every
 leaf-to-root branch is closed.

\cvt

\hgskip

%\cvwew

\nskip

%bbbbb {

It is next helpful to define
the notion of a 
{\bf Z-Based Deduction Tree},
in a context where $~Z~$ represents 
an
 axiom system,
typically different from 
the prior paragraph's
$~\alpha~$.
This object 
will be
defined to be
identical to
a semantic tableaux proof, except for the following
changes:
\bed
%%corbl  \baselineskip = 1.1 \normalbaselineskip 
\itemsep 5pt
\item[   i   ]
{\it Every node} in a $Z-$Based 
deduction tree
must be either an axiom of  $~Z~$
or a deduction from a higher node of the tree via
the 
rules 1-8.
(This
applies
also to the root of
a $Z-$Based  deduction tree.
It will  store an 
axiom of  $~Z~$
 in its root, unlike a semantic tableaux proof
 which had
stored
 Reverse$(\Phi)$ 
in its
root.)
\item[  ii   ]
There will be no requirement that each 
 leaf-to-root branch
be
 closed
in a 
$Z-$Based  deduction tree.
(Indeed, 
some branch will 
automatically
not be closed
if  $~Z~$ is  consistent.) 
\ennd
Items (i) and (ii) make it 
apparent
 that
$Z-$Based  deduction trees
are different 
% properties
from 
semantic tableaux proofs.
It will
turn out,
nevertheless, 
 that the study of
$Z-$Based  deduction trees
will clarify the nature of 
 semantic tableaux proofs.

\nskip
\nskip

\medskip

{\bf Definition 
D.\ref{D2-def}.}
Let $ ~   a ~   $ and $ ~   b ~   $ denote two
integers that are powers of $\,  2 \,  $ 
satisfying
$  a \, > b \,
\geq \, 2 \, $
Then
 an axiom system $ \, Z \, $
(employing $ \, L^* \,$'s language)
will be called
a {\bf Normed(a,b)} formalism  iff:
\begin{enumerate}
\item All  $Z$'s 
axioms  are 
either $\Pi_1^{*}$
or  $\Sigma_1^{*}$ sentences.
\item 
Each  $\Pi_1^{*}$ axiom of $~Z~$
will satisfy Definition \ref{xd+1x5}'s
Good($~$Log$_2a~)$ 
criteria, and
 each
 $\Sigma_1^{*}$ axiom of $~Z~$
will
likewise
  satisfy Good($~$Log$_2b~)$. 
\ene

\nskip
\nskip

{\bf  Clarification about Definition D.\ref{D2-def}
: } The
``Normed(a,b)''  
concept (above)
is obviously 
equivalent to the same-named notion
appearing in
Definition 4 of \cite{ww5}.
It uses,
however,
 a 
different notation
to make it compatible 
with
Section \ref{3uuuu2}'s formalism.
Thus, Item 2's assertion
that the
  $\Pi_1^{*}$ axiom 
$ \, \forall  \,  v_1 \,  \, \forall  \,  v_2 \,  \,  ...
\forall  v_k   \,  \,  \,  \, \phi(v_1,v_2,...v_k) \,  \, $
satisfies
Good($~$Log$_2a~)$  is equivalent to 
\eq{normscopede}'s statement.
The
Good($\,$Log$_2b\,)$ property of 
$\exists \,v_1\,\exists \,v_2\, ...\,\exists \,v_k\, ~~\phi(v_1,v_2,...v_k)~$
is,
likewise,
 equivalent to \eq{normscopedx}.

\nskip

\beq
\label{normscopede}
\forall ~ v_1~ < ~a~~\forall ~ v_2~ < ~a~~ ...
\forall ~ v_k~ < ~a
   ~~~~~ : ~~~~~
\phi(v_1,v_2,...v_k)~~.
\enq

\nskip

\beq
\label{normscopedx}
\exists ~ v_1~ < ~b~~\exists ~ v_2~ < ~b~~ ...
\exists ~ v_k~ < ~b
 ~~~~~ : ~~~~~
\phi(v_1,v_2,...v_k)~~.
\enq

\parskip 5pt

\nskip

Our interests in
this notation
will center around Fact D.3's invariant:

{
%%corbl  \baselineskip = 1.17 \normalbaselineskip 

\medskip

\nop

{\bf Fact D.3$~.$}
Let
$ \, \xi^* \, $ 
denote 
 Definition D.\ref{D1-def}'s
generic configuration, and 
$  Z $ be an extension of $ \, \xi^* \,$'s
base axiom
system  $ B^*  $ 
which
satisfies Definition D.\ref{D2-def}'s
Normed$(a,b)$ constraint.
Then any
$Z-$Based  deduction tree
$\, T \, $that
has a G\"{o}del
number 
smaller than
 $ \, (a/b)^4 \,$
must 
contain
at least one 
root-to-leaf branch, called $ \, \peta \,  \, , \, $
that is not
``closed''. {\rm
(In other words, this path   $ \, \peta \, $
will be contradiction-free, insofar as it
does not contain  both  some
sentence $~\Psi~$ and its formal
negation).}
}

%%notre
  
%%notre-only  {\bf Proof:} $~$ The  justification of
%%notre-only  Fact D.3 is essentially 
%%notre-only  an
%%notre-only  %a direct and
%%notre-only  immediate
%%notre-only   consequence
%%notre-only  of the
%%notre-only  %% no  
%%notre-only  %% no \newpage
%%notre-only  %% no \noindent
%%notre-only  %% no 
%%notre-only  Lemmas 1 and 2 from \cite{ww5}
%%notre-only  after one maps these two lemmas's notation into the
%%notre-only  context of Fact D.3's
%%notre-only  hypothesis. (We provide a detailed justification
%%notre-only  of this 
%%notre-only  technical
%%notre-only  fact in the Footnote 23 of the slightly longer
%%notre-only  unabridged version of this article 
%%notre-only  that
%%notre-only  resides in the Cornell Archives \cite{ww11}.) 

  \smallskip
  
\nskip

%%cccorn-ff

%%cccorn-cc  
%%cccorn-cc  
   {\bf Proof:}
   $~$ The 
   justification of
   Fact D.3 is
% essentially 
%   an
   a direct 
% and    immediate
    consequence
   of 
   %% no  
   %% no \newpage
   %% no \noindent
   %% no 
the   Lemmas 1 and 2
appearing in article \cite{ww5} 
   (see footnote 
   \footnote{\label{newtry}
   \baselineskip = 1.3 \normalbaselineskip  
   A 
   %formal 
   proof of Fact D.3 from first principles would
   be 
   quite
   complicated because there are eight elimination rules employed
   by semantic tableaux deduction, $~$each
   of which needs to
   be 
   % simultaneously 
   examined by such a proof's umbrella
   formalism. Fortunately, we do not need provide such a
   complicated analysis here 
   because  
   a 4-page proof of
   the Lemmas 1 and 2 
   in Section 5.2
   of \cite{ww5}
   had 
    already
   visited these issues. 
   Thus,
    Fact D.3 
   turns out to be an easy consequence 
   of these two lemmas
   %the Lemmas 1 and 2 of \cite{ww5}
   after the following two straightforward issues are addressed: 
   \bee
   \topsep -4pt
   \itemsep -2pt
   \item Section 5.2 of \cite{ww5} 
   had
   defined the $\,$``U-Height''$\,$ of a deduction tree to
   be the 
   largest
   number of 
   U-Grounding function symbols
   that appear in any of its root-to-leaf
   branches. Its Lemma 1
   proved that
   every deduction tree with a U-Height  $\, \leq \, Log_2a  \, - \, Log_2b \, $
   will contain at least one branch
   satisfying a condition,
   which \cite{ww5}
    called ``Positive(a,b)''.
   The Lemma 2 in \cite{ww5} then showed that this
   Positive(a,b) property  implies that the germane 
   deduction tree must contain some branch that is 
   contradiction-free.
   The combination of these two lemmas
   thus amounts to the establishing of
   the following 
   %slightly 
   rephrased 
   hybridized statement:
   %invariant:
   \begin{quote}
   $\bullet~~~$
   If a Z-based deduction
   tree  has 
                   a U-Height  $\, \leq \, Log_2a  \, - \, Log_2b \, $,
   then some branch of it 
   is contradiction-free
   (i.e.  this branch cannot
   contain  both  some
   % $\Delta_0^*$
   sentence $~\Psi~$ and its
   negation).
   \end{quote}
   \item
   %% 
   %% Item $\, \bullet \,$
   %% applies to Fact D.3 's  
   %% deduction trees
   %% because 
   %% 
   Fact D.3 's   hypothesis indicated the
   G\"{o}del number $g$ for its deduction tree
   satisfied the following 
   conditions:
   \bed
   \topsep -7pt
   \itemsep -1pt
   \item[   I.   ]
   $~~g~\leq ~ (a/b)^4~$
   \item[  II.   ]
   The   U-Height of $~g \,$'s deduction tree
   is less than $~\frac{1}{5}~$Log$_2~g~$. (This
   is
   simply
    because
    Fact D.3 presumes that
   the ``Conventional Tableaux Encoding'' methodology
   from Part-iv of 
   Definition D.\ref{D1-def}
   was used to encode $~g$'s
   G\"{o}del number.)
   \ennd
   Items
   I and II imply 
   $g \,$'s  tree has  
    a U-Height  $\, \leq \, Log_2a  \, - \, Log_2b  $.
   The invariant  $\, \bullet \,$ 
   %then, in turn, implies
   %that this deduction tree
   then implies
   this deduction tree
   has at least one branch
   that is contradiction-free (as  Fact D.3
   claimed). $\Box$
   \ene
   We emphasize that the 
   %The
   above 
   % 2-part 
   justification
    of Fact D.3 is 
   {\it much simpler}
   than a proof from first principles.
   %
   %This is  because the 
   %
   The
   latter would
   require examining 
    eight different 
   tableaux 
   elimination rules,
   as the detailed
   proofs of \cite{ww5}'s
    Lemmas 1 and 2 
   actually
   did do.}
   %
   %\noindent
   %-----------------------------------------------------------------}
for more details).  $~~\Box$

%%%% 
%%%% provides   some
%%%%    %some
%%%% %   precise 
%%%% % details 
%%%%    details about how 
%%%%    the   Lemmas 1 and 2
%%%% of \cite{ww5} do
%%%%  imply Fact D.3).
%%%% 

\nop

We will now apply Fact D.3 to prove 
Theorem D.\ref{D.4-theorx}.
Its invariant will, 
interestingly,
collapse entirely
\footnote{\h55 $~$The
difficulty posed by 
%the 
multiplication 
%operation
can be easily understood when one compares
two 
integers sequences  
   $ x_0, x_1, x_2,  ... $ 
and  $ y_0, y_1, y_2,  ... $,
defined 
as follows:
\begin{center}
$x_{i}~~=~~x_{i-1}+x_{i-1}~~~~~$
$~~$ AND  $~~$ 
$~~~~~y_{i}~~=~~y_{i-1}*y_{i-1}$
\end{center}
It turns out that the faster growth rate of  multiplication
under
the series  $~ y_0, y_1, y_2,  ...~ $
enables
one to to construct  tiny 
$Z-$Based  deduction trees
$\, T$ that
violate
the analog of Fact D.3 's paradigm.
(This is 
because 
such trees
 can have
G\"{o}del
numbers 
smaller than
 $ \, (a/b)^4 \,$,
while 
all their 
root-to-leaf branches 
can be simultaneously   
``closed'' via contradictions.)
This 
property of
            multiplication is 
analogous
to Example \ref{ex3-3}'s
observations about
          how        the
differing growth rates of  
$ x_0, x_1, x_2,  ... $ 
and  $ y_0, y_1, y_2,  ... $
are related to the
threshold where the semantic tableaux version of 
 Second Incompleteness Theorem can be evaded. \fend }
%% 
%% It
%% explains
%% intuitively 
%%  why   Fact D.3 and 
%% Theorem D.\ref{D.4-theorx} 
%% use
%% paradigms
%% that recognize addition but not multiplication
%% as total functions
%% during their evasions of
%% the semantic tableaux version of
%% the Second Incompleteness Theorem.  \fend }
%% 
 if one were to
merely  add
a multiplication function symbol to 
the
    U-Grounding language.
This is why
our 
boundary-case exceptions
to the semantic tableaux version of the Second Incompleteness
allow 
a Type-A axiom system 
to
recognize addition
as a total function (but
suppress a similar 
 treatment of
 multiplication).

\hgskip

{\bf Theorem
D.\ref{D.4-theorx}.}
{\it The generic configuration $\xi^*$ 
is both A-stable and
E-stable.} {\rm (This 
implies  many different 
self-justifying formalisms
exist
via Theorems 
\ref{ppp1}, \ref{pqq3}, \ref{pqq4}, \ref{pqq5},
 \ref{ppp6},
 G.2 and G.3.)}

\cvs

\hgskip

Our proof of
Theorem D.\ref{D.4-theorx} will 
separately
show 
 $\xi^*$  is 
A-stable and E-stable. 

\hgskip 

\parskip 3pt

{\bf Proof of 
 $ \, \xi^* \,$'s A-stability :}
Suppose for the sake of establishing a
proof by
 contradiction
that  $ \, \xi^* \, $ was not A-stable.
Then the constraint $ \, * \, $ 
of Definition \ref{astab}
would  be violated
by at least
%%%%%  bbbbb
some
$\theta \, \in \,$RE-Class$(\xi)$.
This violation
 will cause the statement $ \, + \, $ to be true
for 
such a
 $\theta~$:
\begin{description}
\item[  +   ]
There exists a
semantic tableaux
 proof $~p~$
of a
 $~\Pi_1^*~ $ theorem,
called say 
$ \, \Uxp ~,~ $ 
from
 the
axiom system of $~\theta \cup B^\xi~$
such that
 Log$(p)~\leq \zzz \theta)~+1~ $
and where 
 $ \, \Uxp \, $ also  
fails to
satisfy 
 Good$\{~ \, \tftt  \zzz \theta) \, ~\}~$.
\end{description}
Let us recall that
if  $\Uxp$ is $\Pi_1^*\,$ then
Reverse$(\, \Uxp \,)~$ 
is a 
$\Sigma_1^*$ sentence equivalent to $~\neg \Uxp~$.
Thus, 
Reverse$(\, \Uxp \,)~$ will satisfy
Good$\{~ \, \tftt  \zzz \theta) \, ~\}~$
criteria (simply because
it
% Reverse$(\, \Uxp \,)~$
  has the opposite 
goodness
property as
$\,\Uxp \,$ ).
Also, if
$~Z~$ 
denotes the axiom system of
 $~\theta \cup B^\xi \, + \, $Reverse$(\, \Uxp \,)$,
it is easy to verify
\footnote{\g55    \label{jnorm}
The axiom system
 $Z$  must satisfy 
\nor1
because:
\bee
\item
The
quantity
$~2^{ \zzz \theta  \, )}$
is a valid
first component
for $Z$'s   norming constraint 
because
all the axioms of $~B^\xi~$ are 
true in the Standard-M model
and
because
Definition \ref{chg}
implies
all of
$ \theta \,$'s axioms satisfy 
Good $\{\,\zzz \theta) \,\}$.
\item
The
quantity      $ ~\sqrt{~2^{ \zzz \theta \,   )}} ~$
is a valid
second component
for $Z$'s   norming constraint 
because
 Reverse$(\, \Uxp \,)~$ 
is the only $\Sigma_1^*$ sentence belonging to
Z, and  because  Reverse$(\, \Uxp \,)~$ 
satisfies 
 Good$\{~ \, \tftt  \zzz \theta) \, ~\}~$.
\ene}
that $~Z\,$'s axioms 
will
satisfy the 
\nor1
criteria.

\smallskip

It is next helpful to  observe that what is
a proof from one perspective corresponds to
being a
deduction tree from 
a different
perspective.
Thus, Item $~\, + \,$'s proof
 $~p~$
of the 
theorem
$ \, \Uxp \,$ 
from
 the
axiom system of $~\theta \cup B^\xi~$
corresponds to being  a Z-based deduction tree,
with $~Z~$ representing the
axiom system of
$~\theta \cup B^\xi \, + \, $Reverse$(\, \Uxp \,)$.
In this context,
Item  $\, + \,$'s 
inequality of 
 Log$(p)~\leq \zzz \theta  \, )~+1~ $
implies 
\footnote{ \g55 \label{fd7} 
Without loss of generality, we 
may assume that 
every non-trivial proof $~p~$ satisfies
  Log$(p)~\geq ~64~$ (since a 
string with fewer than 64 bits is too short to be a
proof).
Then the footnoted paragraph's 
 Log$(p)~\leq \zzz \theta)~+1~ $ inequality 
trivially implies 
$~p~<~3^{ \zzz \theta  \, )}~$.
In a context where 
 $~Z~$ 
is a
\nor1
axiom system,
the
 latter inequality 
certainly
implies
$~p~$, viewed as a deduction tree for $~Z~,~$
has a small enough G\"{o}del number to
 satisfy the
hypothesis for Fact D.3. 
(This is because if one sets
$~a~=~2^{ \zzz \theta  \, )} ~$ 
and 
$~b~=~\sqrt{2^{\zzz \theta  \, )}} ~$ 
then obviously 
$ ~~p~<~3^{ \zzz \theta  \, )}~ <~4^{ \zzz \theta  \, )}~ =~
(a/b)^4~~~).$  \fend } 
$\,$that 
$~p~$, viewed as a deduction tree for $~Z~,~$
 satisfies the
hypothesis of Fact D.3 . 
Hence, Fact D.3 establishes
that $~p~$ must contain
at least
 one contradiction-free
root-to-leaf branch.

\smallskip

This last observation 
is all that is needed
to confirm
 $~\xi^* \,$'s A-stability,
via a proof-by-contradiction.
This is because the 
definition of a semantic tableaux proof implies
every 
one of its
root-to-leaf branches 
must
end with a pair of
contradicting nodes.
However, the last paragraph showed 
$~p~$
will not satisfy this
required
 property,
if $~\xi^*~$ is not A-stable.
Hence 
our construction has proven the A-stability of
 $~\xi^*~$ by showing that otherwise an infeasible
circumstance will arise.
 $~~\Box$

\hgskip

{\bf  Proof of 
 $  \xi^*  $'s E-stability :}
A
proof-by-contradiction
will verify
 $  \xi^*  $ is
E-stable, analogous
 to  the
 proof of
its
    A-stability.
Thus if
 $  \xi^*  $ was not
E-stable, then
statement $++$ would be true
for some  $\theta$. (This is because 
at least one
$\theta \, \in \,$RE-Class$(\xi)$
would then
 violate
Definition \ref{estab}'s
requirement of  $~**~~$.)
\begin{description}
\item[ ++   ]
There exists a
semantic tableaux
 proof $~p~$
of a $\Sigma_1^\xi$ theorem 
$ \, \Uxp \, $ from the
axiom 
system $\theta \cup B^\xi$
such that
 Log$(p) \,\leq \zzz \theta) +1 $
and  
 $  \Uxp  $  also
fails to
satisfy 
Good$\{  \tftt  \zzz \theta) \,\}.$
\end{description}
Item $++  $ 
 implies 
Reverse$( \Uxp )$ 
satisfies
 Good$\{  \tftt  \zzz \theta) \}$
(because Reverse$( \Uxp )$ again
   has
 the opposite 
goodness property
as 
$ \Uxp $ ).
Let
$Z$ now
denote the formal axiom system of  
 $\, \theta \cup B^\xi  \, +  \, $Reverse$( \Uxp )$.
The footnote \footnote{ \sm55   
The axiom system
 $~Z~$ must satisfy 
\xor2
because:
\bee
\item The first component of its norming constraint 
can be set equal to
 $ ~\sqrt{~2^{ \zzz \theta   \, )}} ~$
because Reverse$(\, \Uxp \,)~$ 
is a 
 Good$\{~ \, \tftt  \zzz \theta  \, ) \, ~\}~$
$\Pi_1^*$ sentence, and all $~Z\,$'s other 
$\Pi_1^*$ sentences satisfy more relaxed constraints.
\item The second component of $Z$'s norming constraint 
is satisfied by the constant of
2 because Definition D.\ref{D2-def} implies  this 
quantity is always permissible
when $~Z~$ 
contains no $\Sigma_1^*$ axiom sentences.
\ene   \fend }
then
   uses
reasoning 
similar
to footnote    \ref{jnorm} 
to
show 
$Z$ 
satisfies 
\xor2

\parskip 3pt

\smallskip

As before 
via a simple change in notation, 
 $~p\,$'s
semantic tableaux proof of
$ \, \Uxp \,$ 
can be viewed 
as 
 a deduction tree
using $~Z\,$'s 
axioms.
Also as before, we may 
use the combination of the facts that 
 $~Z~$ is a
\xor2 system
and that  Item ++ 
indicated 
 Log$(p) \,\leq \zzz \theta) +1 $ to 
deduce\footnote{\sm55 The proof that 
 $~p~$ is
 small enough to satisfy 
Fact D.3 's
hypothesis in the current E-stable case 
is almost
identical to Footnote  \ref{fd7}'s analysis of
the A-stable case. 
Thus as in the earlier case,
Item ++'s inequality of
 Log$(p)~\leq \zzz \theta)~+1~ $ 
trivially implies 
$~p~<~3^{ \zzz \theta  \, )}~$.
Also, we may again assume that
  Log$(p)~\geq ~64~$ (since a sequence with fewer than 64 
bits cannot amount to a proof of any interesting fact under
all normal coding conventions). 
An analog of Footnote  \ref{fd7}'s chain of inequalities
will then allow us to conclude that 
$~p~$
is
 small enough 
proof from a 
\xor2
system
to
 satisfy the
hypothesis for 
Fact D.3.  \fend }  that
 $~p~$ is small enough to satisfy 
Fact D.3 's
hypothesis.
 Hence once again,
Fact D.3  implies that
$~Z~$ must contain
at least
 one contradiction-free
root-to-leaf branch.
As before, the existence 
 of this
contradiction-free
path 
violates the definition of a semantic tableaux proof
and 
enables
 our proof-by-contradiction to
reach its
desired end. $~\Box$ 

\hgskip

{\bf Remark D.5} {\it  (about 
Theorem  D.\ref{D.4-theorx}'s
significance) :}
Part-ii of
 Definition
D.\ref{D1-def}
indicated  $~\xi^* \,$'s base axiom of
$~B^*~$ was a Type-A formalism that recognized addition
as a total function. This is significant because
 \cite{ww0,ww2,ww7,ww9} 
showed
nearly all Type-M formalisms, including 
all the common axiomatizations for
I$\Sigma_0$,
are unable to recognize their 
semantic tableaux consistency.
Thus, the declaration that multiplication is a total
function
is {\it the trigger-point} causing
\footnote{\sm55 
We formally proved in \cite{ww0,ww2,ww7,ww9} 
that
multiplication's totality property
causes
 the semantic tableaux version of the Second
Incompleteness Theorem to become active.
The  
 Example \ref{ex3-3} 
summarizes the main
%%%  
%%%  offers a 
%%%  % nice brief 
%%%  summary of the
%%%  % underlying
%%%  
 intuition behind these
results.}
the semantic tableaux version of the
Second Incompleteness Theorem
to become active. 
This threshold effect is 
significant
% quite tight 
because
Theorem D.4, combined with
Theorems 
 \ref{pqq4}, \ref{pqq5}, G.2 and G.3,
formalize       {\it four different respects} in which 
Type-A
self-justifying 
formalisms can
prove all 
Peano Arithmetic's
 $\Pi_1^*$ theorems
{\it ( after} multiplication's totality axiom is suppressed).

\cvt

\subsection*{D-3. $~$Three Further Examples of
Stable Generic Configurations}

Our second 
example
of an
 EA-stable
        configuration
is
called $~\xi^{**}~$. 
It will be identical to 
 $~\xi^*~$ except that it will replace
semantic tableaux with a stronger
deduction method,
which 
\cite{ww5} 
    called  Tab$-U_1^*$.
The latter 
is          a
revised version of 
semantic tableaux
that permits
 a modus ponens
rule to perform
deductive cut 
 operations on
$\Pi_1^*$ 
and $\Sigma_1^*$ sentences.
(The formal
definition of 
 Tab$-U_1^*$ 
deduction 
had appeared in \cite{ww5}.
It will  be
unnecessary to
 repeat here.)

The Section 5.3 of 
\cite{ww5} 
noted
Tab$-U_1^*$ 
has
 similar self-justification properties
as conventional semantic tableaux.
%It is thus unnecessary to discuss
%Tab$-U_1^*$ in detail here.
All the results that
Section D-2 
proved about
 $~\xi^{*}~$
% ,
% however,
apply also  to  $~\xi^{**}~$,
via their natural generalization
under
\cite{ww5}'s
Tab$-U_1^*$
deduction
method.
Thus,
 $~\xi^{**}~$
is also EA-stable.

%%%%% 
%%%%% 
%%%%% Our second 
%%%%% example
%%%%% of an
%%%%%  EA-stable
%%%%%         configuration
%%%%% is
%%%%% called $~\xi^{**}~$. 
%%%%% It will be identical to 
%%%%%  $~\xi^*~$ except that it will replace
%%%%% semantic tableaux with a stronger
%%%%% deduction method,
%%%%% which 
%%%%% \cite{ww5} 
%%%%%     called  Tab$-U_1^*$.
%%%%% The latter 
%%%%% is          a
%%%%% revised version of 
%%%%% semantic tableaux
%%%%% that permits
%%%%%  a modus ponens
%%%%% rule to perform
%%%%% deductive cut 
%%%%%  operations on
%%%%% $\Pi_1^*$ 
%%%%% and $\Sigma_1^*$ sentences.
%%%%% (The 
%%%%% definition of 
%%%%%  Tab$-U_1^*$ 
%%%%% deduction 
%%%%% appeared in \cite{ww5}.
%%%%% It
%%%%% is unimportant to
%%%%%  repeat here.)
%%%%% 
%%%%% 
%%%%% 
%%%%% 
%%%%% The Section 5.3 of 
%%%%% \cite{ww5} 
%%%%% noted
%%%%% Tab$-U_1^*$ 
%%%%% has
%%%%%  similar self-justification properties
%%%%% as conventional semantic tableaux.
%%%%% It is thus unnecessary to discuss
%%%%% Tab$-U_1^*$ in detail here.
%%%%% %%  
%%%%% %% 
%%%%% %% We will not
%%%%% %%  go  into the full details
%%%%% %% again,
%%%%% %% for the sake of brevity.
%%%%% %% 
%%%%% %% 
%%%%% All the results that
%%%%% Section D-2 had
%%%%% proved about
%%%%%  $~\xi^{*}~$
%%%%% % ,
%%%%% % however,
%%%%% apply also  to  $~\xi^{**}~$,
%%%%% via their natural generalization
%%%%% under
%%%%% \cite{ww5}'s
%%%%% Tab$-U_1^*$
%%%%% deduction
%%%%% method.
%%%%% Thus,
%%%%%  $~\xi^{**}~$
%%%%% is also EA-stable.

A key point is that 
 there is a
non-trivial
distinction between  
$  \xi^{*}  $ 
and  $  \xi^{**}  $, despite
the fact that they
have 
similar technical qualities.
This is because 
 $  \xi^{**}  $ contains  a
Level-1
modus ponens rule
(unlike $  \xi^{*}  $ ).
 If 
it were infeasible to expand  $  \xi^{*}  $ 
into 
a broader 
 $  \xi^{**}  $, $  $then both formalisms could,
perhaps,
 be easily
dismissed as having 
negligible pragmatic 
significance (since
modus ponens is central
to
cogitation). 
However in a context where  $  \xi^{**}  $ does permit
a  Level-1
modus ponens rule, 
it is a tempting formalism
(despite its limited modus ponens rule).

Unlike
 $~\xi^{*}~$ 
and  $~\xi^{**}~$,
our third 
example of an EA-stable
configuration,
called $~\xi^{-}~$,
will support            an unlimited modus ponens
rule.
This will be possible because 
 $~\xi^-\,$'s
language of $~L^-~$ 
will 
be weaker than the languages of
 $~\xi^{*}~$ 
and  $~\xi^{**}~$.
Thus $~L^-$  will
include 
 the six Grounding functions, but
not the Growth functions of addition and doubling.
It will thus treat addition and multiplication as
3-way atomic predicates,
 Add$(x,y,z)$ 
and Mult$(x,y,z)$, rather than as
formal functions.

\cvt
\parskip 2pt

This
perspective
enabled  $~\xi^-~$   to
support an evasion of 
the Second Incompleteness Theorem with
an unlimited  modus ponens rule
present,
in a context where
the other four parts of 
its
generic configuration are defined below:
\bee
\item
The 
  $~\Delta_0^{-}~$
class for $~\xi^-~$
will 
be built in an
essentially natural
 manner from
the Grounding function
set.
It will 
thus 
include
all formulae in $~L^- \,$'s language, 
whose 
quantifiers
are bounded
 in any arbitrary manner 
using
the Grounding function
primitives.
\item
The
base axiom system
$~B^-~$ of $~\xi^-~$
 will
employ  an infinite number of constant symbols,
denoted as 
$    K_1   ,   K_2   ,   K_3   , \,   ...   $ 
where 
 $K_1=1$ 
and where 
 $K_{i+1}$ is
a power of 2
 defined by the axiom of:
\beq
\label{addc}
\mbox{Add}(~K_{i}~,~K_{i}~,~K_{i+1}~)
\enq
Thus, the combination of 
$    K_1   ,   K_2   ,   K_3   , \,   ...   $ 
with the
Grounding function  
of subtraction allows
the language  $L^-$
to encode the value of any
arbitrary natural
number (as 
Part 1 
of Definition \ref{def3.3} 
had required).
Essentially,
$~\xi^- \,$'s 
base axiom system
of $~B^-~$ 
can
be any 
consistent
r.e. set of
$~\Pi_1^-~$ sentences 
that includes \eq{addc}'s axiom schema
and is
able to
prove every 
 $\Delta_0^{-}~$ sentence which is valid
in the Standard-M model.
\item
$~\xi^- \, $'s
deduction method can be any version of a classic
Hilbert-style proof methodology. (Thus, it  will  include
a modus ponens rule with no restrictions.)
\item
$~\xi^- \, $'s
G\"{o}delization method can be essentially
any natural technique.
\ene
An interesting aspect of  $\xi^-$
is it can be 
proven to be
EA-stable
via an analog of
Section D-2's treatment of  $\xi^*$.
Thus,
Theorem \ref{pqq4}
implies
every axiom system
$\alpha$,
whose
 $\Pi_1^-$ theorems 
hold true
in the Standard-M model,
can be mapped
onto an extension of  $\xi^-\,$'s base
axiom system
that can 
 recognize
its own 
Hilbert
 consistency 
and
prove 
$\alpha $'s  $\Pi_1^-$ theorems.
Except for
minor changes in notation,
this  
result
represents a new way of proving
\cite{wwapal}'s Theorem$\, 3.~$

\medskip

The self-justifying features
of
 $ \xi^{*} $,  $ \xi^{**} $  and  $ \xi^{-} $
are of interest
primarily 
 because
the 
Second 
Incompleteness Theorem implies 
that they
cannot be 
improved
much 
further.
This tight fit is
summarized by  Items 1-4.
\begin{enumerate}
\topsep -15pt
\item 
The Theorem 2.1 
(due to  the combined work of 
Nelson,
Pudl\'{a}k,
Solovay and Wilkie-Paris
\cite{Ne86,Pu85,So94,WP87} )
implies no natural axiom system can
prove Successor is a total function and recognize its own
Hilbert consistency. This theorem 
thus explains why
the presence of
growth functions 
must be omitted from
 $\xi^{-}\,$'s
base axiom system of $B^-$.
\item
Moreover, \cite{wwapal}
  proved 
$~\xi^{-} \,$'s method 
for evading  
the Second Incompleteness
Theorem
will 
collapse
if one replaces  
\ep{addc}'s  ``addition-based named sequence'' of
constant symbols
$    K_1   ,   K_2   ,   K_3   , \,   ...   $ 
with a faster 
growing
``multiplicative 
convention'',
where 
the
constant symbols 
$    C_1   ,   C_2   ,   C_3   , \,   ...   $
are formally defined via 
\eq{multc}'s
 schema.
\beq
\label{multc}
\mbox{Mult}(~C_{i}~,~C_{i}~,~C_{i+1}~)
\enq
Thus,  \cite{wwapal} showed that there
exists a 
 $\Pi_1^{-}$ sentence $ \,W \,$ 
(provable from Peano Arithmetic)
such that no consistent  system can
simultaneously prove  $W$,
contain \eq{multc}'s axiom schema and prove the non-existence
of proof of $0=1$ from itself.
There is no space to prove it here, but 
a generalization of the Second Incompleteness
Theorem implies  the modification
of $\xi^-$ that replaces
\eq{addc}'s  axiom schema with 
\eq{multc}'s schema 
{\it is not even 0-stable.}
\item
Similarly,
\cite{ww2,ww7} proved that if 
 $ \, \xi^{*} \, $'s and  $ \, \xi^{**} \, $'s 
 base 
 axiom system of 
$ \, B^* \, $ 
was strengthened
to include the assumption that
multiplication was a total function then \cite{ww5}'s
two
semantic tableaux  evasions of
the Second Incompleteness Theorem would 
both collapse.
\item
Also, \cite{wwlogos} 
proved
that  an analog of $\xi^{**} \,$'s
evasion of the Second Incompleteness Theorem
will collapse if its 
modus ponens rule was expanded to apply to
either $\Pi_2^*$ or $\Sigma_2^*$ sentences.
\ene
The Item 3
is especially interesting because
\cite{ww6} 
proved \cite{ww5}'s evasion of the Second Incompleteness
Theorem 
was compatible with its formalism recognizing
an infinitized generalization of a
computer's floating point multiplication
as a total function.
Thus
 while the semantic tableaux formalisms
of 
 $\xi^{*}$ or  $\xi^{**}$ 
are provably unable \cite{ww2,ww7} to
 recognize integer
multiplication as a total function,
their relationship to floating point multiplication
is 
more subtle.

\smallskip

Our fourth example of an
application of Section \ref{sect64}'s theorems
was stimulated by
some
insightful 
email  we
received from
L. A. Ko{\l}odziejczyk
\cite{Ko5} in 2005.
It noted
there existed a
potential
 exponential gap between the lengths
of semantic tableaux and Herbrand-style proofs 
under some 
circumstances.
Our earlier research \cite{ww2} 
addressed a 
1981 Paris-Wilkie open question \cite{PW81} 
by
generalizing some 
Adamowicz-Zbierski techniques 
\cite{Ad2,AZ1} 
to show
a natural  axiomatization of I$\Sigma_0$ 
 satisfied the semantic tableaux version of the Second
Incompleteness Theorem.
In this context,
Ko{\l}odziejczyk
asked whether this would
apply to       all
plausible axiomatizations for I$\Sigma_0$ ?

We replied in
\cite{ww9}
to
Ko{\l}odziejczyk's
stimulating question
 by 
distinguishing between Example \ref{ex3-1}'s 
$\Delta_0^A$ 
and $\Delta_0^R$  formulae
and by using
the Paris-Dimitracopoulos \cite{PD82}
translation algorithm for 
$\Delta_0$ formulae. (The latter procedure was summarized
earlier by
 Lemma \ref{lex22}. It
demonstrated
 how to
 map
classic arithmetic's
$\Delta_0^A$ formulae
onto  equivalent   
$\Delta_0^R$  formulae
in  the Standard-M
model.) Our reply to 
Ko{\l}odziejczyk's
question, 
thus,
employed this
translation methodology to show that there existed an axiom system,
called Ax-3, which proved the identical set of theorems as
the more common Ax-1 and Ax-2 encodings of 
 I$\Sigma_0$ and which possessed the following 
pair of 
quite
fascinating 
  contrasting properties:
\bed
\topsep -3pt
\item[   A  ]
No
consistent
 superset $~\beta~$ of Ax-3's set of axioms is capable of
proving its
own
 semantic tableaux consistency \cite{ww9}.
\item[   B  ]
In contrast,  
if  ``Herb'' denotes
the next paragraph's
 Herbrand-styled deduction
and if ``SelfRef'' denotes 
the sentence $~\bullet~$ from
Section \ref{secc1},
then 
Ax3$\, + \,$SelfRef(Ax-3,Herb) will
be a self-justifying axiom system.
%%% 
%%% 
%%% In contrast,  
%%% if  ``Herb'' denotes
%%% the next paragraph's
%%%  Herbrand-styled deduction
%%% method
%%% then 
%%% Ax3$\, + \,$SelfRef(Ax-3,Herb) will
%%% satisfy both Parts (i) and (ii) of 
%%% Section \ref{secc1}'s 
%%%  definition of
%%% self-justification.
%%% 
%%% 
\ennd

The intuition behind \cite{ww9}'s proof of Items A and B can be 
easily summarized if 
we define a
 ``Herbrandized-style''  proof of a 
theorem $ \,\Phi \,$ from an axiom system $ \, \alpha \,$
as being an essentially 
2-part
structure where:
\bee
\topsep -3pt
\item Each  of $\alpha$'s axioms and
also the sentence $~ \neg \Phi~$ are
first written as  Skolemized
expressions.
\item
A propositional calculus proof 
is then used to show
 that some formal conjunction
of instances of Item 1's Skolemization schema has no satisfying
truth assignment.
\ene
Such a
 formalism is 
different from the definition
of a semantic tableaux proof (appearing in for example Fitting's
textbook \cite{Fi90} ).
This is because the latter replaces the use
of Skolemization in Items 1 and 2 with an existential quantifier
elimination rule. It turns out that this distinction enables
some semantic tableaux proofs to be exponentially more compressed
than their Herbrandized counterparts,
as Ko{\l}odziejczyk observed 
\cite{Ko5,Ko6}. This fact
 enabled \cite{ww9}
to prove that Herbrandized and semantic tableaux proofs
have the divergent properties summarized by Items A and B.

One reason
Ax-3's evasion of the Second Incompleteness Theorem
is of interest is that I$\Sigma_0$ supports many more generalizations
of the Second Incompleteness Theorem
than evasions 
of it.
 Thus,
Willard  \cite{ww2,ww7,ww9}
proved that the semantic tableaux version of
the Second Incompleteness Theorem 
was valid for three different encodings of
I$\Sigma_0$, and
Adamowicz, Salehi and  Zbierski
have discussed in great detail
\cite{Ad2,AZ1,Sa11}
 various Herbrandized generalizations
of the Second Incompleteness Theorem for
particular encodings of  
I$\Sigma_0$ 
and I$\Sigma_0+\Omega_i$. 
Moreover, an added 
facet of \cite{ww9}'s 
Ax-3 encoding for I$\Sigma_0$ is that most automated
theorem provers use a 
particular
variant of the Resolution method
that causes \cite{ww9}'s 
unusual 
methodology
to apply
also to them \footnote{
\b55 
The main theorems in \cite{ww9} generalize for resolution because
Resolution-based theorem provers
employ skolemization analogously to Herbrand deduction. \fend }.

The reason for our 
%particular 
interest in
\cite{ww9}'s results 
% in the current article
is that it represents a fourth example 
%% of
where the 
 meta-theorems from Sections 
\ref{3uuuu3} and  \ref{sect64} 
can be useful.
Thus,
the footnote
\footnote{\b55  The 
discussion in \cite{ww9}
 did not
technically  
use
 Definition  \ref{estab}'s machinery
to establish 
%that 
there existed an extension of
its ``Ax-3'' encoding for 
I$\Sigma_0$ that could recognize its own Herbrand consistency.
Its formalism,
however,
 could be
easily couched in terms of 
 Definition  \ref{estab}'s machinery, if one uses a 
generic 
configuration $~\xi^R~$ where
\bee
\item  $~\xi^R\,$'s 
base language is the same as the usual language of
arithmetic, 
\item  $~\xi^R\,$'s 
$~\Delta_0^R~$ sub-class is defined by
Item (b) in Example \ref{ex3-1}, 
\item  $~\xi^R\,$'s 
base axiom system 
is \cite{ww9}'s ``Ax-3'' system, 
\item  $~\xi^R\,$'s 
deduction method
is either a Herbrandized styled-method
or a Resolution system that relies upon
Skolemizatin in a similar manner.
\item  $~\xi^R\,$'s G\"{o}del encoding scheme
may be any such 
natural method.
\ene
This approach 
supports a 
% somewhat
  stronger form of
self-justification result
than
had 
appeared in \cite{ww9}.
This is
because $~\xi^R~$ can be proven to be E-stable (by a 
% straightforward 
generalization of \cite{ww9}'s analysis
techniques). 
Thus,
 \phx{ppp2} 
implies
that
Ax-3 has a well-defined self-justifying extension
that
can  recognizes its own
formalized
 Level$(0^R)$ consistency.
(This 
self-justification result
is  stronger than
\cite{ww9}'s main theorem.
The latter
merely established
that some extension of Ax-3
recognized the non-existence of a
Herbrandized deduction of $0=1$ from itself.) \fend }
summarizes how  a fourth 
type of
generic configuration,
called $~\xi^R~$,
can be defined
 that both duplicates
\cite{ww9}'s
main
 self-justification results 
under the above definition of
Herb-deduction,
as well as  strengthens
them. (In particular, 
  $~\xi^R~$ meets
Theorem \ref{ppp2}'s
requirements,
and  self-justifying 
extensions of its Ax-3 system thus
recognize their 
%own
 Level$(0^R)$ consistency.)

\parskip 2pt

The properties of our
four generic configurations of $~\xi^R$,
  $~\xi^*~$,
 $~\xi^{**}~$ and $~\xi^-~$ are  summarized by Table I.
These configurations are
listed in ascending order according to the strength
of their deduction methods $~d~$. As their deduction methods
increase in strength, 
these configurations 
% the associated configurations 
have
their ability
reduced
%weakened
 to recognize 
the totality of 
the addition and multiplication operations.
%% 
%% 
%% (This is because the corresponding generic configurations will
%% violate the requirements of the Second Incompleteness Theorem,
%% if they simultaneously use too strong a deduction method
%% and possess too strong an understanding of their
%% own consistency.)

$~\xi^R\,$ is thus a Type
Almost-M 
system
that
can prove multiplication is a total function
(but 
which
does not contain 
\ep{totdefsymbm}'s
 totality statement {\it as an axiom)}.
On the other hand, $~\xi^-~$ 
uses a stronger Hilbert-styled
 deduction methodology,
which is  incompatible with treating
the totality of addition or multiplication as either
axioms {\it or as derived theorems.} 
%% (This incompatibility is unavoidable because the
%% Theorem 2.1, due to the joint work of
%% Nelson,
%% Pudl\'{a}k,
%% Solovay and Wilkie-Paris
%% \cite{Ne86,Pu85,So94,WP87}, 
%% implies that self-justifying systems cannot simultaneously
%% recognize successor as a total function and prove their
%% own Hilbert consistency.)

Each of Table I's
rows
 $ \, \xi^R$,
 $ \, \xi^{**} \, $ and $ \, \xi^- \, $  
are  maximal (in that an alternate row
improves upon one   column's measurement 
only when it is 
weaker from the perspective of
another 
column).
Only  $ \, \xi^* \, $ 
is
an exception to this rule: It is
strictly
weaker
\footnote{\b55 \label{f28} $  \xi^{**}  $ employs a stronger
deduction method 
than   $  \xi^*  $
because it allows a modus ponens rule for $\Pi_1^*$ and
$\Sigma_1^*$ sentences to be added to semantic tableaux 
deduction
(see \cite{ww5} for the precise definition of this  
``Tab$-U_1^*~$''
 modification of the  semantic tableaux 
deductive method).  \fend }
     than
     $ \,  \xi^{**}  $. 
    This appendix has discussed 
      $ \, \xi^* \, $ 
%% despite its sub-optimality, 
%%    %primarily 
%%    partly 
because it makes
    Theorem D.\ref{D.4-theorx}'s proof simpler
    (and also because 
     semantic tableaux 
    % deduction 
    is a 
    % very
    frequent topic in the logic literature).

    %%%footnnnnnnnnn 

%\smallskip

\begin{center}
{\bf Table I}
\end{center}

\bigskip

\small
\baselineskip = 1.1 \normalbaselineskip 
\noindent
\begin{tabular}{|c|c|c|c|c|c|c|} \hline
Name & Deduction Method  & Type & Almost   & Type  & {\bf Axiom} 
& {Self-Just }  \\
 &     & {\bf A} 
& {\bf M} & {\bf M} &  { Format} &  { Level}
\\ \hline\hline 
 &  Resolution and/or  & &  &  &  & \\
 $\xi^R$ & Herbrandized analogs & Yes$^{35}$
  & Yes  & No &  E-stable &Level $(0^R)$ \\ \hline
 &  & &  &  &  & \\
 $\xi^*$ & Semantic Tableaux & Yes  & No  & No & EA-stable & Level $(1^*)$  \\ \hline
 &  & &  &  &  & \\
 $\xi^{**}$ & Tab$-U_1^*$ Deduction$\, ^{34}$ & Yes  & No  & No &  EA-stable &Level $(1^*)$   \\ \hline
 &  & &  &  &  & \\
 $\xi^{-}$ & Hilbert Deduction & No  & No  & No &  EA-stable &
 Level $(\, \infty^{-} \,)$   \\ \hline
\end{tabular}

%\njp

\bigskip
\normalsize
\cvl

%\cvt

\parskip 0pt

%\nskip
%\nskip

\bigskip
The footnote
\footnote{\b55    For the sake of simplicity,
    the Ax-3 system of \cite{ww9}  did not
    use either 
    Equations
    \eq{totdefxa} or  
    \eq{totdefsymba}'s
  as axiom statements
(since  they were
provable as  theorems).
 All
    \cite{ww9}'s 
     results 
    do,
however,  generalize
%    readily
    when 
\eq{totdefxa}'s
statement  about
addition's totality 
is included as 
%added 
an axiom.
Thus, it is appropriate to attach the designation
of ``Yes'' with a caveat to  the ``Type-A'' entry
in Table I's first row. (This row is called 
``Resolution and/or Herbrandized analogs'' 
because it applies
to 
essentially 
any deduction scheme that relies upon Skolemization
as an alternative to  \cite{Fi90}'s
semantic tableaux
existential quantifier elimination rule.)} ,
 attached to Table I's first row,
explains why  a caveat is attached to its first
``Yes'' entry. 
The theme of Table I is 
%thus 
that self-justifying
axiom systems have some nice redeeming features,
although
the
 Second
Incompleteness Theorem
clearly also
imposes
severe
 limits on their abilities.
This point  will be reinforced 
when
Appendix E introduces a generalization of the
Second Incompleteness Theorem,
that shows 
\thx{ppp6}'s translational reflection principle is close
to being a maximal feasible result, and when Appendix F 
discusses the epistemological  significance
of self justification.

%% 
%% the next two sections discuss
%% hybridizations of
%% Theorems  \ref{pqq4}
%% and \ref{pqq5}
%% and the epistemological  implications
%% of self justifying systems.

%% 
%% an alternate system
%% and when 
%% Appendix F discusses the
%% epistemological  implications
%% of self justification.

%% 
%% the next two sections discuss
%% hybridizations of
%% Theorems  \ref{pqq4}
%% and \ref{pqq5}
%% and the epistemological  implications
%% of self justifying systems.

\nskip
\nskip

%% 
%% an alternate system
%% and when 
%% Appendix F discusses the
%% epistemological  implications
%% of self justification.

\nop

\section*{Appendix E: A Clarification of
Theorem \ref{ppp6}'s Significance}

% An Anti-Reflection Principle for
%\LARGE
% \baselineskip = 1.8 \normalbaselineskip 

It has been known since the time of G\"{o}del that
most conventional 
arithmetic 
axiom systems will satisfy the following two 
invariants:
\bee
\item
They are {\it physically unable} to prove
their own consistency
\item
They are  $\Sigma_1$ 
complete. This means they can
formally prove any  $\Sigma_1$
arithmetic sentence that holds
true in the Standard-M model, and
they can likewise refute any
$\Pi_1^\xi$ sentence 
that is false.
\ene
Let $~\xi~$ denote any generic configuration of the
form \gggcp . This appendix will use the
term {\bf $\xi -$Conventional} to describe any axiom
system that satisfies analogs of the
preceding 
%two 
conditions for generic configurations.
Thus
$~\alpha~$ is
$\xi -$Conventional 
iff it satisfies the following two
criteria:
\bed
\item[ a.  ]
The axiom system $~\alpha~$
will be {\it  unable} to verify
its own consistency under $~\xi\,$'s
deduction
    method  of $~d~$.
\item[ b.  ]
The axiom system $~\alpha~$
will be an extension of $~\xi \,$'s base axiom of $~B^\xi~$.
Part-3 of Definition \ref{def3.3} will thus imply it
is 
$\Sigma^\xi_1$ 
complete. (Hence, $~\alpha~$ can
formally prove any  $\Sigma^\xi_1$ sentence that holds
true in the Standard-M model, and
it can likewise refute any
$\Pi^\xi_1$ sentence 
that is false.)
\ennd

This section will prove no analog of \ep{T-reflect}'s
translational reflection principle is feasible for 
$\xi -$Conventional axiom systems. Thus,
Theorem \ref{ppp6} must be close to being a maximal result,
since it cannot 
plausibly
be further extended to hold under conventional axiom systems.

\medskip

{\bf Theorem E.1}  (A New Type of Version of the Second Incompleteness
Theorem): 
{\it
There exists no 
$\xi -$Conventional 
axiom
system $~\alpha~$ 
that can prove the validity of 
\eq{G-reflect}'s
Translational Reflection Principle
for any translation-mapping T.}
(In other words,
there exists no
 algorithm  $~T~$ that 
maps 
$~\Pi_1^\xi$ sentences $~\Psi~$  
onto alternate  $~\Pi_1^\xi$ 
sentences
$~\Psi^T~,~$ 
which 
 are equivalent 
to $~\Psi~$
in the
Standard-M model
and where $~\alpha~$ can
verify
\eq{G-reflect}'s reflection principle for every
$~\Pi_1^\xi$ 
 sentence $~\Psi .~)$
%
%simultaneously
%for 
% all  $~\Pi_1^\xi$ sentences
% $~\Psi~$.) 
\beq
\label{G-reflect}
\forall ~p~~~[~~ \mbox{Prf}_{\alpha,d}(~\lceil \, \Psi \, \rceil~,~p~)
  ~~~ \Rightarrow ~~~ \Psi^T~~]
\enq

{\bf Proof:} $~$It is  easy to prove Theorem E.1 via a
proof-by-contradiction. Thus consider the possibility
that Theorem E.1's translational mapping $~T~$ did exist.
One can then easily select a  $~\Pi_1^\xi$ sentences $~\Psi~$  
that is false in the 
Standard-M model. 
Then 
$~\Psi^T~$ is also false under the 
Standard-M model 
(since 
 $~\Psi~$  
and
$~\Psi^T~$ 
are equivalent in this model).

%% Moreover, $T$'s definition requires 
%% $\Psi^T$ 
%% to have a $\Pi_1^\xi$ format.
%% Thus,

Hence
Part-b of the definition of 
$\xi -$Conventionality implies
$~\alpha~$ must prove 
$~\neg~\Psi^T$
 (on account of $~\Psi^T\,$'s 
$\Pi_1^\xi$ format).

It is at this juncture that our proof-by-contradiction will
reach its end. This is because if $~\alpha~$ can prove
\eq{G-reflect}'s statement
and
{\it also prove}
 the sentence $~\neg~\Psi^T~$,
 then it
certainly can combine these two facts to prove the
non-existence of a proof of $~\Psi~$.
The latter 
contradicts 
Part-a of the definition of 
$\xi -$Conventionality (because it  shows $~\alpha~$
can verify its own consistency). $~~\Box$

%%  The latter would
%% contradict 
%% Part-a of the definition of 
%% $\xi -$Conventionality (because it would show $~\alpha~$
%% can verify its own consistency). $~~\Box$ 

\bigskip

{\bf Remark E.2.} 
We remind the reader that 
Footnote \ref{imper} 
pointed out that
$T$'s translational mapping 
would lose its main functionality, if 
it did not require $\Psi^T$ 
to have a  $~\Pi_1^\xi~$ format,
similar to $\Psi$.
In essence,
Theorem E.1 is
of interest
because it shows that 
Theorem \ref{ppp6}'s evasion of the Second Incompleteness Theorem 
is close to
being a maximal result. 
%(because 
(It thus shows that \eq{G-reflect}'s
translational reflection principle does not generalize 
to conventional 
axiom systems.)
%settings.)
This 
dichotomy
may explain why self-justifying axiom systems,
along with    Theorem \ref{ppp6}'s
particular
 invariant,
are 
potentially useful results.

%surprising topics.

 % {\bf Remark E.2.} 
% The 
% preceding proof 
% of Theorem E.1 makes it
% clear
% that its generalization of the
% Second Incompleteness Theorem
% is a 
% straightforward
% result. It is 

\cvl

\section*{Appendix F: Epistemological 
Perspective and Speculations}

\cvs

\parskip 6pt

It is desirable to include a short
purely
epistemological
discussion 
within this 
mostly mathematical article so that 
the more subtle nature of our
 results 
cannot be
misconstrued.

Part of the reason 
Self Justification
can lend itself 
to
easy
misinterpretations is that 
the
First Incompleteness Theorem 
demonstrates 
the
impossibility of 
constructing an ideally
optimal axiomatization of number theory. 
For
any initial
r.e$. \,$axiom system $~\alpha~$
and deduction
method$~d$, 
G\"{o}del 
thus
noted
it is
easy
\footnote{  \b55 
Let $\mho(a,d)$ the
classic
 G\"{o}del
sentence 
that  asserts:
{\it ``There is no proof of {\it this sentence} from 
$\alpha$'s
axiom system 
under $\,d$'s deduction method.''}
G\"{o}del
\cite{Go31}
 noted 
$~\alpha+\mho(\alpha,d)~$ 
always proves  more theorems than  $\alpha$.}
 to
develop
 an extension of  $~\alpha~$
that can prove 
strictly more
theorems than $~\alpha~$ under 
 $\,d$'s
deduction method. 
Moreover, a large number of generalizations of the Second Incompleteness
Theorem, starting with its 1939  
Hilbert-Bernays version
\cite{HB39}, are known to be 
%very 
robust results.

Such considerations 
naturally lead to  questions
about whether
any
 r.e$. \,$axiom system can encompass the workings
of the human mind. It may surprise some readers to learn
that this author 
shares 
such
 skepticism.
That is, 
we
doubt
{\it any single ISOLATED} self-justifying
r.e.$\,$logic 
can {\it  fully} approximate
the complex workings of the human mind.

In this short appendix,
let us instead
view  cogitation
as 
{\it roughly} a
process wondering
though some universe $ \,   \nal{U}$, 
comprised of
{\it both} consistent and
inconsistent axiom systems,
with a
trial-and-error evolutionary method 
focusing 
 its attention
over time  
increasingly  
onto
 the members of this universe 
 $~  \nal{U}~$ that are found to be
 consistent.
It is 
%
%%notre
%
%
%%notre-only straightforward
%
%easy
%
%%cccorn-ff  
%%cccorn-cc  
  straightforward\footnote{  \b55 It 
   is trivial from a theoretical perspective
   to design a 
   learning heuristic that
   will utilize all 
   consistent axiom systems
   from
   its available  universe $  \nal{U}$ 
    eventually, and
   it will
   spend only an infinitesimal fraction of its effort on
   inconsistent systems as time runs to infinity.
   (This because
   there exists only
   a countable number of distinct r.e. sets
   belonging to the universe
   $  \nal{U}.~)~$
   Also, 
   this
   learning
    process can 
   presumably be made to
    employ 
   some type of
   smart souped-up
   AI heuristics to enhance its efficiency,
   whose details will not concern us 
   % here
   in this abbreviated 3-page appendix.
   What is 
   central to the current discussion
    is that some type of formally
   {\it non-recursive} 
   and presumably trial-and-error
   method must 
   obviously
   be used 
   by this learning process
    to find
   the consistent elements of
   $~  \nal{U}~,~$ 
   on account of G\"{o}del's
   undecidability results.}
 to define
many 
universes  $~  \nal{U}~$  and
 evolutionary processes that 
fall into this gendre.
Our 
goal 
in this section
will be to
examine
 Section \ref{3uuuu3}'s
``R-View'' $\theta$ and its RE-Class$(\xi)$.

Thus,  $\theta$ will denote an 
R-View
that consists of
an arbitrary
r.e$. \,$set of
$\Pi_1^\xi$ sentences.
Also, 
RE-Class$(\xi)$ will
again
denote the
set of all
$~\theta~$ which 
can be built under $~\xi \,$'s language
of
$\,L^\xi$.
(Section \ref{3uuuu3}
had
allowed
 both valid and invalid
R-Views  $~\theta$ 
 to appear in
RE-Class$(\xi)$ because 
no recursive 
decision
procedure can 
identify
all  
the Standard-M model's true  $\Pi_1^\xi$ sentences.)

%% \nop

The epistemological purpose of this notation was revealed
in Section \ref{sect64}.
For the cases where $k  =  0$
or 1,
Section \ref{sect64}
defined 
$G^\xi_k(\, \theta \,)$
to be  the 
axiom system:
\begin{equation}
\label{f4gedef}
G^\xi_k(~ \theta ~)~~= ~~
\theta~\cup~B^\xi~\cup~\mbox{SelfCons}^k\{~[~\theta \,\cup \,B^\xi~]~,d~\}
\end{equation}
Also, Definition \ref{dap4-1} indicated that
  the function 
$ \, G^\xi_k ~ $
(which maps $ \, \theta \, $ onto 
$G^\xi_k(\, \theta \,)~~~)$ would be
called 
{\bf Consistency Preserving} iff 
$ \, G^\xi_k(\, \theta \,) \, $is 
assured to be consistent whenever 
all the sentences in $~\theta ~$
are  
true  under the Standard-M model.
\thx{pqq3} 
indicated,
in this context,
 that
$~G^\xi_1~$ 
satisfies 
this 
property
whenever $~\xi~$ is
EA-stable.
Likewise,
$~G^\xi_0~$ 
is      consistency preserving whenever $~\xi~$ is
one of A-stable, E-stable or 0-stable.

\nskip

These results indicate 
a trial-and-error experimental process
can, indeed, walk 
{\it in an 
unusually
orderly manner} through an universe
of self-reflecting 
candidate
formalisms, when
RE-Class$(\xi)$  denotes 
$~  \nal{U}\,$'s 
universe  
and
 $~\xi~$ satisfies any of
the  EA-stable, E-stable, 
A-stable or 0-stable conditions.
This is because if 
$~\theta~$ designates
a 
set of
$\Pi_1^\xi$ sentences
holding true 
in
 the Standard-M model, 
then  
$~G^\xi_k(\theta)~$ 
will
{\bf automatically}
satisfy both Parts (i) and (ii)
of Section 1's definition of Self Justification,
according to \thx{pqq3}.

Such 
consistency preservation is
surprising because 
it is simply inapplicable to
%does not apply to
the $\,G^\xi_k \,$ 
% functions of
 functions for
most 
pairs $~(\xi,k).$
\thx{pqq3}'s  first
 contribution 
is,
 thus,  
that it formalizes
how  $G^\xi_k\,$'s  mapping function 
can
represent
 a type of approximation
for
instinctive faith,
under certain well-defined circumstances.

This notion of instinctive faith is, of course, 
less
robust than a conventional proof.
One
obvious
  difficulty
is that a 1-sentence proof,
using an {\it ``I am consistent''} axiom,
is 
less convincing
than a full-length proof from first principles. Also, if
the initial formalism $~\theta~$ contains a false
$~\Pi_1^\xi~$ sentence then  $~B^\xi+\theta~$
and $~G^\xi_k(\theta)~$ 
will be  both inconsistent.

\nop

\cvl

Nevertheless
for $\,k\,$ equals 0 or 1,
 if   $~\theta~$ is comprised of the true sentences
in the Standard-M model, then 
\thx{pqq3} 
will
assure that
$~G^\xi_k(\theta)~$ is 
a consistent system that has an ability
to
use its {\it ``I am consistent''}
axiom sentence
to 
formalize
its
own consistency.
Moreover, the axiom system
$~G^\xi_k(\theta)~$ is helpful because
G\"{o}del's famous centennial paper
% has certainly 
implicitly
raised
the following
bedeviling 
issue:
%dilemma:
\begin{quote}
$\#~~$   How is it that Human Beings
 manage 
to muster 
the physical 
drive
to think (and prove theorems) when the many
generalizations of 
G\"{o}del's Second 
Incompleteness Theorem 
demonstrate
conventional logics
lack knowledge of
their own consistency?  
\end{quote}
While  philosophical paradoxes and ironical 
dilemmas,
similar to 
$~\# ~,~$ 
 never yield
perfect answers, the preceding discussion is helpful
because it 
explores
a
certain
syllogism
% paradigms
whereby a logic
can 
formalize
 at least some
fragmented
operational
 appreciation of its own consistency.

Moreover, 
Part-3 of Appendix D
 indicated that its 
four self-justifying configurations were
 close to being maximal results
that cannot  be 
much
improved,
on account of various
barriers imposed by
 the Second Incompleteness Theorem. 
Thus, these
particular
 positive results,
combined with Theorems \ref{ppp1}
\ref{pqq3},   \ref{pqq4},  
\ref{pqq5},
 \ref{ppp6},
D.4, E.1,
 G.2, G.3 and Remarks
\ref{re4-1} and
 \ref{recc1},  
come 
close to formalizing the
maximal variants of instinctive faith that a 
first-order logic can 
bolster.

The theme of the last two paragraphs
is
thus
 that our approximation of 
 {\it ``instinctive faith''} may be imperfect, but it
is still a useful partial reply to 
$~\# \,$'s puzzling
dilemma
{\it in a context where}
unambiguous 
full 
resolutions to $\, \# \,$
{\it are not permitted by}
the Second Incompleteness
Theorem.
Furthermore,
\ep{T-reflect}'s
translational reflection principle,
together with Theorem \ref{ppp6}
and the Remarks \ref{f88} and \ref{remhappy},
illustrate how the notion of 
 an instinctive faith 
about the usefulness of $\Pi_1^\xi$ theorems
 can
be
almost physically
{\it hard-wired} into self-justifying formalisms.

%%% instinctive faith 
%%% about the validity of $\Pi_1^\xi$ theorems
%%%  can
%%% be
%%% almost physically
%%% {\it hard-wired} into self-justifying formalisms.

%%cccorn-ff 

  \bigskip        
     {\bf A Yet Further
      Facet
       of this Unusual Epistemological Interpretation: } 
      Let
      % us use 
      the term 
      {\it Epistemological Bundle Theory} 
      refer to
       the underlying
      theory, advanced in this appendix, which 
      speculates about a 
       Thinking Agent
      % as 
      walking
      through   
      RE-Class$(\xi)$'s
      bundled universe of valid and invalid
      collections of $\Pi_1^\xi$ sentences
      and 
      then applying  some heuristic  to
      attempt to
      identify
      %locate
      % locate (via  heuristics)
       those
      % particular elements
       $\theta \, \in \,$RE-Class$(\xi)$
      whose sentences are true under the Standard-M model.
      
      Such a 
      theory has a second virtue, aside from 
      % its
      addressing  $ \# \,$'s
       paradoxical question
      about the nature  
      of {\it ``instinctive faith''. }
      It also clarifies
      the meaning of 
      our main theorems
      %%  
      %%   
      %%   and simplifies
      %%  the 
      %%  % mathematical
      %%   structure of 
      %%  Sections
      %%  \ref{3uuuu1}-\ref{sect64}'s theorems  
      %%  % formal
      %%  
      %%  
      and the related
           E-stability, A-stability,
       EA-stability and 
      RE-Class$(\xi)$ constructs.

      %%
      %%help 
      %%analyze 
      %%%prove their theorems about 
      %%self-justification.
      %%It turns out that
      %%%A nice aspect about the  
      %% epistemological bundling
      %%can explain the motivation behind these
      %%theorems.
      %%
      %%
      %%*is that it can explain much of the intuition behind
      %%* these theorems.

      This is because the
      Items $\, * \,$ and  $\, ** \,$
      from the
      definitions of 
      A-stability and E-stability
      in Section \ref{3uuuu3} 
      formalize
      how a thinking agent $~T~$ can view short
      proofs from a {\it technically inconsistent} axiom system
      of $~B^\xi \cup \theta~$ as containing
      % some 
      pragmatically
      useful information
      {\it under the assumption} that the
      lengths 
      %   FIXED ALREADY length's length's length's
      of $~T\,$'s proofs {\it are shorter}
      than  
      the errors in
      $~\theta\,$'s $\Pi_1^\xi$ styled-statements.
      The pleasing aspect 
       about this
      observation, illustrated by 
      % for example
      Remark \ref{re3-1},
      %
      % epistemological bundling
      is that those same invariants,
      $\, * \,$ and  $\, ** \,$, 
      which 
%%%%%      can 
      tempt a
       thinking agent $~T~$ 
      to engage in a trial-and-error walk through
      RE-Class$(\xi)$'s bundled universe, 
      %% are 
      also 
      %% the invariants that 
      make
      viable
      %% the operating prerequisites for making
       \thx{ppp1}'s self-justifying formalisms.

      %active.
      
      Thus aside from 
      addressing
       $\, \# \,$'s dilemma about the nature of 
      instinctive faith, 
      the 
      % mathematical
      meta-formalism in this appendix
      is
      %was 
      % also
      useful 
      %in 
      %venting
      % for
      in  explaining
      the 
      % underlying
       motivation behind the 
      %very
      %quite
      % fairly
      elaborate
      network 
      % labyrinth
      of
      theorems, proofs and definitions 
      %raised
      that were introduced 
      % had 
      %appeared 
      in this 
      paper.
In summary,
EA-stable logics are thus
 interesting   both
in their own right
(as a vehicle 
 enabling a Thinking Being to partially tolerate
its own errors), and 
because
they are  useful in explaining
how a Thinking Being 
can possess a type of instinctive
faith in its own consistency
(via the 
reflection
principles of
 Theorem \ref{ppp6}
and 
of  Remarks \ref{f88} and \ref{remhappy}).

\cvl

\section*{Appendix G:  Improvements upon Theorems  \ref{pqq4} 
and \ref{pqq5}  }

Let us recall that 
Remark \ref{re4-n} indicated that there was a
subtle
trade-off between 
Theorems  \ref{pqq4} 
and \ref{pqq5},
where neither result was 
strictly
 better than the other.
This section
will introduce
two hybrid 
methodologies, using 
Definition G.1's
formalism, that 
improve upon \thx{pqq5}
while  retaining a large part of \thx{pqq4}'s
nice features.

\hgskip

% \cvt
\parskip 2pt

{\bf Definition G.1}
Let $~\xi~$ denote 
the  generic configuration, 
whose base axiom system is again denoted as 
$ \, B^\xi \,$,
$\,~ \Phi~$ denote any          $\Pi_1^\xi$ 
sentence 
that   is true in    the Standard-M model
and 
$~j~$ denote an index that represents some
predicate 
 Test$^\xi_j \,$ 
lying 
in Definition \ref{gsim}'s
 TestList$^\xi$ sequence. Then
a
$\,\Pi_1^\xi \,$ sentences $\Psi $ 
 will be said to be a 
{\bf Braced}$^\xi( \, \Phi \, , \, j \, )$ expression when
$~ B^\xi \, + \, \Phi~$ can prove:
\beq
\label{punch}
\{~~~ \forall ~x~~~
\mbox{Test}_j^\xi(~\lceil~\Psi~\rceil~,~x~) ~~~\}
 ~~~~\longrightarrow ~~~~ \Psi 
\enq

\medskip

{\bf Theorem G.2}
{\it $~$Let  $ \,\xi \,$ again 
denote an arbitrary  generic configuration  
  \gggcp, 
and let
 $(  \nal{B},D)$ again denote any second axiom system and deduction
method whose $\Pi_1^\xi$ theorems are true under the
Standard-M model.
Then for any 
integer $~j~$ and for any 
  $\Pi_1^\xi$ sentence
$~\Phi~$ that is  true  in the
Standard-M model,
 the following  invariants
do hold:}
\bed
\item[   i ] 
{\it If $ \, \xi \, $ is EA-stable\sss  
then there 
will exist a self-justifying
$~\beta_j~\supset~B^\xi$ that 
can recognize its 
Level$(1^\xi$)
 consistency,
contains
only {\bf a finite number} of additional axioms 
beyond those appearing in 
$~B^\xi$,
 and which
can 
prove
all  of $(  \nal{B},D)$'s $\Pi_1^\xi$ theorems that 
are Braced$^\xi(  \Phi  ,j)$ expressions.}
\item[  ii ] 
{\it Likewise, 
if $\xi$ is 
  E-stable, A-stable  or 0-stable 
then 
 a
self-justifying
$\beta_j \supset B^\xi$ 
will exist
with the same properties except 
that it
recognizes its own
Level$(0^\xi$)
 consistency.}
\ennd

\medskip

{\bf Proof.}
To justify
Theorem G.2, we must first define
the axiom system  $~\beta_j~,~$ 
whose existence is claimed by
Items (i) and (ii).
It will be defined to
consist of the union of
the initial base axiom system
 $B^\xi$ 
with the following three added axiom-sentences.
\bed
\item[   1  ]
The $\Pi_1^\xi$ sentence $~\Phi~$
used
by  Definition G.1's
 Braced$^\xi(  \Phi  ,j  )$ formula.

\smallskip

\item[   2  ]
A GlobSim$^D_{  \nal{B}} \,(\xi,j)$ sentence whose indexing
integer $ \, j \,$ 
is defined by  Definition G.1.
This global simulation sentence is
thus the statement:
\beq
\label{glob2}
\forall ~t~~
\forall ~q~~
\forall ~x~~\{~~
[~~\mbox{Prf}_{  \nal{B}}^D \,(   t   ,   q   )~~ \wedge ~~
\mbox{Check}^\xi(t)~~]~~~
\longrightarrow ~~~ 
\mbox{Test}_j^\xi(t,x)  ~~~\}
\enq
\item[   3  ]
A  $\Pi_1^\xi$ sentence
 of the form
$\mbox{SelfCons}^k\{~[~\theta \,\cup \,B^\xi~]~,d~\}~$
where:
\bed
\item[   a   ] 
$\theta~$ is an R-view
consisting of the two
 $\Pi_1^\xi$ sentences defined by
Items 1 and 2.
\item[   b   ] 
$B^\xi\,$ is $~\xi \,$'s base axiom system,
and
\item[   c   ] 
$k~$ equals respectively 
 1  and 0 under formalisms  (i) and (ii).
\ennd
\ennd
Thus,
the  system $\beta_j$
uses identical definitions
 under
formalisms (i) and (ii),      
except 
that
its third sentence will use a different value for
$~k~$.
Our proof of 
Theorem G.2
will require        first confirming
the following fact:
\begin{description}
\item[  Claim *  ]
The axiom system
$\,\beta_j \,$ (which consists of the union of
$B^\xi$  with the
  sentences 1-3)
will have a capacity
to prove every 
Braced$^\xi(  \Phi  ,j)$ 
sentence $~\Psi~$
that is a
 $\Pi_1^\xi$
  theorem of
$(  \nal{B},D)$.
\end{description}
The proof of Claim * is quite simple.
It will rest on
the following
three
observations:
\bed
\item[   a  ]
For each  $\Pi_1^\xi$ sentence $\Psi$,
the  system $~\beta_j~$
must certainly
have a capacity to prove 
\eq{glob21}'s sentence (which 
states 
that $~\Psi\,$'s
G\"{o}del number 
formally
encodes
a $\Pi_1^\xi \,$  statement). This  
is because 
\eq{glob21}
 is true
 in the Standard-M
model
 and
because Part 3 of Definition \ref{def3.3}
indicated
that
 the 
$~B^\xi~$ 
sub-component of $~\beta_j~$
has a capacity to prove 
 every
$\Delta_0^\xi$ sentence
 that is true.
\beq
\label{glob21}
\mbox{Check}^\xi( ~ \lceil \, \Psi \, \rceil ~)
\enq
\item[   b  ]
Since
Claim $\,* \,$ specifies
 $~\Psi~$ is a theorem of 
 $(  \nal{B},D)$, 
there must certainly exist some integer
$~N~$ that is the G\"{o}del number of its proof from
 $(  \nal{B},D)$. This implies that 
\eq{glob22} must be a true
$\Delta_0^\xi$ sentence
 under the Standard-M model.
As was the case with \ep{glob21}, this implies
that
it must be provable from
$~B^\xi~$ (because it is a valid 
$\Delta_0^\xi$ sentence).
\beq
\label{glob22}
\mbox{Prf}_{  \nal{B}}^D \,(   ~ \lceil \, \Psi \, \rceil ~  , ~  N ~ )
\enq
\item[   c  ]
It is apparent that 
Equations \eq{glob2}, \eq{glob21} and \eq{glob22}
imply the validity of \eq{glob23}.
Moreover, Part 4 of Definition \ref{def3.3} indicated that
the generic configuration $~\xi\,$'s deduction method
does satisfy  G\"{o}del's Completeness Theorem.
This fact assures that 
$~\beta_j~$ must be able to prove 
\eq{glob23} because it contains
\eq{glob2}
as an axiom and 
 \eq{glob21} and \eq{glob22} 
as derived theorems \footnote{\label{fcomp} \sm55 
Every deduction method $\,d$,
$\,$satisfying G\"{o}del's Completeness Theorem,
will be automatically able to prove
a theorem $~Z~$
when it contains
$X$, $Y$ and $~(X \wedge Y)~\rightarrow ~Z~$ as
theorems, irregardless of whether or not it contains
an explicit built-in
modus
ponens rule.
Thus $~d~$ can prove
\eq{glob23} because of its knowledge about 
\eq{glob2}--\eq{glob22}'s validity.}.
\ennd
\vspace*{- 0.1 em}
\beq
\label{glob23}
\forall ~x~~~
\mbox{Test}_j^\xi(   ~ \lceil \, \Psi \, \rceil ~  , ~  x ~ )
\vspace*{- 0.1 em}
\enq
Claim $*$ is a
consequence of 
Observations a-c. This is  because
$ \, \Phi \, $ 
is one of $ \, \beta_j \,$'s defined axioms,
and
Definition G.1 
indicated 
 $ \, B^\xi \, + \, \Phi \, $ was capable of proving
\eq{punch}'s statement
 for every  Braced$^\xi(  \Phi  ,j  )$ sentence $ \, \Psi \, $.
These facts corroborate Claim $*$ 
because they imply
that $ \, \beta_j \, $
must be able to verify
Claim $\,* \,$'s
 sentence
 $ \, \Psi \, $ (because
 $ \, \beta_j \, $
can verify statements
\eq{punch} and \eq{glob23}).

\medskip

The remainder of 
Theorem G.2's proof is 
analogous
to
\phx{pqq5}'s proof.
This is because the prior paragraph 
established
that $~\beta_j~$  can prove 
every 
 Braced$^\xi(  \Phi  ,j  )$ 
theorem of
 $(  \nal{B},D)$
(as was required by
 Claims i and ii ). 
The only 
remaining task is to
show that
 $~\beta_j~$  is a self-justifying formalism that can
recognize its 
Level($1^\xi$) and Level($0^\xi$)
consistencies,
as specified by
 Claims i and ii.  
This part of 
Theorem G.2's
 verification
is identical
to
the methods used
to prove
Theorems \ref{pqq3} and \ref{pqq5}.
It
will
thus
 not be repeated
here.
  $~~\Box$

\medskip

The last part of this appendix will
require the
following additional
 notation to formalize
the main intended application
of
Theorem G.2's 
formalism.
\bee
\item
\topsep -7pt
 Count$(  \Psi )$ will denote
the number of 
quantifiers appearing in the  sentence $\Psi$
(including both its bounded and unbounded quantifiers).
\item 
Size$^\xi(c)$ 
will
denote the  set of $\Pi_1^\xi$ sentences
$\Psi$ where  Count$(  \Psi ) \, \leq \, c \,$.
\ene
Our next theorem will be   a specialized
variant of
Theorem G.2, using the 
Size$^\xi(c)$  construct.
It  will explain the 
 intended application
of
this
formalism:

\medskip

{\bf Theorem G.3.}
{\it 
$~$Let $~\xi~$ denote any one of Appendix D's four sample
generic configurations of $~\xi^*~$, 
$~\xi^{**}~$,  $~\xi^-~$  or   $~\xi^R~$.
Then 
 for any  $~c>0~$, 
Theorem G.2's axiom systems
of $~\beta_j~$
can be arranged so that 
they can prove all of 
  $  (    \nal{B} ,  D  ) $'s
Size$^\xi(c)$ 
 $\Pi_1^\xi$ 
theorems while
simultaneously also
 recognizing their:
\bee
\item Level(1)
consistency  for the cases
when  $ \, \xi \, $ is one of   $ \, \xi^* \, $,
$ \, \xi^{**}\,   $ or   $ \, \xi^-$.
\item
 Level(0) consistency when $~\xi~$ is
$~\xi^R~$.
\ene}

\cvmew

{\bf Proof Sketch:}
%% 
%%  There is 
%% insufficient
%%  space to prove
%% Theorem G.3 here, but its intuition 
%% is easy to summarize.
%% 
%% 
 The  intuition behind
Theorem G.3's proof is 
quite easy to summarize.
For 
arbitrary $     c>0      $
and any
of Appendix D's configurations of
 $   ~  \xi^*  ~   $, 
$   ~  \xi^{**}  ~   $,  $  ~   \xi^-  ~   $  and   $  ~   \xi^R  ~   $, 
it is 
routine to
 construct an ordered pair $  ~   (\Phi,j) ~    $ where every
$\Pi_1^\xi$ sentence of
Size$^\xi(c)$ is a 
 Braced$^\xi(  \Phi  ,j  )$ expression.
Theorem G.3's 
first claim is,
thus, a 
consequence of 
Part (i)
of Theorem G.2
 and 
the fact that each of
 $   ~  \xi^*  ~   $, 
$   ~  \xi^{**}  ~   $ and  $  ~   \xi^-  ~   $  
are EA-stable.
Likewise, Theorem G.3's 
second claim follows from
Part (ii)
of Theorem G.2
 and the fact that 
$ ~  \xi^R  ~ $ is E-stable,
$~~\Box$  

\bigskip
\medskip

{\bf Remark G.4.}
The  Theorems G.2 and G.3 are
of interest  because the set of 
$\Pi_1^\xi$ sentences of
Size$^\xi(c)$ is a natural class to examine.
It is,
thus, tempting to consider
a system that 
recognizes
its own
formal 
 consistency, uses only a finite
number of axiom sentences beyond those in $~B^\xi~,~$ and 
which
can
 prove all of
  $  (    \nal{B} ,  D  ) $'s
 $\Pi_1^\xi$ theorems 
of Size$^\xi(c)$.
Such a system 
replies to 
Remark 
\ref{re4-n}'s challenge by
% nicely
 hybridizing the properties of
Theorems  \ref{pqq4} 
and \ref{pqq5},
%%  .
in a seemingly pragmatic manner.

\end{document}